\documentclass{amsart}
\usepackage{amsmath,amssymb,stmaryrd}
\usepackage{amsrefs}
\usepackage{color}
\usepackage{enumerate}

\def \rr {\mathbb{R}}
\def \nn {\mathbb{N}}
\def \rn {\rr^n}

\def \rnm {\rr^n_-}
\def \rnp {\rn\setminus\{0\}}

\def \huno {H_0^1(\Omega)}
\def \dundeuxrn {D^{1,2}(\rn)}

\def \eps {\epsilon}
\def \ue {u_\epsilon}

\def \Ue {U_\epsilon}

\def \ua {u_\alpha}
\def \tua {\tilde{u}_\alpha}
\def \aa {a_\alpha}
\def \la {\lambda_\alpha}
\def \xa {x_\alpha}
\def \ma {\mu_\alpha}
\def \va {v_\alpha}
\def \ya {y_\alpha}
\def \za {z_\alpha}
\def \ra {r_\alpha}
\def \da {d_\alpha}

\newtheorem{proposition}{Proposition}

\newtheorem{lem}{Lemma}

\def \bp {\beta_+(\gamma)}
\def \bm {\beta_-(\gamma)}
\def \bps {\beta_+}
\def \bms {\beta_-}

\def \Omegabar {\overline{\Omega}}

\def \crit {2^{\star}}
\def \crits {2^\star(s)}

\newtheorem{theorem}{Theorem}

\title[Hardy-Schr\"odinger operator with interior singularity]{The Hardy-Schr\"odinger operator with interior singularity: The remaining cases  
}
\author{Nassif Ghoussoub}
\address{Nassif Ghoussoub, Department of Mathematics, 1984 Mathematics Road, The University of British Columbia, BC, Canada V6T 1Z2}
\email{nassif@math.ubc.ca}

\author{Fr\'ed\'eric Robert}
\address{Fr\'ed\'eric Robert, Institut \'Elie Cartan, Universit\'e de Lorraine, BP 70239, F-54506 Vand{\oe}uvre-l\`es-Nancy, France}
\email{frederic.robert@univ-lorraine.fr}
\date{August 24, 2017}

\thanks{This work was carried out while N. Ghoussoub was visiting l'Institut \'Elie Cartan, Universit\'e de Lorraine in May 2015. He  was partially supported by a research grant from the Natural Science and Engineering Research Council of Canada (NSERC). The paper was completed  while Fr\'ed\'eric Robert was visiting the University of British Columbia in July 2016. }

\begin{document}

\begin{abstract}
We consider the remaining unsettled cases in the problem of existence of energy minimizing solutions for the Dirichlet value problem $L_\gamma u-\lambda u=\frac{u^{2^*(s)-1}}{|x|^s}$ on a smooth bounded domain $\Omega$  in $\rn$ ($n\geq 3$) having the singularity $0$ in its interior.  
Here $\gamma <\frac{(n-2)^2}{4}$, $0\leq s <2$, $2^*(s):=\frac{2(n-s)}{n-2}$ and $0\leq \lambda <\lambda_1(L_\gamma)$, the latter being the first eigenvalue of the Hardy-Schr\"odinger operator $L_\gamma:=-\Delta -\frac{\gamma}{|x|^2}$. 
There is a threshold $\lambda^*(\gamma, \Omega) \geq 0$ beyond which the minimal energy is achieved, but below which, it is not.  It is well known that $\lambda^*(\Omega) = 0$ in higher dimensions, for example if $0\leq \gamma \leq \frac{(n-2)^2}{4}-1$. Our main objective in this paper is to show that this threshold is strictly positive in ``lower dimensions" such as when $ \frac{(n-2)^2}{4}-1<\gamma <\frac{(n-2)^2}{4}$, to identify the critical dimensions (i.e., when the situation changes), and to characterize it in terms of $\Omega$ and $\gamma$. 
If either $s>0$ or if $\gamma > 0$, i.e., in {\it the truly singular case}, we show that  in low dimensions, a solution is guaranteed by the positivity of the ``Hardy-singular internal mass" of $\Omega$, a notion that we introduce herein. On the other hand, and just like the case wnen $\gamma=s=0$ studied by Brezis-Nirenberg \cite{bn} and completed by Druet \cite{d2}, $n=3$ is the critical dimension, and the classical positive mass theorem is sufficient for the {\it merely singular case}, that is when $s=0$, $\gamma \leq 0$.
\end{abstract}

\maketitle

\tableofcontents
\section{\, Introduction}

Let $\Omega$ be a smooth bounded domain in $\rn$ ($n\geq 3$) such that $0\in \Omega$  and consider the following Dirichlet boundary value problem:
 \begin{eqnarray} \label{one}
\left\{ \begin{array}{llll}
-\Delta u-\gamma \frac{u}{|x|^2}-\lambda u&=&\frac{u^{2^*(s)-1}}{|x|^s}  \ \ &\text{on } \Omega,\\
\hfill u&>&0 &\text{on } \Omega, \\
\hfill u&=&0 &\text{on }\partial \Omega, 
\end{array} \right.
\end{eqnarray}
where $\gamma <\frac{(n-2)^2}{4}$, $0\leq s <2$, $\crits:=\frac{2(n-s)}{n-2}$ and $0\leq \lambda <\lambda_1(L_\gamma)$, the latter being the first eigenvalue of the Hardy-Schr\"odinger operator $L_\gamma:=-\Delta -\frac{\gamma}{|x|^2}$, that is 
\begin{equation*}
\lambda_1(L_\gamma, \Omega):=\inf\left\{\frac{\int_{\Omega} |\nabla u|^2\, dx-\gamma \int_{\Omega}\frac{u^2}{|x|^2}dx}{\int_{\Omega} u^2\, dx};\, u\in \huno\setminus \{0\}\right\}.
\end{equation*}
Equation (\ref{one}) is essentially the Euler-Lagrange equation corresponding to the following energy functional on $\huno$,
\[
J_{\gamma, s, \lambda}^\Omega(u)=\frac{\int_{\Omega} |\nabla u|^2\, dx-\gamma \int_{\Omega}\frac{u^2}{|x|^2}dx-\lambda \int_{\Omega} |u|^2\, dx }{\left(\int_{\Omega}\frac{u^{2^*(s)}}{|x|^s}dx\right)^{\frac{2}{2^*(s)}}},
\]
where $\huno$ is the completion of $C^\infty_c(\Omega)$ for the norm $u\mapsto \Vert \nabla u\Vert_2$. We shall therefore study whether the following minimization problem 
\begin{equation*} 
\mu_{\gamma, s, \lambda}(\Omega):=\inf\left\{J_{\gamma, s, \lambda}^\Omega(u);\, u\in \huno\setminus\{0\}\right\},
\end{equation*}
is attained, that is if $\mu_{\gamma, s, \lambda}(\Omega)=J_{\gamma, s, \lambda}^\Omega(u_0)$ for some $u_0\in \huno$. For convenience, we define $\dundeuxrn:=H_0^1(\rn)$, that is the completion of $C^\infty_c(\rn)$ for the norm $u\mapsto \Vert\nabla u\Vert_2$.

Note that the fact that $\mu_{\gamma, s, 0}(\rn) >0$ is equivalent to the critical case of the Caffarelli-Kohn-Nirenberg inequalities \cite{ckn}. In particular, see for instance Ghoussoub-Robert \cite{gr5},
$$\mu_{\gamma,s,0}(\rn)\hbox{ is achieved iff }\{s>0\}\hbox{ or }\{s=0\hbox{ and } \gamma\geq 0\}.$$
It is also standard that $\mu_{\gamma, s, 0}(\Omega)= \mu_{\gamma, s, 0}(\rn)$ whenever  $\Omega$ is a domain containing $0$ in its interior, and hence  $\mu_{\gamma, s, 0}$ is not attained if $\Omega$ is bounded. 

 The idea of restoring compactness  by considering non-trivial negative linear perturbations was pioneered by Brezis-Nirenberg \cite{bn} in the case when $\gamma=0$, $s=0$ and $0<\lambda <\lambda_1(\Omega)$, the latter being the first eigenvalue of the Laplacian on $\huno$. They showed that in this case (\ref{one}) 
has a solution for $n\geq 4$. 
The case $n=3$ is special and involves a ``positive mass" condition introduced by Druet \cites{d2, d3}, and inspired by the work of  Schoen \cite{schoen1} on the Yamabe problem.  
The bottom line is that --at least for $\gamma=0$-- the geometry of $\Omega$ need not be taken into account in dimension $n\geq 4$, while in dimension $n=3$, the existence depends on the domain $\Omega$ via ``a positive mass condition". We shall elaborate further on this theme. 

In this paper, we consider the case when the Laplacian is replaced by the Hardy-Schr\"odinger  operator $L_\gamma$. Here, the position of the singularity $0$ within $\bar \Omega$ matters a great deal. In \cite{gr4}, we considered the case where $0$ belongs to the boundary $\partial \Omega$ of the domain $\Omega$. In this sequel, we deal with the case when $0\in \Omega$, which was first considered by Janelli \cite{Jan} in the case $s=0$. 

It is already well known that there is a threshold $\lambda^*$ beyond which the infimum $\mu_{\gamma,s,\lambda}(\Omega)$ is achieved, and below which, it is not. It can be characterized as
\begin{equation}\label{def:star}
\lambda^*(\Omega):
=\sup\{\lambda;\, \mu_{\gamma,s,\lambda}(\Omega)= \mu_{\gamma,s,0}(\rn)\}.
\end{equation}
It is easy to see that $0\leq \lambda^*(\Omega) < \lambda_1(L_\gamma,\Omega)$. It is also part of the folklore --that we sketch below-- that $\lambda^*(\Omega) = 0$ in higher dimensions. Our main objective in this paper is to show that this threshold is strictly positive in ``lower dimensions," to identify the critical dimensions (i.e., when the situation changes), and to try to characterize it in terms of $\Omega$ and $\gamma$.

As opposed to Brezis-Nirenberg \cite{bn} and Druet \cite{d2}, we are dealing here with the case where $0$ is an interior singularity, and our analysis below shows that the identification of $\lambda^*$ differ according to two distinct singularity regimes: 
\begin{itemize}
\item The {\it truly singular case}, which corresponds to when either $s>0$ or $\gamma>0$. We note that in this case $\mu_{\gamma,s,0}(\rn)$ is achieved.
\item The {\it merely singular case}, which corresponds to the case when $s=0$ and $\gamma \leq 0$, a case where $\mu_{\gamma,s,0}(\rn)$ is not achieved, unless $s=\gamma=0$.
\end{itemize}
The following three theorems are the main results of this paper. The first is rather standard. It deals with high dimensions and is included for completeness and comparison purposes. The second deals with the low dimensional cases, i.e., the remaining cases which are yet to be addressed in the literature. 

\begin{theorem} \label{main.interior} {\rm \bf (The higher dimensional case)} Let $\Omega$ be a smooth bounded domain in $\rn$ ($n\geq 3$) such that $0\in  \Omega$. Assume that we are in the following situation:
\begin{itemize}
\item either in the truly singular case and $\gamma <\frac{(n-2)^2}{4}-1$,
\item or in the merely singular case and $n\geq 4$.
\end{itemize}
Then $\mu_{\gamma,s,\lambda}(\Omega)$ is achieved if and only if $\lambda>\lambda^*(\Omega)$. Moreover, 
\begin{enumerate}[i)]
\item  In the truly singular case (i.e, when either $s>0$ or $\gamma > 0$), and    if  $\gamma\leq\frac{(n-2)^2}{4}-1$, then 
\begin{equation}\label{*}
\lambda^*(\Omega)=0.
\end{equation}
\item In the merely singular case (i.e, when $s=0$ and $\gamma \leq 0$), and if $n\geq 4$,  then 
\begin{equation}\label{**}
\lambda^*(\Omega)=\inf \left\{\frac{|\gamma|}{|x|^2}; \, x\in \Omega\right\}>0\hbox{ if }\gamma<0.
\end{equation}
\end{enumerate}  
\end{theorem}
Part (i) of Theorem \ref{main.interior} was proved by Janelli \cite{Jan} in the case when $s=0$. The case when $s>0$ is not much different and was noted in several works such as  \cites{RW, KP1, KP2, KP3, CHa, CHP, CP}.  Part (ii) of Theorem \ref{main.interior}, that is the case when $s=0$ and $\gamma <0$, in dimension $n\geq 4$ was also tackled by Janelli \cite{Jan} and  Ruiz-Willem \cite{RW}. Their proof, though not complete, essentially gives the above result.

Janelli \cite{Jan} also considered the lower dimensional case, that is
 \begin{equation*}
 \frac{(n-2)^2}{4}-1<\gamma <\frac{(n-2)^2}{4},
 \end{equation*}
 when $\Omega$ is the  ball $B$ centered at $0$. He gave the following explicit value for $\lambda^*$: 
\begin{equation} \label{lambda.star.0}
 \lambda^*(B)=\inf\left\{\frac{\int_{B}\frac{|\nabla u|^2}{|x|^{2\beta_+}}dx}{\int_{B}\frac{u^2}{|x|^{2\beta_+}}dx};\, u\in H^1_0(B)\setminus \{0\}\right\}>0, 
\end{equation}
 where 
\begin{equation*}
\beta_{\pm}(\gamma):=\frac{n-2}{2}\pm\sqrt{\frac{(n-2)^2}{4}-\gamma}.
\end{equation*}
Note that the radial function $x\mapsto |x|^{-\beta}$ is a solution of  $(-\Delta -\frac{\gamma}{|x|^2})u=0$ on $\rn\setminus\{0\}$ if and only if $\beta\in \{\beta_-(\gamma),\beta_+(\gamma)\}$. 

In order to characterize the threshold $\lambda^*(\Omega)$ for a general domain $\Omega$, we need to define the notion of  {\it Hardy-singular interior mass} associated to the operator $-\Delta -\frac{\gamma}{|x|^2} -\lambda$ on a bounded domain $\Omega$ in $\rn$ containing  $0$.

\begin{theorem} \label{interior.mass} {\rm \bf (The Hardy singular internal mass)} Let $\Omega$ be a smooth bounded domain in $\rn$  ($n\geq 3$) such that $0\in \Omega$. Suppose $h$ is a $C^2$-potential on $\Omega$ so that the operator $-\Delta-(\frac{\gamma}{|x|^2}+h(x))$ is coercive. 
\begin{enumerate}[i)] 
\item There exists then  $H\in C^\infty(\overline{\Omega}\setminus\{0\})$ such that 
$$(E)\qquad\qquad \left\{\begin{array}{ll}
-\Delta H-\left(\frac{\gamma}{|x|^2}+h(x)\right)H=0 &\hbox{ in } \Omega\setminus\{0\}\\
\hfill H>0&\hbox{ in } \Omega\setminus\{0\}\\
\hfill H=0&\hbox{ on }\partial\Omega.
\end{array}\right.$$
These solutions are unique up to a positive multiplicative constant, and there exists $c>0$ such that 
$H(x)\simeq_{x\to 0}\frac{c}{|x|^{\beta_+(\gamma)}}.$
\item If either $h$ is sufficiently small around $0$ or if $\frac{(n-2)^2}{4}-1<\gamma<\frac{(n-2)^2}{4}$, then for any solution $H\in C^\infty(\overline{\Omega}\setminus\{0\})$ of $(E)$, there exist $c_1>0$ and $c_2\in\rr$ such that 
$$H(x)=\frac{c_1}{|x|^{\bp}}+\frac{c_2}{|x|^{\bm}}+o\left(\frac{1}{|x|^{\bm}}\right) \qquad \hbox{ as } x\to 0.$$
The uniqueness implies that the ratio $c_2/c_1$ is independent of the choice of $H$, hence the  ``Hardy-singular internal mass" of $\Omega$ associated to the operator $L_\gamma - h(x)I$ can be defined unambigously as 
$$m_{ \gamma, h}(\Omega):=\frac{c_2}{c_1}\in\rr.$$
\end{enumerate}
\end{theorem} 
For the merely singular case ($s=0$ and $\gamma \leq 0$) and the critical dimension $n=3$, we need a more standard 
notion of mass associated to the operator $L_\gamma$ at an internal point $x_0\in \Omega$, which is reminiscent of Schoen-Yau's approach to complete the solution of the Yamabe conjecture in low dimensions. For that, one considers for a given $\gamma\leq 0$, the corresponding {\it Robin function} or the regular part of the Green function with pole at $x_0\in \Omega\setminus \{0\}$. One shows that for $n=3$, any solution $G$ of 
$$\left\{\begin{array}{ll}
-\Delta G-\frac{\gamma}{|x|^2} G-\lambda G=0 &\hbox{ in } \Omega \setminus\{x_0\}\\
\hfill G>0 &\hbox{ in }\Omega \setminus\{x_0\}\\
\hfill G=0 &\hbox{ on }\partial\Omega,
\end{array}\right.$$
 is unique up to multiplication by a constant, and that there exists $R_{\gamma, \lambda}(\Omega, x_0)\in \rr$ and $c_{\gamma, \lambda}(x_0)>0$ such that
\begin{equation}\label{def:R}
G(x)=c_{\gamma, \lambda} (x_0)\left(\frac{1}{|x-x_0|^{n-2}}+R_{\gamma, \lambda}(\Omega, x_0)\right)+o(1)\quad \hbox{ as }x\to x_0.
\end{equation}
The quantity $R_{\gamma, \lambda}(\Omega, x_0)$ is then well defined and will be called the {\it internal mass} of $\Omega$ at $x_0$. We then define
\[
\hbox{$R_{\gamma, \lambda}(\Omega)=\sup\limits_{x\in \Omega\setminus\{0\}}R_{\gamma, \lambda}(\Omega, x)$.}
\]
These will allow us to give an explicit value for $\lambda^*(\gamma, \Omega)$ as follows.
 
\begin{theorem} \label{gamma.star} {\rm \bf (The lower dimensional case)} Let $\Omega$ be a smooth bounded domain in $\rn$ ($n\geq 3$) such that $0\in  \Omega$. 
\begin{enumerate}[i)]
\item Assume we are
\begin{itemize}
\item either in the truly singular case and $\frac{(n-2)^2}{4}-1  <\gamma <\frac{(n-2)^2}{4}$,
\item or in the merely singular case and $n=3$.
\end{itemize}
Then, there exists $\lambda^*(\Omega)>0$ such that 
$\mu_{\gamma, s, \lambda}(\Omega)$ is not achieved for $\lambda<\lambda^*$
and $\mu_{\gamma, s, \lambda}(\Omega)$ is  achieved for $\lambda>\lambda^*.$
\item Moreover, in the truly singular case, 
with $\frac{(n-2)^2}{4}-1  <\gamma <\frac{(n-2)^2}{4}$, and under the assumption that $\mu_{\gamma, s, \lambda^*}(\Omega)$ is not achieved, we have that $m_{\gamma, \lambda^*}(\Omega)= 0$, and 
 \begin{equation}\label{lambda.0}
 \lambda^*(\Omega)=\sup \{\lambda;\,  m_{\gamma, \lambda}(\Omega)\leq 0\}.
 \end{equation}  
\item In the merely singular case, and with $n=3$,  then $\mu_{\gamma,s,\lambda^\star}(\Omega)$ is not achieved and
\begin{equation}\label{***}
\lambda^*(\Omega)=\sup \{\lambda;\,  R_{\gamma, \lambda}(\Omega)\leq 0\}.
\end{equation}
\end{enumerate}
\end{theorem}
We conjecture that in all cases, $\mu_{\gamma, s, \lambda^*}(\Omega)$ is never achieved, which means that (\ref{lambda.0}) must hold unconditionally. Note that $\mu_{\gamma, s, \lambda^*}(\Omega)=\mu_{\gamma, s, 0}(\rn)$, but we don't know whether this suffices to conclude that $\mu_{\gamma, s, \lambda^*}(\Omega)$ is not achieved.
 When $s=\gamma=0$ and $n=3$, Druet \cite{d2} proved  that this is indeed the case by using a very elegant geometric argument. This extends to the merely singular case. In the truly singular case, the conjecture holds in the radially symmetric case, i.e., when $\Omega$ is a ball. This was verified by Janelli \cite{Jan}.

 Finally, we note that the above analysis lead to the following definition of a {\it critical dimension for the operator $L_\gamma$}. It is the largest scalar $n_\gamma$ such that for $n < n_\gamma$, there exists a bounded smooth domain $\Omega \subset \rn$ and a $\lambda \in (0, \lambda_1(L_\gamma, \Omega))$ such that there is a non-trivial minimiser satisfying \eqref{one}. \\
$\mu_{\gamma, s, \lambda}(\Omega)$ is not attained.
\begin{equation*}
n_\gamma =\left\{\begin{array}{ll} 2\sqrt{\gamma +1}+2\quad & \, {\rm if}\, \gamma \geq -1\\
2 \quad  & \, {\rm if}\,  \gamma <-1.
\end{array}\right.
\end{equation*}
Note that $n < n_\gamma$ is exactly when $\bp-\bm < 2$, which is the threshold where
the radial function $x\to |x|^{-\beta_+(\gamma)}$ is locally $L^2$-summable.

\medskip\noindent{The proofs of Theorems \ref{main.interior} and \ref{gamma.star} rely on a refined blow-up analysis for certain families of solutions of equation \eqref{one}. We give --in Theorems \ref{th:blowup:1} and \ref{th:blowup:2} below-- a complete description of how such blowups may occur. In particular, we show that in the truly singular case, the solutions necessarily concentrate at the singularity $0$, while in the merely singular case, they do so at a point $x_0\neq 0$ of the domain $\Omega$. In the appendices, we establish several important properties of the Green function associated to the operator $-\Delta-\gamma|x|^{-2}$, that are crucial for the proofs of Theorems \ref{th:blowup:1} and \ref{th:blowup:2}.

\section{\, The higher dimensional case}\label{sec:high}

We recall the following facts, which by now are standard. 

\begin{enumerate}[i)]

\item $\mu_{\gamma,s,\lambda}(\Omega)>0 $ if and only  $\lambda< \lambda_1(L_\gamma,\Omega)$.
\item $\mu_{\gamma,s,\lambda}(\Omega)= \mu_{\gamma,s,0}(\rn)$ for all $\lambda\leq 0$.
\item $\mu_{\gamma,s,\lambda}(\Omega)$ is attained if $\mu_{\gamma,s,\lambda}(\Omega)< \mu_{\gamma,s,0}(\rn)$.
\item The function $\lambda\mapsto \mu_{\gamma,s,\lambda}(\Omega)$ is continuous and nonincreasing.
\item If $\mu_{\gamma,s,\lambda}(\Omega)$ is attained, 
then $\mu_{\gamma,s,\lambda'}(\Omega)<\mu_{\gamma,s,\lambda}(\Omega)$ for any $\lambda'>\lambda$.
\end{enumerate}
Writing $\lambda^*=\lambda^*(\Omega)$ for short, where $\lambda^*(\Omega)$ is defined in \eqref{def:star}, it follows from the above that
\begin{enumerate}[vi)]
\item  $\mu_{\gamma,s,\lambda}(\Omega)= \mu_{\gamma,s,0}(\rn)$ for all $\lambda\leq \lambda^*$ and $\mu_{\gamma,s,\lambda}(\Omega)< \mu_{\gamma,s,0}(\rn)$ for all $\lambda> \lambda^*$.
\item $\mu_{\gamma,s,\lambda}(\Omega)$ is not achieved for all $\lambda<\lambda^*$.
\item $\mu_{\gamma,s,\lambda}(\Omega)$ is  achieved for all $\lambda>\lambda^*$.
\end{enumerate}

It is clear that $\lambda^* \geq 0$. This section is devoted to show that  $\lambda^*=0$ in ``high dimensions," which in our case will depend on $\gamma$.  The calculations are standard, and we include them for the convenience on the reader and for comparison to the other cases. As mentioned above in (iii),  in order to show that extremals exist for $\mu_{\gamma, s, \lambda}(\Omega)$, 
it suffices to prove that 
$\mu_{\gamma, s, \lambda}(\Omega) <\mu_{\gamma, s,0}(\rn),$
where $\mu_{\gamma, s,0}(\rn):=\mu_{\gamma, s, 0}(\rn)$. This kind of condition
is now standard when dealing with borderline variational problems. See also Aubin \cite{aubin}, Br\'ezis-Nirenberg \cite{bn}. The condition limits the energy level of minimizing sequences, prevents the creation of ``bubbles" and hence insures compactness. 

To show the strict inequality, one needs to test the functional $J_{\gamma, s, \lambda}^\Omega$
on minimizing sequences of the form $\eta U_\epsilon$, where $U_\epsilon$ is an extremal for $\mu_{\gamma, s,0}(\rn)$ 
and $\eta\in C^\infty_c(\Omega)$ is a cut-off function equal to $1$ in a neigbourhood of $0$.
It is therefore important to know for which parameters $\gamma$ and $s$, the best constant $\mu_{\gamma, s,0}(\rn)$ is attained. A proof of the following can be found in \cite{gr5}. For explicit extremals, we refer to Beckner \cite{beck} or Dolbeault et al. \cite{DELT}. 
 
\begin{proposition} Assume $\gamma <\frac{(n-2)^2}{4}$, $n\geq 3$ and $0\leq s <2$. Then,
\begin{enumerate}[i)]
\item  $\mu_{\gamma, s,0}(\rn)$ is attained if either $s>0$ or if $\{s=0$ and $\gamma \geq 0\}$.
\item If $0\leq \gamma < \frac{(n-2)^2}{4}$, then the extremals for  $\mu_{\gamma, s,0}(\rn)$ are explicit and take the form  $u_\eps (x)=c\cdot\eps^{-\frac{(n-2)}{2}}U (\frac{x}{\eps})$, where  $c\neq 0$, $\eps>0$ and
\begin{equation}\label{J}
U(x):=\frac{1}
{\left(|x|^{\frac{(2-s)\beta_-(\gamma)}{n-2}}+|x|^{\frac{(2-s)\beta_+(\gamma)}{n-2}}\right)^{\frac{n-2}{2-s}}}
\qquad \hbox{ for }x\in\rn\setminus\{0\},
\end{equation}

\item  On the other hand, if $s=0$ and $\gamma <0$, then $\mu_{\gamma, 0}(\rn)$ is not attained and is equal to $\mu_{0, 0}(\rn)$, which is the best constant in the Sobolev inequality.  
\end{enumerate}
\end{proposition}
\vskip 5pt
\noindent {\bf Subsection 2.1: The truly singular case}\\

We now give a proof of Theorem \ref{main.interior}. Assuming $\gamma\leq \frac{(n-2)^2}{4}-1$, we construct a minimizing sequence $u_\epsilon$ in $ \huno\setminus\{0\}$ for the functional 
$J_{\gamma,s, \lambda}$
in such a way that $\mu_{\gamma, s, \lambda}(\Omega)<\mu_{\gamma, s,0}(\rn)$.

\medskip\noindent Since either $s>0$ or $\gamma\geq 0$, then the infimum $\mu_{\gamma,s,0}(\rn)$ is achieved by 
the function 
$$U(x):=\frac{1}{\left(|x|^{\frac{(2-s)\bm}{n-2}}+|x|^{\frac{(2-s)\bp}{n-2}}\right)^{\frac{n-2}{2-s}}}
\hbox{ for }x\in\rn\setminus\{0\}.$$
In other words, $U\in \dundeuxrn$ and 
$J_{\gamma,s,0}^{\rn}(U)=\inf_{u\in \dundeuxrn\setminus\{0\}}J_{\gamma,s,0}^{\rn}(u),$
where
$$J_{\gamma,s,0}^{\rn}(u):=\frac{\int_{\rn}\left(|\nabla u|^2-\frac{\gamma}{|x|^2}u^2\right)\, dx}{\left(\int_{\rn}\frac{|u|^{\crits}}{|x|^s}\, dx\right)^{\frac{2}{\crits}}}\quad \hbox{for $u\in \dundeuxrn\setminus\{0\}$.}$$ 
In particular, there exists $\chi>0$ such that
\begin{equation}\label{eq:U}
-\Delta U-\frac{\gamma}{|x|^2}U=\chi \frac{U^{\crits-1}}{|x|^s}\hbox{ in }\rn\setminus\{0\}.
\end{equation}
For convenience, we will write in the sequel, 
$\bps:=\bp\hbox{ and }\bms:=\bm.$
Note that the assumption that $\gamma\leq \frac{(n-2)^2}{4}-1$ is equivalent to
$\bp-\bm\geq 2.$
 Define a scaled version of $U$ by 
\begin{equation*}
\Ue(x):=\eps^{-\frac{n-2}{2}}U\left(\frac{x}{\eps}\right)=\left(\frac{\eps^{\frac{2-s}{n-2}\cdot\frac{\bps-\bms}{2}}}{\eps^{\frac{2-s}{n-2}\cdot(\bps-\bms)}|x|^{\frac{(2-s)\bms}{n-2}}+ |x|^{\frac{(2-s)\bps}{n-2}}}\right)^{\frac{n-2}{2-s}}\quad \hbox{for $x\in\rn\setminus\{0\}$.}
\end{equation*}
\medskip\noindent Fix now a function $h\in C^{0,\theta}(\overline{\Omega})$, $\theta\in (0,1)$,  consider a cut-off function $\eta\in C^\infty_c(\Omega)$ such that $\eta(x)=1$ for $x$ in a neighborhood of $0$ contained in $\Omega$, and define for $\eps>0$ the test-function $\ue\in\huno$ by
$$\ue(x):=\eta(x)\Ue(x) \quad \hbox{for $x\in\Omegabar\setminus\{0\}$.}$$
We now estimate $J_{\gamma,s,h}^{\Omega}(\ue)$, where
$$J_{\gamma,s,h}^{\Omega}(u):=\frac{\int_{\Omega}\left(|\nabla u|^2-\left(\frac{\gamma}{|x|^2}+h(x)\right)u^2\right)\, dx}{\left(\int_{\Omega}\frac{|u|^{\crits}}{|x|^s}\, dx\right)^{\frac{2}{\crits}}}.$$

Note first that $0\leq \ue(x)\leq C\eps^{\frac{\bps-\bms}{2}}|x|^{-\bps}$ for all $x\in \Omega\setminus\{0\}$. Therefore, since $\bp-\bm\geq 2$, we have as $\eps\to 0$, 
\begin{eqnarray*}
\int_\Omega \frac{\ue^{\crits}}{|x|^s}\, dx&=& \int_{B_\delta(0)} \frac{\Ue^{\crits}}{|x|^s}\, dx+ o\left(\eps^2\right)\\
&=& \int_{B_{\eps^{-1}\delta}(0)} \frac{U^{\crits}}{|x|^s}\, dx+ o\left(\eps^2\right)\\
&=& \int_{\rn} \frac{U^{\crits}}{|x|^s}\, dx+ o\left(\eps^2\right).
\end{eqnarray*}
We also have
\begin{eqnarray*}
\int_{\Omega}\left(|\nabla \ue|^2-\frac{\gamma}{|x|^2}\ue^2\right)\, dx &=&  \int_{B_\delta(0)}\left(|\nabla \Ue|^2-\frac{\gamma}{|x|^2}\Ue^2\right)\, dx+ O\left(\eps^{\bps-\bms}\right)\\
&=& \int_{B_{\eps^{-1}\delta}(0)}\left(|\nabla U|^2-\frac{\gamma}{|x|^2}U^2\right)\, dx+ O\left(\eps^{\bps-\bms}\right)\\
&=&\chi\int_{B_{\eps^{-1}\delta}(0)}\frac{U^{\crits}}{|x|^s}\, dx+ O\left(\eps^{\bps-\bms}\right)\\
&=&\chi\int_{\rn}\frac{U^{\crits}}{|x|^s}\, dx+ O\left(\eps^{\bps-\bms}\right).
\end{eqnarray*}
Finally, we estimate the last term as $\eps\to 0$, 
\begin{eqnarray*}
\int_{\Omega}h(x)\ue^2\, dx &=&  \int_{B_\delta(0)}h(x)\Ue^2\, dx+ O\left(\eps^{\bps-\bms}\right)\\
&= &\eps^2\int_{B_{\eps^{-1}\delta}(0)}h(\eps x)U^2\, dx+ O\left(\eps^{\bps-\bms}\right).
\end{eqnarray*}
If $\gamma<\frac{(n-2)^2}{4}-1$ and $\bps-\bms>2$, the extremal $U\in L^2(\rn)$ and therefore
\begin{equation*}
\int_{\Omega}h(x)\ue^2\, dx = h(0)\int_{\rn}U^2\, dx \eps^2+o(\eps^2)\quad \hbox{as $\eps\to 0$. }
\end{equation*}
 If now $\gamma=\frac{(n-2)^2}{4}-1$, then $U(x)\simeq_{x\to \infty} |x|^{-\frac{n}{2}}$ and $\bps-\bms=2$. Therefore
\begin{equation*}
\int_{\Omega}h(x)\ue^2\, dx = h(0)\omega_{n-1} \eps^2\ln\left(\frac{1}{\eps}\right)+o(\eps^2\ln\eps)\quad \hbox{as $\eps\to 0$}, 
\end{equation*}
where $\omega_{n-1}$ is the volume of the canonical $(n-1)-$sphere. Combining the above estimates as $\eps\to 0$ yields  
\begin{equation*}
J_{\gamma,s,h}^\Omega(\ue)= J_{\gamma,s,0}^{\rn}(U)-\left\{\begin{array}{ll}
\frac{h(0)\Vert U\Vert_2^2}{\Vert U\Vert_{\crits, |x|^{-s}}^2} \eps^2+o(\eps^2)&\hbox{ if }\gamma<\frac{(n-2)^2}{4}-1\\
\frac{h(0)\omega_{n-1}}{\Vert U\Vert_{\crits, |x|^{-s}}^2} \eps^2\ln\left(\frac{1}{\eps}\right)+o(\eps^2\ln\eps)&\hbox{ if }\gamma=\frac{(n-2)^2}{4}-1.
\end{array}\right.
\end{equation*}
In either case, if $h(0)=\lambda >0$, we get that 
$$\mu_{\gamma, s, \lambda}(\Omega)\leq J_{\gamma,s,\lambda}^\Omega(\ue) <J_{\gamma,s,0}^{\rn}(U)= \mu_{\gamma, s,0}(\rn),$$
and we are done. \medskip

\newpage
\noindent {\bf Subsection 2.2: The merely singular case}\\

We now prove the second part of Theorem \ref{main.interior}. Assuming that $s=0$, $\gamma <0$ and $n\geq 4$, we shall prove that $\mu_{\gamma, s, \lambda}(\Omega)$ is attained if and only if $\lambda^\star(\gamma,\Omega)<\lambda$, where  
$$\lambda^\star(\gamma, \Omega)=\inf \left\{\frac{|\gamma|}{|x|^2}; \, x\in \Omega\right\}< \lambda_1(L_\gamma).$$
Note that in this case, we have $\mu_{0, 0,0}(\rn)=\mu_{\gamma, 0,0}(\rn)$ as noted in \cite{gr4}, that is 
$$\inf_{U\in\dundeuxrn\setminus\{0\}}J_{\gamma,0,0}^{\rn}(U)=\frac{1}{K(n,2)^2}:=\inf_{U\in\dundeuxrn\setminus\{0\}}\frac{\int_{\rn}|\nabla U|^2\, dx}{\left(\int_{\rn}|U|^{\crit}\, dx\right)^{\frac{2}{\crit}}},$$
and the infimum of $\mu_{\gamma, 0}(\rnp)$ is not achieved. 
Consider now the following known extremal for $\mu_{0,0}(\rn)$, 
$$U(x):=\frac{1}{\left(1+|x|^{2}\right)^{\frac{n-2}{2}}} \hbox{ for }x\in\rn.$$
Fix $x_0\in\Omega$, $x_0\neq 0$, and define the test-function
\begin{equation}\label{scaled}
u_\eps(x):=\eta(x) \eps^{-\frac{n-2}{2}}U(\eps^{-1}(x-x_0))\hbox{ for all }x\in \Omega,
\end{equation}
where $\eta\in C^\infty_c(\Omega)$ is such that $\eta(x)=1$ around $x_0\in\Omega$. A straightforward computation yields
$$J_{\gamma,0, \lambda}^\Omega(u_\eps)=\mu_{0,0}(\rn)+o(1)\hbox{ as }\eps\to 0,$$
which yields that $\mu_{\gamma, 0, \lambda}(\Omega)\leq \mu_{0,0}(\rn)$.
Classical computations in the spirit of Aubin \cite{aubin}, which can be done by replacing $-\lambda$ with a more general function $h$, yield that for $n\geq 4$, there exists   $c_n>0$, such that as $\eps\to 0$, 
\begin{equation*}
J_{\gamma,h}^\Omega(\ue)=\frac{1}{K(n,2)^2}+\left\{\begin{array}{ll}
c_n \left(\frac{|\gamma|}{|x_0|^2}-h(x_0)\right)\eps^2+o(\eps^2)& \hbox{ if }n\geq 5\\
c_4 \left(\frac{|\gamma|}{|x_0|^2}-h(x_0)\right)\eps^2\ln(\eps^{-1})+O(\eps^2)& \hbox{ if }n=4.
\end{array}\right.
\end{equation*}
Therefore, if $n\geq 4$ and assuming there exists $x_0\in\Omega\setminus\{0\}$ such that $h(x_0)>\frac{|\gamma|}{|x_0|^2}$, we obtain that $\inf_{u\in\huno\setminus\{0\}}J_{\gamma,h}^{\Omega}(u)<\inf_{U\in\dundeuxrn\setminus\{0\}}J_{\gamma,0}^{\rn}(U)$, and $\mu_{\gamma, h}(\Omega)$ is attained. 

\medskip\noindent Conversely, if $h(x)\leq \frac{|\gamma|}{|x|^2}$ for all $x\in\Omega\setminus\{0\}$, then $-\left(\frac{\gamma}{|x|^2}+h(x)\right)\geq 0$ for all $x\in\Omega\setminus\{0\}$, hence
 $\mu_{\gamma, 0, h}(\Omega)\geq \mu_{0,0,0}(\Omega)=\mu_{0,0,0}(\rn)$. We therefore have equality, and there is no extremal  for $\mu_{\gamma, 0, h}(\Omega)$ since the extremals on $\rn$ are rescaled and translated versions of $U$.

\section{\, The Hardy-singular interior mass of a domain}\label{sec:mass}
This section is devoted to the construction of the singular interior mass, as stated in Theorem \ref{interior.mass}. We start with the following key result. 

\begin{proposition}\label{prop:sol:sing} Assume $\Omega$ is a smooth bounded domain in $\rn$ and let $h\in C^{0,\theta}(\overline{\Omega})$ with $\theta\in (0,1)$. If the operator $-\Delta -\frac{\gamma}{|x|^2}-h(x)$ is coercive, then there exists a solution $H\in C^\infty(\overline{\Omega}\setminus \{0\})$ for the linear problem
\begin{equation}\label{eq:H}
\left\{\begin{array}{ll}
-\Delta H-\left(\frac{\gamma}{|x|^2}+h(x)\right)H=0 &\hbox{ in }\Omega\setminus\{0\}\\
\hfill H>0 &\hbox{ in }\Omega\setminus\{0\}\\
\hfill H=0 &\hbox{ on }\partial\Omega.
\end{array}\right.
\end{equation}
Moreover, there exists $c>0$ such that 
\begin{equation}\label{estimate1}
H(x)\simeq_{x\to 0}\frac{c}{|x|^{\bp}}.
\end{equation}
If $H'\in C^\infty(\overline{\Omega}\setminus \{0\})$ is another solution for   \eqref{eq:H}, then there exists $\lambda>0$ such that $H'=\lambda H$. \end{proposition}

\smallskip\noindent {\it Proof:}  
First, we prove existence of a solution. For that, 
let $\eta_1\in C^\infty(\rr)$ be such that $\eta_1(t)=0$ for $t<1$ and $\eta_1(t)=1$ for $t>2$. For $\epsilon>0$, set $\eta_\eps(x):=\eta_1(|x|/\eps)$ for all $x\in \rn$. Then let $\bar{H}_\eps\in C^\infty(\overline{\Omega}\setminus\{0\})$ be the Green's function for the operator $-\Delta-\eta_\eps(x)\left(\frac{\gamma}{|x|^2}+h(x)\right)$
that is singular at $0$. In particular, we have that
$$\left\{\begin{array}{ll}
-\Delta \bar{H}_\eps-\eta_\eps(x)\left(\frac{\gamma}{|x|^2}+h(x)\right)\bar{H}_\eps=0 &\hbox{ in }\Omega\setminus\{0\}\\
\hfill \bar{H}_\eps>0 &\hbox{ in }\Omega\setminus\{0\}\\
\hfill \bar{H}_\eps=0 &\hbox{ on }\partial\Omega.
\end{array}\right.
$$
Fix $x_0\in\Omega\setminus \{0\}$ and define $H_\eps(x):=\frac{\bar{H}_\eps(x)}{\bar{H}_\eps(x_0)}$ for all $x\in \overline{\Omega}\setminus\{0\}$. For $\delta>0$ such that $B_{4\delta}(0)\subset \Omega$ and $\delta<|x_0|/4$, we take $\epsilon\in (0,\delta/2)$. We then have
$$\left\{\begin{array}{ll}
-\Delta H_\eps-\left(\frac{\gamma}{|x|^2}+h(x)\right)H_\eps=0 &\hbox{ in }\Omega\setminus B_\delta(0)\\
\hfill H_\eps>0 &\hbox{ in }\Omega\setminus B_{\delta}(0)\\
\hfill  H_\eps=0 &\hbox{ on }\partial\Omega.
\end{array}\right.
$$
It follows from the boundary Harnack inequality (see for instance Ghoussoub-Robert \cite{gr4}, Proposition 7.2) that there exists $C_\delta>0$ such that 
$$\frac{H_\eps(x)}{d(x,\partial\Omega)}\leq C_\delta \frac{H_\eps(x_0)}{d(x_0,\partial\Omega)}=\frac{C_\delta}{d(x_0,\partial\Omega)}\hbox{ for all }x\in \Omega\setminus B_{2\delta}(0).$$
Since this is valid for any $\delta>0$ small enough, it then follows from standard ellitpic theory that there exists $H\in C^\infty(\overline{\Omega}\setminus \{0\})$ such that $H_\eps\to H$ in $C^k_{loc}(\overline{\Omega}\setminus \{0\})$ as $\epsilon\to 0$ for all $k\in\nn$. In particular, we have 
$$\left\{\begin{array}{ll}
-\Delta H-\left(\frac{\gamma}{|x|^2}+h(x)\right)H=0 &\hbox{ in }\Omega\setminus\{0\}\\
\hfill H\geq 0 &\hbox{ in }\Omega\setminus\{0\}\\
\hfill H=0 &\hbox{ on }\partial\Omega.
\end{array}\right.$$
Since $H(x_0)=1$, it follows from the strong maximum principle that $H>0$, hence it satisfies \eqref{eq:H}. \\
 It then follows from Theorem \ref{th:sing} that there exists $c>0$ such that
$$\hbox{either }H(x)\simeq_{x\to 0}\frac{c}{|x|^{\bp}}\hbox{ or }H(x)\simeq_{x\to 0}\frac{c}{|x|^{\bm}}.$$
If \eqref{estimate1} does not hold, we the second case holds and $H\in \huno$: since $-\Delta-\gamma|x|^{-2}-h$ is coercive, equation \eqref{eq:H} then yield $H\equiv 0$, contradicting $H>0$. This proves \eqref{estimate1}.

\medskip\noindent To prove uniqueness, let $H'\in C^\infty(\Omegabar\setminus\{0\})$ be such that \eqref{eq:H} holds. Set $\lambda_0:=\max\{\lambda\geq 0;\,  H'\geq \lambda H\}$. This is well defined and we set $\tilde{H}:=H'-\lambda_0 H$. Then $\tilde{H}\geq 0$ satisfies $-\Delta \tilde{H}-\left(\frac{\gamma}{|x|^2}+h(x)\right)\tilde{H}=0$ in $\Omega\setminus\{0\}$. Therefore, if $\tilde{H}\not\equiv 0$, it follows from the maximum principle that $\tilde{H}>0$. Then the asymptotic control \eqref{estimate1} and Hopf's boundary comparison principle yield the existence of $\epsilon_0>0$ such that $\tilde{H}\geq \eps_0 H$ in $\Omega\setminus\{0\}$, contradicting the definition of $\lambda_0$. Therefore $\tilde{H}\equiv 0$ and $H'=\lambda_0 H$, which proves the uniqueness statement.
\hfill$\Box$

We now proceed with the proof of Theorem \ref{interior.mass}.

\begin{proposition}\label{prop:def:mass} Let $\Omega$ be a smooth bounded domain in $\rn$ and fix $h \in C^{0,\theta}(\overline{\Omega})$, $\theta\in (0,1)$. Assume that the operator $-\Delta -\left(\frac{\gamma}{|x|^2}+ h(x)\right)$ is coercive and that $\gamma>\frac{(n-2)^2}{4}-1$. If $H\in C^\infty(\Omegabar\setminus\{0\})$ is a solution to \eqref{eq:H}, then there exist $c_1,c_2\in \rr$ with $c_1>0$ such that
\begin{equation}\label{exp:H}
H(x)=\frac{c_1}{|x|^{\bp}}+\frac{c_2}{|x|^{\bm}}+o\left(\frac{1}{|x|^{\bm}}\right)\hbox{ as }x\to 0.
\end{equation}
The ratio $\frac{c_2}{c_1}\in\rr$ is independent of the choice of $H$. We can therefore define the {\it mass} as $m_{\gamma, h}(\Omega):=\frac{c_2}{c_1}$. 
\end{proposition}

\smallskip{\it Proof:} Let $\eta\in C^\infty_c(\Omega)$ be such that $\eta(x)\equiv 1$ around $0$. Our first objective is to write $H(x):=\frac{\eta(x)}{|x|^{\bp}}+f(x)$ for some $f\in \huno$. For that, we consider the function
$$g(x):=-\left(-\Delta -\left(\frac{\gamma}{|x|^2}+h(x)\right)\right)(\eta |x|^{-\bp})\hbox{ in }\Omega\setminus\{0\}.$$
Since $\eta(x)=1$ around $0$, we have that
\begin{equation}\label{up:g}
|g(x)|=|\frac{h(x)}{|x|^{\bp}}|\leq C|x|^{-\bp}\hbox{ around }0.
\end{equation}
Therefore, $g\in L^{\frac{2n}{n+2}}(\Omega)$ if and only if $\bp<\frac{n+2}{2}$, which holds if and only if $\gamma>\frac{(n-2)^2}{4}-4$. The latter is guaranteed by our assumption on $\gamma$. \\
Since $L^{\frac{2n}{n+2}}(\Omega)=L^{\frac{2n}{n-2}}(\Omega)^\prime\subset \huno^\prime$, there exists $f\in \huno$ such that
$$-\Delta f -\left(\frac{\gamma}{|x|^2}+h(x)\right) f=g\hbox{  in }\huno.$$
By regularity theory, we have that $f\in C^2(\Omegabar\setminus\{0\})$. We now show that 
\begin{equation} \label{finitelimit}
\hbox{$|x|^{\bm} f(x)$ has a finite limit as $x\to 0$.}
\end{equation}
Define
$$H(x):=\frac{\eta(x)}{|x|^{\bp}}+f(x)\hbox{ for all }x\in\Omegabar\setminus\{0\}.$$
and note that $H\in C^2(\Omegabar\setminus\{0\})$ and is a solution to $\left(-\Delta -\frac{\gamma+h(x)}{|x|^2}\right)H=0$. 

\medskip\noindent Write $g_+(x):=\max\{g(x),0\}$ and $g_-(x):=\max\{-g(x),0\}$ so that $g=g_+-g_-$, and let $f_1,f_2\in \huno$ be weak solutions to
\begin{equation}\label{eq:f1}
-\Delta f_1 -\left(\frac{\gamma}{|x|^2}+h(x)\right) f_1=g_+\hbox{ and } -\Delta f_2-\left(\frac{\gamma}{|x|^2}+h(x)\right) f_2=g_-\hbox{ in }\huno.
\end{equation}
In particular, uniqueness, coercivity and the maximum principle yields $f=f_1-f_2$ and $f_1,f_2\geq 0$.

Assume that $f_1\not\equiv 0$, so that $f_1>0$ in $\Omega\setminus \{0\}$,  fix $\alpha\in (\bm,\bp)$, choose $\mu\in\rr$ such that $\mu(\bm(n-2-\bm)-\gamma)<0$, and define $u_-(x):=|x|^{-\bm}+\mu|x|^{-\alpha}$ for all $x\neq 0$. As in the proof of the previous proposition, we get that for some $\delta>0$ small, 
$$\left(-\Delta-\frac{\gamma+h(x)}{|x|^2}\right)u_-<0 \quad \hbox{ for $x\in B_\delta(0)\setminus \{0\},$}
$$
that is $u_-$ is a subsolution on $B_\delta(0)\setminus \{0\}$.

 Fix now $C>0$ such that $f_1\geq C u_-$ on $\partial B_\delta(0)$. Since $f_1$ and $Cu_-\in \huno$ are respectively super- and sub-solutions to $\left(-\Delta-\frac{\gamma+h(x)}{|x|^2}\right)u=0$, it follows from the comparison principle (via coercivity) that $f_1\geq C u_-$ in $B_\delta(0)\setminus\{0\}$, and therefore $f_1\geq C' |x|^{-\bm}$ in $B_\delta(0)\setminus\{0\}$. It then follows from \eqref{up:g} that $g_+(x)\leq |g(x)|\leq C_1 |x|^{-2+ (2-(\bp-\bm))}f_1$, and therefore \eqref{eq:f1} yields
\begin{equation}\label{eq:f12}
\left(-\Delta -\frac{\gamma+O(|x|^{2-(\bp-\bm)})}{|x|^2}\right) f_1=0\hbox{ weakly in }\huno.
\end{equation}
Since $\gamma>\frac{(n-2)^2}{4}-1$ if and only if $\tau:=2-(\bp-\bm)>0$, we can  argue as in the proof of Proposition \ref{prop:sol:sing} (see also the regularity Theorem \ref{th:hopf}) and get that $|x|^{\bm} f_1(x)$ has a finite limit as $x\to 0$.  
Similarly, $|x|^{\bm} f_2(x)$ has also a finite limit as $x\to 0$, and therefore (\ref{finitelimit}) is verified.\\
 It follows that there exists $c_2\in\rr$ such that 
$$H(x)=\frac{1}{|x|^{\bp}}+\frac{c_2}{|x|^{\bm}}+o\left(\frac{1}{|x|^{\bm}}\right)\hbox{ as }x\to 0,$$
which proves the existence of a solution $H$ to the problem with the relevant asymptotic behavior. The uniqueness result of Proposition \ref{prop:sol:sing} then yields the conclusion. \hfill $\Box$

The following proposition summarizes the properties of the mass.

\begin{proposition} Let $\Omega$ be a smooth bounded domain in $\rn$ and fix $h \in C^{0,\theta}(\overline{\Omega})$ with $\theta\in (0,1)$. Assume that the operator $-\Delta -\left(\frac{\gamma}{|x|^2}+ h(x)\right)$ is coercive and that $\gamma>\frac{(n-2)^2}{4}-1$. The mass 
 $m_{\gamma, h}(\Omega)$ then satisfies the following properties:
\begin{enumerate}[i)]
\item $m_{\gamma, 0}(\Omega)<0$,
\item If $h\leq h'$ and $h\not\equiv h'$, then $m_{\gamma, h}(\Omega)<m_{\gamma, h'}(\Omega)$,
\item If $\Omega'\subsetneq \Omega$, then $m_{\gamma, h}(\Omega')<m_{\gamma, h}(\Omega)$.
\item The function $h\mapsto m_{\gamma, h}(\Omega)$ is continuous for the  $C^0(\overline{\Omega})$ norm.
\end{enumerate}
\end{proposition}
\medskip\noindent{\it Proof:} For any such $h\in C^{0,\theta}(\overline{\Omega})$, we let $H_h$ be the unique solution to \eqref{eq:H} such that \eqref{exp:H} holds with $c_1=1$. In other words,
$$H_h(x)=\frac{1}{|x|^{\bp}}+\frac{m_{\gamma,h}(\Omega)}{|x|^{\bm}}+o\left(\frac{1}{|x|^{\bm}}\right)\hbox{ as }x\to 0.$$
Since $(-\Delta -\left(\frac{\gamma}{|x|^2}+ h(x)\right))(H_0(x)-|x|^{-\bp})=0$ and is negative on $\partial\Omega$, it follows from the maximum principle that $H_0(x)-|x|^{-\bp}<0$ on $\Omega$. It then follows from Theorem \ref{th:hopf} that $m_{\gamma,0}(\Omega)<0$. This prove property (i) of the proposition. Property (iii) goes similarly.

\smallskip\noindent For (ii), we define $g:=H_h-H_{h'}$. We have that $g\in \huno$ and $-\Delta g-\left(\frac{\gamma}{|x|^2}+ h(x)\right)g=(h-h')H_{h'}\leq 0$, but $\not\equiv 0$. Therefore $g<0$, and it follows from Theorem \ref{th:hopf} in Appendix C that there exists $K>0$ such that $g(x)|x|^{\bm}\to -K$ as $x\to 0$, and therefore $m_{\gamma, h}(\Omega)-m_{\gamma, h'}(\Omega)=-K<0$, which proves the second part of the proposition.\qed

\section{\, Positive mass and the existence of extremals in lower dimensions}

In this section, we show how the positivity of the Hardy-singular mass $m_{\gamma, \lambda}(\Omega)$ in the truly singular case (resp., the mass in the merely singular case) yields that $\mu_{\gamma,s, \lambda}(\Omega)$ is attained in the corresponding low dimensions, i.e., $\frac{(n-2)^2}{4}-1  <\gamma <\frac{(n-2)^2}{4}$ in the truly singular case, and $n=3$ in the merely singular case.\vskip 10pt

\noindent {\bf Subsection 4.1: The truly singular case}

\begin{proposition} \label{positivemass:1} Let $\Omega$ be a smooth bounded domain in $\rn$ ($n\geq 3$) such that $0\in  \Omega$.  
 Assume either  $0<s<2$ or that \{$s=0$, $\gamma > 0$\}.  If $\frac{(n-2)^2}{4}-1  <\gamma <\frac{(n-2)^2}{4}$ and $0<\lambda <\lambda_1(L_\gamma)$ is such that the mass $m_{\gamma, \lambda}(\Omega)$ is positive,  then  $\mu_{\gamma,s, \lambda}(\Omega)$ is attained.
\end{proposition}

\medskip\noindent{\it Proof:} Assuming that $\frac{(n-2)^2}{4}-1<\gamma <\frac{(n-2)^2}{4}$, we know that the mass $m_{\gamma, \lambda}(\Omega)$ is defined. 
We need to show that if $m_{\gamma, \lambda}(\Omega) >0$, then  $\mu_{\gamma,s, \lambda}(\Omega)<\mu_{\gamma,s,0}(\rn)$. \\
Consider again for each $\epsilon>0$ the extremals
\begin{equation}\label{eq:Ue}
\Ue(x):=\eps^{-\frac{n-2}{2}}U\left(\frac{x}{\eps}\right)=\left(\frac{\eps^{\frac{2-s}{n-2}\cdot\frac{\bps-\bms}{2}}}{\eps^{\frac{2-s}{n-2}\cdot(\bps-\bms)}|x|^{\frac{(2-s)\bms}{n-2}}+ |x|^{\frac{(2-s)\bps}{n-2}}}\right)^{\frac{n-2}{2-s}}. 
\end{equation}
We shall first replace $\lambda$ with any function $h\in C^{0,\theta}(\overline{\Omega})$, where $\theta\in (0,1)$ and the operator $-\Delta -\left(\frac{\gamma}{|x|^2}+ h(x)\right)$ is coercive . Consider again a test function $\eta\in C^\infty_c(\Omega)$  such that $\eta(x)=1$ for $x\in\Omega$ in a neighborhood of $0$. Since $\gamma>\frac{(n-2)^2}{4}-1$, it follows from Proposition \ref{prop:def:mass} that there exists $\beta\in \huno$ such that 
\begin{equation}\label{asymp:beta}
\beta(x)\simeq_{x\to 0}\frac{m_{\gamma, h}(\Omega)}{|x|^{\bms}},
\end{equation}
and the function $H(x):=\frac{\eta(x)}{|x|^{\bps}}+\beta(x)$ satisfies 
\begin{equation}\label{eq:H:bis}
\left\{\begin{array}{ll}
-\Delta H-\left(\frac{\gamma}{|x|^2}+h(x)\right)H=0 &\hbox{ in }\Omega\setminus\{0\}\\
\hfill H>0 &\hbox{ in }\Omega\setminus\{0\}\\
\hfill H=0 &\hbox{ on }\partial\Omega.
\end{array}\right.
\end{equation}
\medskip\noindent Define now
$$\ue(x):=\eta(x)\Ue(x)+\eps^{\frac{\bps-\bms}{2}}\beta(x)\hbox{ \quad for }x\in \Omegabar\setminus\{0\}. $$
It is clear that $\ue\in\huno$ for all $\eps>0$. We now estimate $J_{\gamma,s,h}^{\Omega}(\ue)$, where again $J_{\gamma,s,h}^{\Omega}$ is the functional on $\huno$ defined by
$$J_{\gamma,s,h}^{\Omega}(u):=\frac{\int_{\Omega}\left(|\nabla u|^2-\left(\frac{\gamma}{|x|^2}+h(x)\right)u^2\right)\, dx}{\left(\int_{\Omega}\frac{|u|^{\crits}}{|x|^s}\, dx\right)^{\frac{2}{\crits}}}.$$
Thereafter, the notation $``o_\delta(1)"$ will mean
$\lim_{\delta\to 0}\lim_{\eps\to 0}o_\delta(1)=0.$

\medskip\noindent {\bf Step 1: Estimating $\int_{\Omega}\left(|\nabla \ue|^2-\left(\frac{\gamma}{|x|^2}+h(x)\right)\ue^2\right)\, dx$}\\

\noindent Letting $\delta\in (0, {\rm dist}(0,\partial\Omega))$, we start by estimating $\int_{\Omega\setminus B_\delta(0)}\left(|\nabla \ue|^2-\left(\frac{\gamma}{|x|^2}+h(x)\right)\ue^2\right)\, dx$.\\
First note that 
\begin{equation}\label{lim:ue:H}
\lim_{\eps\to 0}\frac{\ue}{\eps^{\frac{\bps-\bms}{2}}}=H\hbox{ in }C^2_{loc}(\Omegabar\setminus\{0\}).
\end{equation}
Therefore,
\begin{equation}\label{lim:num:1}
\lim_{\eps\to 0}\frac{\int_{\Omega\setminus B_\delta(0)}\left(|\nabla \ue|^2-\left(\frac{\gamma}{|x|^2}+h(x)\right)\ue^2\right)\, dx}{\eps^{\bps-\bms}}=\int_{\Omega\setminus B_\delta(0)}\left(|\nabla H|^2-\left(\frac{\gamma}{|x|^2}+h(x)\right)H^2\right)\, dx.
\end{equation}
Integrating by parts and using equation \eqref{eq:H:bis} yields
\begin{eqnarray}
\int_{\Omega\setminus B_\delta(0)}\left(|\nabla H|^2-\left(\frac{\gamma}{|x|^2}+h(x)\right)H^2\right)\, dx &=& \int_{\Omega\setminus B_\delta(0)}H\left(-\Delta H-\left(\frac{\gamma}{|x|^2}+h(x)\right)H\right)\, dx\nonumber \\
&&+\int_{\partial (\Omega\setminus B_\delta(0))}H\partial_\nu H\, d\sigma\nonumber\\
&=&-\int_{\partial  B_\delta(0)}H\partial_\nu H\, d\sigma. \label{est:int:1}
\end{eqnarray}
Since $\bps+\bms=n-2$, using elliptic estimates, and the definition of $H$ yields
$$H\partial_\nu H=-\bps|x|^{-2\bps-1}-(n-2)m_{\gamma,h}(\Omega)|x|^{-(n-1)}+o(|x|^{-(n-1)})\quad \hbox{as $x\to 0$.} $$
Therefore, plugging this expansion into \eqref{lim:num:1} and \eqref{est:int:1} yields
\begin{equation}\label{lim:num:2}
\int_{\Omega\setminus B_\delta(0)}\left(|\nabla \ue|^2-\frac{\gamma+h(x)}{|x|^2}\ue^2\right)\, dx= \eps^{\bps-\bms}\omega_{n-1}\left(\frac{\bps}{\delta^{\bps-\bms}}+(n-2)m_{\gamma,h}(\Omega)+o_\delta(1)\right)
\end{equation}
We now deal with the expression
$\int_{ B_\delta(0)}\left(|\nabla \ue|^2-\left(\frac{\gamma}{|x|^2}+h(x)\right)\ue^2\right)\, dx.$\\
Take $\delta>0$ small enough such that $\eta(x)=1$ for $x\in B_\delta(0)\subset \Omega$. Therefore, $\ue(x)=\Ue(x)+\eps^{\frac{\bps-\bms}{2}}\beta(x)$ for $x\in B_\delta(0)$ and then
\begin{eqnarray*}
\int_{B_\delta(0)}\left(|\nabla \ue|^2-\left(\frac{\gamma}{|x|^2}+h(x)\right)\ue^2\right)\, dx&=& \int_{B_\delta(0)}\left(|\nabla \Ue|^2-\left(\frac{\gamma}{|x|^2}+h(x)\right)\Ue^2\right)\, dx\nonumber\\
&&+2\eps^{\frac{\bps-\bms}{2}}\int_{B_\delta(0)}\left(\nabla\Ue\nabla\beta-\frac{\gamma+h(x)}{|x|^2}\Ue\beta\right)\, dx \nonumber\\
&&+\eps^{\bps-\bms}\int_{B_\delta(0)}\left(|\nabla \beta|^2-\frac{\gamma+h(x)}{|x|^2}\beta^2\right)\, dx.
\end{eqnarray*}
Since $\Ue,\beta\in \huno$ and $\Ue$ is explicit, we integrate by parts the first and second term of the right-hand-side and we neglect the third term to get
\begin{eqnarray}
\int_{B_\delta(0)}\left(|\nabla \ue|^2-\left(\frac{\gamma}{|x|^2}+h(x)\right)\ue^2\right)\, dx&=& \int_{B_\delta(0)}\left(-\Delta \Ue-\frac{\gamma}{|x|^2}\Ue\right)\Ue\, dx\nonumber\\
&&+\int_{\partial B_\delta(0)}\Ue\partial_\nu\Ue\, d\sigma -\int_{B_\delta(0)}h(x)\Ue^2\, dx \nonumber\\
&&+2\eps^{\frac{\bps-\bms}{2}}\int_{B_\delta(0)}\left(-\Delta \Ue-\left(\frac{\gamma}{|x|^2}+h(x)\right)\Ue\right)\beta\, dx\nonumber\\
&&+2\eps^{\frac{\bps-\bms}{2}}\int_{\partial B_\delta(0)}\beta\partial_\nu\Ue\, d\sigma
+o_\delta(\eps^{\bps-\bms}).\label{exp:0}
\end{eqnarray}
We estimate the terms of the right-hand-side separately. 

Note first that \eqref{eq:U} and \eqref{eq:Ue} yield that as $\eps\to 0$, 
\begin{eqnarray}
\int_{B_\delta(0)}\left(-\Delta \Ue-\frac{\gamma}{|x|^2}\Ue\right)\Ue\, dx&=& \lambda \int_{B_\delta(0)}\frac{\Ue^{\crits}}{|x|^s}\, dx=\chi \int_{B_{\frac{\delta}{\eps}}(0)}\frac{U^{\crits}}{|x|^s}\, dx\nonumber \\
&=&\chi \int_{\rn}\frac{U^{\crits}}{|x|^s}\, dx +O\left(\int_{\rn\setminus B_{\eps^{-1}\delta}(0)}|x|^{-\bps \crits-s}\, dx\right)\nonumber \\
&=&\chi \int_{\rn}\frac{U^{\crits}}{|x|^s}\, dx +o\left(\eps^{\bps-\bms}\right).\label{exp:1}
\end{eqnarray}
The explicit expression of $\Ue$ in \eqref{eq:Ue} yields
\begin{equation}
\int_{\partial B_\delta(0)}\Ue\partial_\nu\Ue\, d\sigma = -\bps\omega_{n-1}\delta^{-(\bps-\bms)}\eps^{\bps-\bms}+o_\delta(\eps^{\bps-\bms}) \quad \hbox{as $\eps\to 0$.}\label{exp:2}
\end{equation}
Since $\gamma>\frac{(n-2)^2}{4}-1$, we have that $2\bps<n$ and therefore 
\begin{equation}\label{exp:3}
\int_{B_\delta(0)}h(x)\Ue^2\, dx=O\left(\eps^{\bps-\bms}\int_{B_\delta(0)}|x|^{-2\bps}\, dx\right)=o_\delta\left(\eps^{\bps-\bms}\right).
\end{equation}
Since $\bps+\bms=n-2<n$, we also have that 
\begin{eqnarray}
\eps^{\frac{\bps-\bms}{2}}\int_{B_\delta(0)}h(x)\Ue\beta\, dx&=&O\left(\eps^{\bps-\bms}\int_{B_\delta(0)}|x|^{-\bps}|x|^{-\bms}\, dx\right)\nonumber\\
&=&o_\delta\left(\eps^{\bps-\bms}\right).\label{exp:4}
\end{eqnarray}
It follows from \eqref{eq:U} and \eqref{eq:Ue} that
\begin{eqnarray}
\int_{B_\delta(0)}\left(-\Delta \Ue-\frac{\gamma}{|x|^2}\Ue\right)\beta\, dx&=&  \chi \int_{B_\delta(0)}\frac{\Ue^{\crits-1}}{|x|^s}\beta\, dx\nonumber\\
&=& \chi \eps^{\frac{\bps-\bms}{2}}\int_{B_{\eps^{-1}\delta}(0)}\frac{U^{\crits-1}}{|x|^s}\eps^{\bms}\beta(\eps x)\, dx=O\left(\eps^{\frac{\bps-\bms}{2}}\right).\label{exp:other}
\end{eqnarray}
Finally, using the expression \eqref{eq:Ue} of $\Ue$ and the asymptotics \eqref{asymp:beta} of $\beta$, we get that
\begin{eqnarray}
\int_{\partial B_\delta(0)}\beta\partial_\nu\Ue\, d\sigma&=& \eps^{\frac{\bps-\bms}{2}}\int_{\partial B_\delta(0)}\beta\partial_\nu |x|^{-\bps}\, d\sigma+ o\left(\eps^{\frac{\bps-\bms}{2}}\right)\nonumber\\
&=& \eps^{\frac{\bps-\bms}{2}}m_{\gamma,h}(\Omega)\int_{\partial B_\delta(0)}|x|^{-\bms}\partial_\nu |x|^{-\bps}\, d\sigma+ o\left(\eps^{\frac{\bps-\bms}{2}}\right)\nonumber\\ &=& -\eps^{\frac{\bps-\bms}{2}}m_{\gamma,h}(\Omega)\bps\omega_{n-1}+o\left(\eps^{\frac{\bps-\bms}{2}}\right).\label{exp:6}
\end{eqnarray}
Plugging together \eqref{lim:num:2}, \eqref{exp:0}, \eqref{exp:1}, \eqref{exp:2}, \eqref{exp:3}, \eqref{exp:4}, \eqref{exp:other} and \eqref{exp:6} yields
\begin{eqnarray}
&&\int_{\Omega}\left(|\nabla \ue|^2-\left(\frac{\gamma}{|x|^2}+h(x)\right)\ue^2\right)\, dx= \chi \int_{\rn}\frac{U^{\crits}}{|x|^s}\, dx\nonumber\\
&& +(n-2-2\bps) m_{\gamma,h}(\Omega)\omega_{n-1} \eps^{\bps-\bms}+ \chi \eps^{\frac{\bps-\bms}{2}} \int_{B_\delta(0)}\frac{\Ue^{\crits-1}}{|x|^s}\beta\, dx +o(\eps^{\bps-\bms}).\label{lim:num:3}
\end{eqnarray}
\noindent {\bf Step 2: Estimating $\int_\Omega \frac{\ue^{\crits}}{|x|^s}\, dx$}\\

\noindent From \eqref{lim:ue:H} and the definition \eqref{eq:Ue} of $\Ue$, we have as $\eps\to 0$ that
\begin{eqnarray*}
\int_\Omega \frac{\ue^{\crits}}{|x|^s}\, dx&=&\int_{B_\delta(0)} \frac{\ue^{\crits}}{|x|^s}\, dx+o(\eps^{\bps-\bms})\nonumber\\
&=& \int_{B_\delta(0)} \frac{|\Ue+\eps^{\frac{\bps-\bms}{2}}\beta|^{\crits}}{|x|^s}\, dx+o(\eps^{\bps-\bms})\nonumber\\
&=& \int_{B_\delta(0)} \left(\frac{\Ue^{\crits}}{|x|^s}+\crits\frac{\Ue^{\crits-1}}{|x|^s}\beta\right)\, dx\nonumber\\
&&+ \int_{B_\delta(0)}O\left(\eps^{\bps-\bms}\frac{\Ue^{\crits-2}}{|x|^s}\beta^2+\eps^{\frac{\crits}{2}(\bps-\bms)}|\beta|^{\crits}\right)\, dx+ o(\eps^{\bps-\bms})\nonumber\\
&=& \int_{B_{\eps^{-1}\delta}(0)} \frac{U^{\crits}}{|x|^s}\, dx+\crits \eps^{\frac{\bps-\bms}{2}}\int_{B_\delta(0)} \frac{\Ue^{\crits-1}}{|x|^s}\beta\, dx+ o(\eps^{\bps-\bms}).
\end{eqnarray*}
Using the expression 
of $U$, we get that
\begin{eqnarray}
\int_\Omega \frac{\ue^{\crits}}{|x|^s}\, dx&=& \int_{\rn} \frac{U^{\crits}}{|x|^s}\, dx+\crits \eps^{\frac{\bps-\bms}{2}}\int_{B_\delta(0)} \frac{\Ue^{\crits-1}}{|x|^s}\beta\, dx+ o(\eps^{\bps-\bms}).\label{lim:num:4}
\end{eqnarray}
Therefore, plugging \eqref{lim:num:3} and \eqref{lim:num:4} into $J_{\gamma,h}^\Omega(\ue)$ and using the equation \eqref{eq:U} satisfied by $U$ yields
\begin{equation*}
J_{\gamma,s,h}^\Omega(\ue)=J_{\gamma,s,0}^{\rn}(U)\left(1-\frac{2\omega_{n-1}\left(\bp-\frac{n-2}{2}\right)}{\chi\int_{\rn}\frac{U^{\crits}}{|x|^s}\, dx}m_{\gamma,h}(\Omega)\eps^{\bps-\bms}+ o(\eps^{\bps-\bms})\right).
\end{equation*}

\noindent This readily shows that if $h(x)=\lambda$, where $0 <\lambda <\lambda_1(L_\gamma)$, and if $m_{\gamma,\lambda}(\Omega)>0$, then $J_{\gamma,s,\lambda}^\Omega(\ue) < J_{\gamma,s,0}^{\rn}(U)=\mu_{\gamma, s,0}(\rn)$, and  therefore $\mu_{\gamma, s, \lambda}(\Omega)$ is attained. 
This completes the proof of Proposition \ref{positivemass:1}. \hfill $\Box$ \vskip 10pt

 \noindent {\bf Subsection 4.2: The merely singular case}

\begin{proposition} \label{positivemass:2} Let $\Omega$ be a smooth bounded domain in $\rn$, $n=3$, such that $0\in  \Omega$.  
 Assume  that $s=0$ and $\gamma < 0$.  If $0<\lambda <\lambda_1(L_\gamma)$ is such that the mass $R_{\gamma, \lambda}(\Omega)$ is positive,  then  $\mu_{\gamma,0, \lambda}(\Omega)$ is attained.
\end{proposition}

\medskip\noindent{\it Proof:} This is by now classical, so we shall sketch a proof. For any $x_0\in\Omega\setminus\{0\}$, we let $G_{x_0}\in C^\infty(\overline{\Omega}\setminus\{0\})$ be the Green's function for the operator $-\Delta-\left(\frac{\gamma}{|x|^2}+h(x)\right)$ at $x_0$ with Dirichlet boundary condition. Since $n=3$, then up to multiplying by a constant, we have 
$$G_{x_0}(x)=\frac{1}{4\pi}\left(\frac{\eta(x)}{|x-x_0|}+\beta(x)\right)\hbox{ for all }x\in\Omega\setminus\{x_0\},$$
where $\beta\in \huno\cap C^0(\Omega)$, and  there exists $R_{\gamma,h}(\Omega,x_0)\in\rr$ such that
$$G_{x_0}(x)=\frac{1}{4\pi}\left(\frac{1}{|x-x_0|}+R_{\gamma,h}(\Omega,x_0)\right)+o(1)\hbox{ as }x\to x_0.$$
Note that $\beta(x_0)=R_{\gamma,h}(x_0)$ is the Robin function at $x_0$. 

Set now $\tilde{u}_\eps(x):=\ue(x)+\sqrt{\eps}\tilde{\beta}(x)$ for all $x\in\Omega\setminus\{0\}$, where $\ue$ are the functions defined in (\ref{scaled}). Then, classical computations in the spirit of Schoen \cite{schoen1} yield
$$
J_{\gamma,h}^\Omega(\tilde{u}_\eps)=\frac{1}{K(n,2)^2}-c_3R_{\gamma,h}(\Omega,x_0)\eps^{n-2}+o\left(\eps^{n-2}\right) \quad \hbox{as $\eps\to 0$.}
$$
 If now $x_0\in\Omega\setminus \{0\}$ and $0<\lambda <\lambda_1(L_\gamma)$ are such $R_{\gamma,\lambda}(x_0)>0$, then 
\begin{equation*}
J_{\gamma,h}^\Omega(\tilde{u}_\eps)=\mu_{\gamma,0}(\rn)-c_3R_{\gamma,\lambda}(x_0)\eps +o\left(\eps\right) \quad \hbox{ as }\eps\to 0.
\end{equation*}
\noindent This implies that  $\mu_{\gamma, 0, \lambda}(\Omega)$ is attained.

\section{\, Blow-Up analysis in the truly singular case}

Let $\Omega$ be a smooth bounded domain of $\rn$, $n\geq 3$, such that $0\in \Omega$ is an interior point. Fix $\gamma<(n-2)^2/4$ and recall that 
$$\mu_{\gamma,s,0}(\rn):=\inf\left\{\frac{\int_{\rn}\left(|\nabla u|^2-\frac{\gamma}{|x|^2}u^2\right)\, dx}{\left(\int_{\rn}\frac{|u|^{\crits}}{|x|^s}\, dx\right)^{\frac{2}{\crits}}};\,  u\in \dundeuxrn\setminus\{0\}\right\},$$
where $0\leq s<2$ and $\crits:=\frac{2(n-s)}{n-2}$.\\
 Let $(a_\alpha)_{\alpha\in\nn}\in C^{1}(\Omegabar)$ be such that there exists $a_\infty\in C^{1}(\Omegabar)$ with 
\begin{equation}\label{lim:a:alpha}
\lim_{\alpha\to +\infty}a_\alpha=a_\infty\hbox{ in }C^{1}(\Omegabar).
\end{equation}
Consider $(\lambda_\alpha)_\alpha\in (0,+\infty)$ such that
\begin{equation}\label{lim:la}
\lim_{\alpha\to +\infty}\lambda_\alpha=\mu_{\gamma,s,0}(\rn).
\end{equation}
Suppose $(u_\alpha)_\alpha\in \huno$ is a sequence of weak solutions to
\begin{equation}\label{eq:ua}
\left\{\begin{array}{ll}
-\Delta \ua -\left(\frac{\gamma}{|x|^2}+\aa\right) \ua = \la \frac{\ua^{\crits-1}}{|x|^s}& \hbox{ in }\Omega\\
\ua\geq 0& \hbox{ a.e. in }\Omega\\
\ua=0& \hbox{ on }\partial \Omega
\end{array}\right.
\end{equation}
with
\begin{equation}\label{mass:ua}
\int_\Omega \frac{\ua^{\crits}}{|x|^s}\, dx=1.
\end{equation}
and such that
\begin{equation}\label{lim:weak}
\ua\rightharpoonup 0\hbox{ as }\alpha\to +\infty\hbox{ weakly in }\huno.
\end{equation}
We shall assume uniform coercivity, that is there exists $c>0$ such that
\begin{equation*}
\int_\Omega\left(|\nabla\varphi|^2-\left(\frac{\gamma}{|x|^2}+\aa\right)\varphi^2\right)\, dx\geq c\int_\Omega \varphi^2\, dx\hbox{ for all }\varphi\in \huno.
\end{equation*}
Note that this is equivalent to the coercivity of $-\Delta-(\gamma|x|^{-2}+a_\infty)$. The two following sections are devoted to the analysis of the Blow-up behavior of $(\ua)$ as $\alpha\to +\infty$. The present section deals mostly with the case $\{s>0\hbox{ or }\gamma>0\}$, for which there are extremals for $\mu_{\gamma,s,0}(\rn)$. The case $\{s=0\hbox{ and }\gamma< 0\}$ will be dealt with in the next section. The case $s=\gamma=0$ has been extensively studied in the litterature, see for instance \cites{d2,dhr} and the references therein.
\begin{theorem}\label{th:blowup:1} Let $\Omega$ be a smooth bounded domain of $\rn$, $n\geq 3$, such that $0\in \Omega$ is an interior point. Fix $\gamma<(n-2)^2/4$, and assume that either 
$s>0$ or $\gamma> 0.$
Let $(\aa)_\alpha\in C^1(\Omegabar)$, $(\la)_\alpha\in (0,+\infty)$ and $(\ua)_\alpha\in \huno$ such that \eqref{lim:a:alpha}, \eqref{lim:la}, \eqref{eq:ua}, \eqref{mass:ua} and  \eqref{lim:weak} hold. Then:
\begin{enumerate}[i)]
\item If $\gamma\leq\frac{(n-2)^2}{4}-1$, then $a_\infty(0)=0$;
\item If $\frac{(n-2)^2}{4}-1<\gamma<\frac{(n-2)^2}{4}$, then $m_{\gamma,a_\infty}(\Omega)=0$.
\end{enumerate}
In addition, there exists $C>0$ such that
\begin{equation}\label{est:c0:bis}
\ua(x)\leq C\frac{\ma^{\frac{\bp-\bm}{2}}}{\ma^{\bp-\bm}|x|^{\bm}+|x|^{\bp}}\hbox{ for all } x\in 
\Omega\setminus\{0\}, 
\end{equation}
where $\ma\to 0$ as $\alpha\to 0$ is defined in \eqref{def:ma} below.
\end{theorem}
The rest of this section is devoted to the proof of this theorem. We shall make frequent use of the following Pohozaev identity.

\begin{proposition}
Let $\Omega\subset\rn$ be a smooth bounded domain and let $u\in C^2(\overline{\Omega})$, $u\geq 0$. For any $p\in \rn$, we have \begin{eqnarray}\label{gen:p:id}
 &&  \int_{\Omega}\left((x-p)^i\partial_i u+\frac{n-2}{2}u\right)\left(-\Delta u-\frac{\gamma}{|x|^2}u-c \frac{u^{\crits}}{|x|^s}\right)\, dx  \\
&&\qquad =\int_{\partial \Omega}\left[(x-p,\nu)\left(\frac{|\nabla u|^2}{2}-\frac{\gamma u^{2}}{2|x|^2}-\frac{c u^{\crits}}{\crits|x|^s}\right)-\left((x-p)^i\partial_i u+\frac{n-2}{2}u\right)\partial_\nu u\right]\, d\sigma \nonumber\\
&& \qquad \quad+\int_\Omega \frac{(p,x)}{|x|^2} \left(\gamma \frac{u^2}{|x|^2}+c \frac{s u^{\crits}}{\crits|x|^{s}}\right)\, dx. \nonumber
\end{eqnarray}
\end{proposition}
\noindent{\it Proof:} The classical Pohozaev identity yields
\begin{eqnarray*}
&&-\int_{\Omega}\left((x-p)^i\partial_i u+\frac{n-2}{2}u\right)\Delta u\, dx\\
&& \qquad =\int_{\partial \Omega}\left[(x-p,\nu)\frac{|\nabla u|^2}{2}-\left((x-p)^i\partial_i u+\frac{n-2}{2}u\right)\partial_\nu u\right]\, d\sigma. \nonumber
\end{eqnarray*}
For any $t\in [0,2]$, integration by parts yields
\begin{eqnarray*}
&&\int_{\Omega}\left((x-p)^i\partial_i u+\frac{n-2}{2}u\right)\frac{u^{\crit(t)-1}}{|x|^t}\, dx\\
&&\qquad \qquad = -\frac{t}{\crit(t)}\int_\Omega \frac{(p,x)}{|x|^{2+t}}u^{\crit(t)}\, dx+\int_{\partial \Omega}\frac{(x-p,\nu)u^{\crit(t)}}{\crit(t)|x|^t}\, d\sigma. \nonumber
\end{eqnarray*}
Putting together these two equalities yields the general identity claimed in the proposition. \qed \\

To prove Theorem \eqref{th:blowup:1}, we start by noting that regularity theory and Theorem \ref{th:hopf} yield that for any $\alpha$, there exists $C_\alpha>0$ such that $\ua\in C^{2,\theta}(\Omegabar\setminus\{0\})$, and 
\begin{equation}\label{asymp:ua:0}
\ua(x)\sim_{x\to 0}\frac{C_\alpha}{|x|^{\bm}}\hbox{\,\, and\,\, }|\nabla\ua(x)|\leq C'_\alpha|x|^{-\bm-1}\hbox{ for }x\in \Omega\setminus\{0\}.
\end{equation}
Fix $\tau\in\rr$ such that
\begin{equation*}
\bm<\tau<\frac{n-2}{2}.
\end{equation*}
It follows from \eqref{asymp:ua:0} that for any $\alpha\in\nn$, there exists $\xa\in\Omega\setminus\{0\}$ such that
\begin{equation}\label{def:xa}
\sup_{x\in \Omega\setminus\{0\}}|x|^\tau \ua(x)=|\xa|^\tau \ua(\xa).
\end{equation}

We now prove the following proposition, which is valid for any $\gamma<(n-2)^2/4$.
\begin{proposition}\label{prop:useful}
Let $\Omega$ be a smooth bounded domain of $\rn$, $n\geq 3$, such that $0\in \Omega$ is an interior point. Fix $\gamma<(n-2)^2/4$, and let $(\aa)_\alpha\in C^1(\Omegabar)$, $(\la)_\alpha\in (0,+\infty)$ and $(\ua)_\alpha\in \huno$ be such that \eqref{lim:a:alpha}, \eqref{lim:la}, \eqref{eq:ua} and \eqref{mass:ua} hold. Let $(\xa)_\alpha\in\Omega\setminus\{0\}$ be as in \eqref{def:xa} and set 
\begin{equation}\label{def:ma}
\ma:=\ua(\xa)^{-\frac{2}{n-2}}.
\end{equation}
Then,
\begin{equation}\label{lim:ua:infty}
\lim_{\alpha\to +\infty}\sup_{x\in \Omega\setminus\{0\}}|x|^\tau \ua(x)=+\infty,
\end{equation}
and therefore $\lim_{\alpha\to+\infty}\ma=0$. In addition,
\begin{equation}\label{lim:bndy}
\lim_{\alpha\to+\infty}\frac{d(\xa,\partial\Omega)}{\ma}=+\infty.
\end{equation}
\end{proposition}
\noindent{\it Proof of Proposition \ref{prop:useful}:}  
If \eqref{lim:ua:infty} does not hold, then there exists $C>0$ such that, up to a subsequence, we have that $|x|^{\tau}\ua(x)\leq C$ for all $x\in\Omega\setminus\{0\}$. Since $\tau<\frac{n-2}{2}$, we then have that
\begin{equation}\label{lim:int:0}
\lim_{\delta\to 0}\lim_{\alpha\to +\infty} \int_{B_\delta(0)} \frac{\ua^{\crits}}{|x|^s}\, dx=0.
\end{equation}
Since $(\ua)$ is bounded uniformly in $L^\infty$ outside $0$, it then follows from \eqref{eq:ua} and \eqref{lim:weak} that $\ua\to 0$ in $C^0_{loc}(\Omegabar\setminus\{0\})$. This limit and \eqref{lim:int:0} yield $\int_{\Omega} |x|^{-s}\ua^{\crits}\, dx\to 0$ as $\alpha\to +\infty$, contradicting \eqref{mass:ua}. This proves \eqref{lim:ua:infty}.

\smallskip\noindent As a remark, note that when $s>0$, the subcriticality $\crits<\crit$ and \eqref{lim:weak} yield $\ua\to 0$ in $C^0_{loc}(\Omegabar\setminus\{0\})$.

\smallskip\noindent We now prove \eqref{lim:bndy}. Assume that $d(\xa,\partial\Omega)=O(\ma)$ as $\alpha\to +\infty$, the above remark then yields $s=0$. We let $x_\infty:=\lim_{\alpha\to +\infty}\xa$ such that $x_\infty\in\partial\Omega$. Since $\Omega$ is smooth, we let $\delta>0$ and $\varphi\in C^\infty(B_\delta(0),\rn)$ be a smooth diffeomorphism onto its image such that $\varphi(0)=x_\infty$, $\varphi(B_\delta(0)\cap\{x_1<0\})=\varphi(B_\delta(0))\cap\Omega$ and $\varphi(B_\delta(0)\cap\{x_1=0\})=\varphi(B_\delta(0))\cap\partial\Omega$. Up to a rotation and a rescaling, we can assume that $d\varphi_0=Id$. Let $((x_\alpha)_1,x'_\alpha)\in (-\infty,0)\times \rr^{n-1}\cap B_\delta(0)$ be such that $\varphi((x_\alpha)_1,x'_\alpha)=\xa$. In particular, $\lim_{\alpha\to 0}|(x_\alpha)_1|+|x'_\alpha|=0$. Define
$$\tilde{u}_\alpha(x):=\ma^{\frac{n-2}{2}}\ua\circ\varphi((0,x'_\alpha)+\ma x)\hbox{ for }x\in \{x_1<0\}\cap B_{\ma^{-1}\delta/2}(0).$$
We then have that 
$$-\Delta_{g_\alpha}\tilde{u}_\alpha-\ma^2\left(\frac{\gamma}{|\varphi((0,x'_\alpha)+\ma x)|^2}+\aa(\varphi((0,x'_\alpha)+\ma x))\right)=\lambda_\alpha\tilde{u}_\alpha^{\crit-1},$$
where $g_\alpha:=(\varphi^\star \hbox{Eucl})((0,x'_\alpha)+\ma x)$. We have that $\tilde{u}_\alpha>0$ on $\{x_1<0\}\cap B_{\ma^{-1}\delta/2}(0)$, $\tilde{u}_\alpha=0$ on $\{x_1=0\}\cap B_{\ma^{-1}\delta/2}(0)$, and $\tilde{u}_\alpha(\ma^{-1}(x_\alpha)_1,0)=1$. Therefore, standard elliptic theory yields the existence of $\tilde{u}\in C^\infty(\{x_1\leq 0\})\cap \dundeuxrn$ such that $-\Delta \tilde{u}=\mu_{\gamma,s,0}(\rn) \tilde{u}^{\crit-1}$ and $\tilde{u}>0$ in $\{x_1<0\}$ and $\tilde{u}=0$ on $\{x_1=0\}$. It follows from Theorem 1.3, chapter III in \cite{st} that this is a contradiction. This proves \eqref{lim:bndy} and ends the proof of Proposition \ref{prop:useful}.\qed

In addiction to the hypothesis of Proposition \ref{prop:useful}, we now assume that either $s>0$ or $\gamma>0$. We claim that
\begin{equation}\label{lim:xa:ma}
\lim_{\alpha\to +\infty}\frac{|\xa|}{\ma}=c>0.
\end{equation}
For that, we first show that
\begin{equation}\label{prelim}
|\xa|=O(\ma)\hbox{ as }\alpha\to +\infty.
\end{equation}
Indeed, if not we can assume that $\ma^{-1}|\xa|\to+\infty$ as $\alpha\to +\infty$. By defining $\tilde{u}_\alpha:=\ma^{\frac{n-2}{2}}\ua(\xa+\ma x)$, it follows from our assumption and Proposition \ref{prop:useful} that for any $R>0$, and for $\alpha$ large enough, $\tilde{u}_\alpha$ is defined on $B_R(0)$ and 
$$-\Delta\tilde{u}_\alpha-\left(\frac{\gamma}{\left|\frac{\xa}{\ma}+ x\right|^2}+\ma^2\aa(\xa+\ma x)\right)\tilde{u}_\alpha=\lambda_\alpha\frac{\tilde{u}_\alpha^{\crit-s}}{\left|\frac{\xa}{\ma}+ x\right|^s}\hbox{ in }B_R(0).$$
It follows from \eqref{def:xa} and the assumption that $\ma^{-1}|\xa|\to+\infty$ as $\alpha\to +\infty$,  that there exists $C(R)>0$ such that $\tilde{u}_\alpha\leq C(R)$ on $B_R(0)$ and that $\tilde{u}_\alpha(0)=1$. It then follows from standard elliptic theory that $\tilde{u}_\alpha\to \tilde{u}$ in $C^2_{loc}(\rn)$ where $0<\tilde{u}\leq 1$ and
\begin{equation}\label{eq:tu}
-\Delta\tilde{u}=\mu_{\gamma,s,0}(\rn)u^{\crit-1}\hbox{ if }s=0\hbox{ and }\Delta \tilde{u}=0\hbox{ if }s>0.
\end{equation}
By the Sobolev embedding, we have that
\begin{equation}\label{bnd:crit}
\int_{B_R(0)} \tilde{u}_\alpha^{\crit}\, dx=\int_{B_{R\ma}(x_\alpha)} u_\alpha^{\crit}\, dx\leq \int_{\Omega} u_\alpha^{\crit}\, dx\leq C \Vert\ua\Vert_{\huno}\leq C,
\end{equation}
where we used the fact that $B_{R\ma}(x_\alpha)\subset \Omega$ since \eqref{lim:bndy} holds. Therefore, by first passing to the limit as $\alpha\to +\infty$ and then as $R\to +\infty$, we get that $\tilde{u}\in L^{\crit}(\rn)$. \\
Assuming that $s>0$, and since $0<\tilde{u}\leq 1$, it  follows from \eqref{eq:tu} and Liouville's theorem that $\tilde{u}\equiv 1$, contradicting that $\tilde{u}\in L^{\crit}(\rn)$. In other words, \eqref{prelim} is proved when $s>0$.

\smallskip\noindent Assuming now that $s=0$ but $\gamma>0$, then by letting $\alpha\to +\infty$ and $R\to +\infty$ in \eqref{bnd:crit} and using\eqref{mass:ua}, we get that  $\int_{\rn} \tilde{u}^{\crit}\, dx\leq 1$. Equation \eqref{eq:tu} then yields
$$\mu_{0,0}(\rn)\leq \frac{\int_{\rn}|\nabla \tilde{u}|^2\, dx}{\left(\int_{\rn} \tilde{u}^{\crit}\, dx\right)^{\frac{2}{\crit}}}=\mu_{\gamma,0}(\rn)\left(\int_{\rn} \tilde{u}^{\crit}\, dx\right)^{\frac{2}{n}}\leq \mu_{\gamma,0}(\rn).$$
Since $\gamma>0$, it follows from classical estimates (see \cite{gr5}) that $\mu_{\gamma,0}(\rn)<\mu_{0,0}(\rn)$, yielding again a contradiction.  In other words, \eqref{prelim} is proved when $s=0$.

 \medskip\noindent We now prove \eqref{lim:xa:ma}. We argue again by contradiction and  assume that $\xa=o(\ma)$ as $\alpha\to +\infty$. We define $\tilde{u}_\alpha(x):=\ma^{\frac{n-2}{2}}\ua(|\xa|x)$ for $x\in B_{|\xa|^{-1}\delta}(0)$ and $\delta>0$ small enough. The definition \eqref{def:xa} yields $(|\xa|\cdot|x|)^\tau\ua(|\xa|x)\leq |\xa|^\tau\ua(\xa)$, and therefore $|x|^\tau\tilde{u}_\alpha(x)\leq 1$ for all $x\in B_{|\xa|^{-1}\delta}(0)\setminus\{0\}$. Equation \eqref{eq:ua} rewrites
$$-\Delta\tilde{u}_\alpha-\left(\frac{\gamma}{|x|^2}+|\xa|^2\aa(|\xa|x)\right)\tilde{u}_\alpha=\la\left(\frac{|\xa|}{\ma}\right)^{2-s}\frac{\tilde{u}_\alpha^{\crits-1}}{|x|^s}\hbox{ in }B_{|\xa|^{-1}\delta}(0)\setminus\{0\}.$$
In addition, we have that $\tilde{u}_\alpha>0$ and $\tilde{u}_\alpha(|\xa|^{-1}\xa)=1$. These estimates and standard elliptic theory then yield the existence of $\tilde{u}\in C^\infty(\rnp)$ such that $\tilde{u}_\alpha\to \tilde{u}$ in $C^2_{loc}(\rnp)$ where
$$-\Delta\tilde{u}-\frac{\gamma}{|x|^2}\tilde{u}=0\hbox{ in }\rnp\; ; \; \tilde{u}>0\; ; \; |x|^\tau\tilde{u}(x)\leq 1\hbox{ in }\rnp.$$ 
The classification of Proposition \ref{prop:funda:2} yields the existence of $A,B\geq 0$ such that $\tilde{u}(x)=A|x|^{-\bp}+B|x|^{-\bm}$ in $\rnp$. The pointwise control $|x|^\tau\tilde{u}(x)\leq 1\hbox{ in }\rnp$ yields $A=B=0$, contradicting $\tilde{u}>0$. This completes the proof of \eqref{lim:xa:ma}.

\medskip\noindent We now define $\va(x):=\ma^{\frac{n-2}{2}}\ua(\ma x)$ for $x\in \ma^{-1}\Omega\setminus\{0\}$, and claim that there exists $U\in \huno\cap C^2(\rnp)$ such that
\begin{equation}\label{lim:va:1}
\lim_{\alpha\to +\infty}\va =U\hbox{ in }H_{1,loc}^2(\rn)\cap C^2_{loc}(\rnp).
\end{equation}
For that, we first note that
\begin{equation*}
-\Delta \va-\left(\frac{\gamma}{|x|^2}+\ma^2\aa(\ma x)\right)\va=\la\frac{\va^{\crits-1}}{|x|^s}\hbox{ in }\ma^{-1}\Omega\setminus\{0\}.
\end{equation*}
Moreover, $\va>0$ and $|x|^\tau\va(x)\leq C$ for all $x\in  \ma^{-1}\Omega\setminus\{0\}$. It then follows from standard elliptic theory that there exists $U\in C^{\infty}(\rnp)$, $U\geq 0$, such that $\lim_{\alpha\to +\infty}\va=U$ in $C^{2,\theta}_{loc}(\rnp)$ and
\begin{equation}\label{eq:U:bis}
-\Delta U-\frac{\gamma}{|x|^2}U=\mu_{\gamma,s,0}(\rn)\frac{U^{\crits-1}}{|x|^s}\hbox{ in }\rnp.
\end{equation}
Since $\va(\ma^{-1}\xa)=1$, it then follows 
 that $U\not\equiv 0$, and therefore $U>0$. Moreover, we have that
$$\int_{B_R(0)\setminus B_\delta(0)}\frac{U^{\crits}}{|x|^s}\, dx=\lim_{\alpha\to +\infty} \int_{B_{R\ma}(0)\setminus B_{\delta\ma}(0)}\frac{\ua^{\crits}}{|x|^s}\, dx\leq 1.$$
Therefore, letting $R\to +\infty$ and $\delta\to 0$ yields $\int_{\rn}\frac{U^{\crits}}{|x|^s}\, dx\leq 1$. Similarly, $\int_{\rn}\frac{U^2}{|x|^2}\, dx<+\infty$ and $\int_{\rn}|\nabla U|^2\, dx<+\infty$. Therefore $U\in \dundeuxrn$, and by integrating by parts, we obtain that

%
%
$$\mu_{\gamma,s,0}(\rn)\leq \frac{\int_{\rn}\left(|\nabla U|^2-\frac{\gamma}{|x|^2}U^2\right)\, dx}{\left(\int_{\rn}\frac{U^{\crits}}{|x|^s}\, dx\right)^{\frac{2}{\crits}}}=\mu_{\gamma,s,0}(\rn)\left(\int_{\rn}\frac{U^{\crits}}{|x|^s}\, dx\right)^{\frac{2-s}{n-s}}\leq \mu_{\gamma,s,0}(\rn).$$
Therefore $\int_{\rn}\frac{U^{\crits}}{|x|^s}\, dx=1$ and $U\in \dundeuxrn$ is an extremal for $\mu_{\gamma,s,0}(\rn)$. 

\medskip\noindent We now show that
\begin{equation}\label{lim:ua:0}
\lim_{\alpha\to +\infty}\ua=0\hbox{ in }C^0_{loc}(\Omegabar\setminus\{0\}).
\end{equation}
\smallskip\noindent Indeed, when $s>0$, we have already noted that the result follows from subcriticality. If however $s=0$, it then follows from the convergence to $U$  that
\begin{eqnarray}\label{bnd:norm:residu}
&&\lim_{R\to +\infty}\lim_{\alpha\to +\infty} \int_{B_{R\ma}(0)\setminus B_{R^{-1}\ma}(0)}\frac{\ua^{\crits}}{|x|^s}\, dx\\
&&\qquad \qquad \qquad =\lim_{R\to +\infty}\lim_{\alpha\to +\infty} \int_{B_{R}(0)\setminus B_{R^{-1}}(0)}\frac{\va^{\crits}}{|x|^s}\, dx\nonumber\\
&&\qquad \qquad \qquad =\lim_{R\to +\infty}\int_{B_{R}(0)\setminus B_{R^{-1}}(0)}\frac{U^{\crits}}{|x|^s}\, dx=1.\nonumber
\end{eqnarray}
Therefore, for any $\delta>0$, we have that $\lim_{\alpha\to +\infty}\int_{\Omega\setminus B_\delta(0)}\frac{\ua^{\crits}}{|x|^s}\, dx=0$. We then rewrite \eqref{eq:ua} as $-\Delta \ua=f_\alpha \ua$ in $\Omega\setminus B_\delta(0)$ where $\lim_{\alpha\to 0}\Vert f_\alpha\Vert_{n/2}=0$. It then follows from the classical deGiorgi-Nash-Moser iterative scheme that $(\ua)$ is uniformly bounded in $C^0_{loc}(\Omegabar\setminus\{0\})$. Elliptic theory and \eqref{lim:weak} then yield the convergence to $0$. This proves \eqref{lim:ua:0}.\qed

\medskip\noindent We now claim that there exists $C>0$ such that
\begin{equation}\label{weak:est:1}
|x|^{\frac{n-2}{2}}\ua(x)\leq C\hbox{ for all }x\in\Omega\setminus\{0\} \hbox{ and }\alpha\in\nn.
\end{equation}
We argue by contradiction and we let $(\ya)_\alpha\in\Omega\setminus\{0\}$ be such that
\begin{equation}\label{def:ya}
\sup_{x\in \Omega\setminus\{0\}}|x|^{\frac{n-2}{2}} \ua(x)=|\ya|^{\frac{n-2}{2}} \ua(\ya)\to +\infty\hbox{ as }\alpha\to +\infty.
\end{equation}
Note that it follows from \eqref{asymp:ua:0} that $(\ya)_\alpha$ is well-defined, and from \eqref{def:xa}, \eqref{def:ma}, \eqref{lim:xa:ma} and \eqref{lim:ua:0} that
\begin{equation}\label{ppte:ya}
\lim_{\alpha\to +\infty}\ya=0\hbox{ , }\lim_{\alpha\to +\infty}\frac{|\ya|}{\ma}=+\infty \hbox{ and }\lim_{\alpha\to +\infty}\frac{|\ya|}{\nu_\alpha}=+\infty,
\end{equation}
where $\nu_\alpha:=\ua(\ya)^{-\frac{2}{n-2}}\to 0$ as $\alpha\to +\infty$. We define $\tilde{u}_\alpha(x):=\nu_\alpha^\frac{n-2}{2}\ua(\ya+\nu_\alpha x )$ for $x\in \nu_\alpha^{-1}\Omega\setminus\{0\}$. Equation \eqref{eq:ua} rewrites
\begin{equation}\label{eq:tua:2}
-\Delta \tilde{u}_\alpha-\left(\frac{\gamma}{|\frac{\ya}{\nu_\alpha}+x|^2}+\nu_\alpha^2\aa(\ya+\nu_\alpha x)\right)\tilde{u}_\alpha=\la\frac{\tilde{u}_\alpha^{\crits-1}}{|\frac{\ya}{\nu_\alpha}+x|^s}\hbox{ in }\nu_\alpha^{-1}\Omega\setminus\{0\}.
\end{equation}
It follows from the definition \eqref{def:ya} that for any $R>0$, $\tilde{u}_\alpha\leq 2$ in $B_R(0)$ for $\alpha>0$ large enough. Since $\tilde{u}_\alpha(0)=1$, elliptic theory yields the existence of $\tilde{u}\in C^2(\rn)$ such that $\tilde{u}_\alpha\to \tilde{u}>0$ in $C^2_{loc}(\rn)$ as $\alpha\to +\infty$. In addition, for all $R>0$, we have with Sobolev's inequality that
$$\int_{B_R(0)}\tilde{u}_\alpha^{\crit}\, dx=\int_{B_{R\nu_\alpha}(\ya)}\ua^{\crit}\, dx\leq \int_{\Omega}\ua^{\crit}\, dx\leq C$$
and therefore, letting $\alpha\to +\infty$ and $R\to +\infty$, we get that $\tilde{u}\in L^{\crit}(\rn)$. We now distinguish two cases:

\smallskip\noindent If $s>0$, then passing to the limit in \eqref{eq:tua:2}, we get that $\Delta\tilde{u}=0$ in $\rn$ and $\tilde{u}>0$ is bounded. Liouville's theorem then yields $\tilde{u}\equiv \tilde{u}(0)=1$, contradicting $\tilde{u}\in L^{\crit}(\rn)$.

\smallskip\noindent If $s=0$, then it follows from \eqref{mass:ua} and \eqref{bnd:norm:residu} that   
$$\lim_{R\to +\infty}\lim_{\alpha\to +\infty} \int_{\Omega\setminus B_{R\ma}(0)}\ua^{\crit}\, dx=0.$$
It follows from \eqref{ppte:ya} that for $\alpha>0$ large enough, then $B_{R\nu_\alpha}(\ya)\cap B_{R\ma}(0)\neq\emptyset$, and therefore, we have that $\lim_{R\to +\infty}\lim_{\alpha\to +\infty} \int_{B_{R\nu_\alpha}(\ya)}\ua^{\crit}\, dx=0$, which yields $\tilde{u}\equiv 0$, contradicting $\tilde{u}(0)=1$.
\medskip\noindent This proves \eqref{weak:est:1}.  

\medskip\noindent
We now claim that 
\begin{equation}\label{weak:est:2}
\lim_{R\to +\infty}\lim_{\alpha\to +\infty}\sup_{x\in \Omega\setminus B_{R\ma}(0)}|x|^{\frac{n-2}{2}}\ua(x)=0.
\end{equation}
We just sketch the proof which is very similar to the proof of \eqref{weak:est:1}. Arguing by contradiction and letting $(\ya)_\alpha\in\Omega$ be such that $\ma^{-1}|\ya|\to +\infty$ as $\alpha\to +\infty$ and $|\ya|^{\frac{n-2}{2}}\ua(\ya)\to c>0$. We rescale at $\ya$ and we get that our hypothesis yields the persistence of some energy outside $B_{R\ma}(0)$ for $R$ and $\alpha$ large, which is a contradiction.\qed 

\medskip\noindent 
We now prove that for any $\eps>0$ small, there exists $C_\eps>0$ such that
\begin{equation}\label{est:eps}
\ua(x)\leq C_\eps\frac{\ma^{\frac{\bp-\bm}{2}-\eps}}{|x|^{\bp-\eps}}\hbox{ for all }x\in \Omega\setminus B_{\ma}(0).
\end{equation}
Note first that in view of \eqref{weak:est:1}, it is enough to prove \eqref{est:eps} in $\Omega\setminus B_{R\ma}(0)$ for $R>0$ large.\\
For that, fix $\gamma'$ such that $\gamma<\gamma'<\frac{(n-2)^2}{4}$, and let $\Omega'$ be a smooth bounded domain of $\rn$ such that $\Omega\subset\subset \Omega'$ is relatively compact in $\Omega'$. We extend $(\aa)_\alpha$ and $a_\infty$ on $\Omega'$ such that \eqref{lim:a:alpha} holds on $\Omega'$. Let $G_{\alpha}$ be the Green's function on $\Omega'$ at $\xa$ of $-\Delta-\left(\frac{\gamma'}{|x|^2}+a_\infty+\nu\right)$, where $\nu>0$ and Dirichlet boundary condition. Up to taking $\gamma'$ close to $\gamma$, $\nu$ small enough and $\Omega'$ close to $\Omega$, the operator is coercive and the Green's function is well defined on $\overline{\Omega'}\setminus\{0,\xa\}$. Theorem \ref{th:green:gamma:domain} in Appendix A then  yields a $C>0$ such that for any $\alpha\in\nn$
\begin{equation}\label{control:G}
0<G_{\alpha}(x)\leq C\left(\frac{\max\{|\xa|,|x|\}}{\min\{|\xa|,|x|\}}\right)^{\beta_-(\gamma')}|x-\xa|^{2-n}\hbox{ for all }x\in \Omega\setminus\{0,\xa\}.
\end{equation}
Define the operator
$$L_\alpha:=-\Delta-\left(\frac{\gamma}{|x|^2}+a_\alpha\right)-\la\frac{\ua^{\crits-2}}{|x|^s}.$$
It follows from \eqref{lim:xa:ma} that there exists $R_0>1$ such that $|\xa|\leq (R_0-1)\ma$ for all $\alpha\in\nn$. It is easy to check that for $x\in \Omega\setminus B_{R_0\ma}(0)$, 
\begin{equation*}
\frac{L_\alpha G_\alpha}{G_\alpha}(x)=\left(\frac{\gamma'-\gamma}{|x|^2}+(a_\infty(x)-a_\alpha(x))+\nu\right)-\la\frac{\ua^{\crits-2}(x)}{|x|^s}.
\end{equation*}
It follows from \eqref{lim:a:alpha} that there exists $\alpha_0>0$ such that $a_\infty(x)-a_\alpha(x)\geq -\nu/2$ for all $\alpha>\alpha_0$ and all $x\in \Omega$. For a fixed  $\delta>0$,  \eqref{weak:est:2} yields $R>R_0$ such that for $\alpha>0$ large enough, we have that $\ua(x)\leq \delta|x|^{-(n-2)/2}$ for $x\in \Omega\setminus B_{R\ma}(0)$. Therefore, with \eqref{lim:la}, we get that for $x\in \Omega\setminus B_{R\ma}(0)$, 
\begin{equation*}
\frac{L_\alpha G_\alpha}{G_\alpha}(x)>\frac{1}{|x|^2}\left(\gamma'-\gamma-\mu_{\gamma,s,0}(\rn)\delta^{\crits-2s}+o(1)\right).
\end{equation*}
Up to taking $\delta>0$ small enough, we then get that $L_\alpha G_\alpha>0$ in $\Omega\setminus B_{R\ma}(0)$. It follows from \eqref{lim:va:1} and \eqref{est:G:low} that there exists $c(R)>0$ such that
\begin{equation*}
\ua(x)\leq c(R) \ma^{\frac{n-2}{2}}G_\alpha(x)\hbox{ for all }x\in \partial B_{R\ma}(0)\hbox{ and }\alpha\in\nn.
\end{equation*}
Therefore, defining $h_\alpha:=c(R) \ma^{\frac{n-2}{2}}G_\alpha-\ua$, we get that $L_\alpha h_\alpha>0$ in $\Omega\setminus B_{R\ma}(0)$ and $h_\alpha\geq 0$ in $\partial(\Omega\setminus B_{R\ma}(0))$. Since $G_\alpha>0$ in $\overline{\Omega\setminus B_{R\ma}(0)}$ and $L_\alpha G_\alpha>0$, it follows from the comparison principle of Berestycki-Nirenberg-Varadhan \cite{bnv} that $L_\alpha$ satisfies the comparison principle on $\Omega\setminus B_{R\ma}(\xa)$. Therefore,  $\ua\leq c(R) \ma^{\frac{n-2}{2}}G_\alpha$ in $\Omega\setminus B_{R\ma}(0)$. With the pointwise control \eqref{control:G}, we then get that
$$\ua(x)\leq C(R)\frac{\ma^{\frac{\beta_+(\gamma')-\beta_-(\gamma')}{2}}}{|x|^{\beta_+(\gamma')}}\hbox{ for all }x\in \Omega\setminus B_{R\ma}(0)
$$
Since this is valid for any $\gamma'>\gamma$ close to $\gamma$, with the remark made at the beginning of the proof, we get \eqref{est:eps}. 

\medskip\noindent
We now claim that there exists $C>0$ such that
\begin{equation}\label{est:c0}
\ua(x)\leq C\frac{\ma^{\frac{\bp-\bm}{2}}}{|x|^{\bp}}\hbox{ for all }x\in \Omega\setminus B_{\ma}(0).
\end{equation}
Indeed, as argued above, the result holds on $B_{R\ma}(0)\setminus B_{\ma}(0)$ for any $R>1$. In order to establish \eqref{est:c0}, we will prove it for any sequence $(z_\alpha)_\alpha\in\Omega$ such that 
\begin{equation}\label{lim:z}
\lim_{\alpha\to +\infty}\frac{|z_\alpha|}{\ma}=+\infty.
\end{equation}
Let $G_\alpha$ be the Green's function of $-\Delta-(\gamma|x|^{-2}+a_\alpha)$ on $\Omega$ with Dirichlet boundary condition. Green's representation formula in Appendix A, and the pointwise control \eqref{est:eps} yield
\begin{eqnarray*}
\ua(\za) &=& \int_\Omega G_\alpha(\za,x) \la\frac{\ua^{\crits-1}(x)}{|x|^s}\, dy\\
&\leq & C \int_\Omega\left(\frac{\max\{|\za|,|x|\}}{\min\{|\za|,|x|\}}\right)^{\bm}|x-\za|^{2-n}\frac{\ua^{\crits-1}(x)}{|x|^s}\, dx.\nonumber
\end{eqnarray*}
We split $\Omega$ into four subdomains. On $D_{1,\alpha}(R):=B_{R\ma}(0)$, we have from \eqref{lim:z}, \eqref{est:eps} and \eqref{lim:va:1} that
\begin{eqnarray*}
&&\left|\int_{D_{1,\alpha}}\left(\frac{\max\{|\za|,|x|\}}{\min\{|\za|,|x|\}}\right)^{\bm}|x-\za|^{2-n}\frac{\ua^{\crits-1}(x)}{|x|^s}\, dx\right|\\
&&\qquad \qquad \leq C |\za|^{-\bp}\int_{B_{\ma}(0)}\frac{\ua^{\crits-1}}{|x|^{s+\bms}}\, dx \leq C |\za|^{-\bp}\ma^{\frac{\bp-\bm}{2}}.\nonumber
\end{eqnarray*}
Let $D_{2,\alpha}(R):=\{R\ma<|x|< \frac{1}{2}|\za|\}$, and note that $|x-\za|>\frac{1}{2}|\za|$ for all $x\in D_{2,\alpha}(R)$. Therefore, taking $\eps>0$ sufficienty small in \eqref{est:eps}, we have that
\begin{eqnarray*}
&& \qquad \left|\int_{D_{2,\alpha}}\left(\frac{\max\{|\za|,|x|\}}{\min\{|\za|,|x|\}}\right)^{\bm}|x-\za|^{2-n}\frac{\ua^{\crits-1}(x)}{|x|^s}\, dx\right|\\
&&\qquad \leq C|\za|^{-\bp}\ma^{(\frac{\bp-\bm}{2}-\eps)(\crits-1)}\int_{D_{2,\alpha}}|x|^{-s-\bm-(\crits-1)(\bp-\eps)}\, dx\nonumber\\
&&\qquad \leq \theta(R)|\za|^{-\bp}\ma^{\frac{\bp-\bm}{2}},\nonumber
\end{eqnarray*}
as $\alpha\to +\infty$, where $\lim_{R\to +\infty}\theta(R)=0$.

\medskip\noindent Set $D_{3,\alpha}:=\{ \frac{1}{2}|\za|<|x|< 2|\za|\}$, and by using again \eqref{est:eps} with $\eps>0$ sufficiently small, we get that
\begin{eqnarray*}
&&\left|\int_{D_{3,\alpha}}\left(\frac{\max\{|\za|,|x|\}}{\min\{|\za|,|x|\}}\right)^{\bm}|x-\za|^{2-n}\frac{\ua^{\crits-1}(x)}{|x|^s}\, dx\right|\\
&&\leq C\ma^{(\frac{\bp-\bm}{2}-\eps)(\crits-1)}{|\za|^{-s-(\bp-\eps)(\crits-1)}}\int_{D_{3,\alpha}}|x-\za|^{2-n}\, dx\nonumber\\
&&\leq C\ma^{\frac{\bp-\bm}{2}}|\za|^{-\bp}\left(\frac{\ma}{|\za|}\right)^{\frac{\bp-\bm}{2}(\crits-2)-\eps(\crits-1)}.\nonumber
\end{eqnarray*}
Finally, let $D_{4,\alpha}:=\{|x|\geq 2|\za|\}\cap\Omega$. Since $|x-\za|\geq |x|/2$, then using  \eqref{est:eps} with $\eps>0$ sufficiently small, we get that
\begin{eqnarray*}
&&\left|\int_{D_{4,\alpha}}\left(\frac{\max\{|\za|,|x|\}}{\min\{|\za|,|x|\}}\right)^{\bm}|x-\za|^{2-n}\frac{\ua^{\crits-1}(x)}{|x|^s}\, dx\right|\\
&&\leq C|\za|^{-\bm}\ma^{(\frac{\bp-\bm}{2}-\eps)(\crits-1)}\int_{D_{4,\alpha}}|x|^{-s-\bp\crits +\eps(\crits-1)}\, dx\nonumber\\
&&\leq C\ma^{\frac{\bp-\bm}{2}}|\za|^{-\bp}\left(\frac{\ma}{|\za|}\right)^{\frac{\bp-\bm}{2}(\crits-2)-\eps(\crits-1)}. \nonumber
\end{eqnarray*}
Plugging together these estimates yields \eqref{est:c0}.

\medskip\noindent Since $U$ is a positive solution to \eqref{eq:U:bis} and $U\in\dundeuxrn$, it follows from the regularity Theorem \ref{th:hopf} that there exists $C_1>0$ such that $U(x)\simeq C_1|x|^{-\bms}$ as $s\to 0$. Taking the Kelvin transform $\tilde{U}(x):=|x|^{2-n}U(x|x|^{-2})$, we get that $\tilde{U}\in \dundeuxrn$ is also a positive solution to \eqref{eq:U:bis}, and enjoys a similar behavior at $0$. Transforming back yields the existence of $C_1,C_2>0$ such that
\begin{equation}\label{asymp:U}
U(x)\simeq_{x\to 0} \frac{C_1}{|x|^{\bm}}\hbox{ and }U(x)\simeq_{|x|\to \infty} \frac{C_1}{|x|^{\bp}}.
\end{equation}
We now show that there exists $H\in C^2(\overline{\Omega}\setminus\{0\})$ such that
\begin{equation}\label{lim:ua:H}
\lim_{\alpha\to +\infty}\frac{\ua}{\ma^{\frac{\bp-\bm}{2}}}=H\hbox{ in }C^2_{loc}(\overline{\Omega}\setminus\{0\}),
\end{equation}
and $H$ is a solution to 
\begin{equation}\label{lim:H}
\left\{\begin{array}{ll}
-\Delta H-\left(\frac{\gamma}{|x|^2}+a_\infty\right) H=0 & \hbox{ in }\Omega\setminus\{0\}\\
H>0 &\hbox{ in }\Omega\setminus\{0\}\\
H=0 &\hbox{ on }\partial\Omega.
\end{array}\right.
\end{equation}
\medskip\noindent Define $w_\alpha:=\ma^{-\frac{\bp-\bm}{2}}\ua$. Equation \eqref{eq:ua} then rewrites as
\begin{equation}\label{eq:va}
\left\{\begin{array}{ll}
-\Delta w_\alpha -\left(\frac{\gamma}{|x|^2}+\aa\right) w_\alpha = \la \ma^{(\crits-2)\frac{\bp-\bm}{2}}\frac{w_\alpha^{\crits-1}}{|x|^s}& \hbox{ in }\Omega\\
w_\alpha\geq 0& \hbox{ a.e. in }\Omega\\
w_\alpha=0& \hbox{ on }\partial \Omega,
\end{array}\right.
\end{equation}
and \eqref{est:c0} yields that $w_\alpha(x)\leq C |x|^{-\bp}$ for all $x\in\Omega\setminus\{0\}$ and $\alpha\in\nn$. It then follows from elliptic theory that there exists $H\in C^2(\Omegabar\setminus\{0\})$ such that $\lim_{\alpha\to +\infty}w_\alpha= H$ in $C^2_{loc}(\Omegabar\setminus\{0\})$. Passing to the limit in \eqref{eq:va} yields $H\geq 0$ and
$$-\Delta H-\left(\frac{\gamma}{|x|^2}+a_\infty\right)H=0\hbox{ in }\Omega\setminus\{0\}\hbox{ and }H=0\hbox{ on }\partial\Omega.$$
Fix $x\in\Omega\setminus\{0\}$. Green's representation formula, the positivity of $G_\alpha$ and a change of variable yields
\begin{eqnarray*}
\ua(x) &=& \int_\Omega G_\alpha(x,y) \la\frac{\ua^{\crits-1}(y)}{|y|^s}\, dy\\
&\geq& \int_{B_{2\ma}(0)\setminus B_{\ma}(0)} G_\alpha(x,y) \la\frac{\ua^{\crits-1}(y)}{|y|^s}\, dy\nonumber\\
&\geq & \ma^{\frac{n-2}{2}}\int_{B_{2}(0)\setminus B_{1}(0)} G_\alpha(x,\ma y) \la\frac{\va(y)^{\crits-1}(y)}{|y|^s}\, dy.\nonumber
\end{eqnarray*}
The asymptotics \eqref{asymp:G} in Appendix A yields $G_\alpha(x,z)\geq c_x|z|^{-\bm}$ for all $\alpha\in\nn$ and all $z\in B_{|x|/2}(0)$. Therefore, we get that for all $\alpha\in\nn$,  
$$\ua(x)\geq c_x \ma^{\frac{\bp-\bm}{2}}\int_{B_{2}(0)\setminus B_{1}(0)} |y|^{-\bm} \frac{\va(y)^{\crits-1}(y)}{|y|^s}\, dy$$
Passing to the limit as $\alpha\to +\infty$ and using \eqref{lim:va:1} yields $H(x)>0$, which proves our claim in \eqref{lim:H}.

\medskip\noindent 
Let now $\delta>0$ be such that $B_\delta(0)\subset \Omega$. For any $0<\eps<\delta$, the Pohozaev identity \eqref{gen:p:id} with $p=0$, and equation \eqref{eq:ua} yield
\begin{eqnarray}\label{id:poho:1}
&&-\int_{B_\delta(0)\setminus B_\eps(0)}\left(a_\alpha+\frac{x^i\partial_ia_\alpha}{2}\right)\ua^2\, dx = \int_{\partial (B_\delta(0)\setminus B_\eps(0))}B_\alpha(x)\, dx,
\end{eqnarray}
where
$$B_\alpha(x):=(x,\nu)\left(\frac{|\nabla \ua|^2}{2}-\left(\frac{\gamma}{|x|^2}+a_\alpha\right)\frac{\ua^2}{2}-\frac{\la\ua^{\crits}}{\crits |x|^s}\right)-\left(x^i\partial_i \ua+\frac{n-2}{2}\ua\right)\partial_\nu\ua.$$
Using the asymptotics \eqref{asymp:ua:0}, we pass to the limit as $\eps\to 0$ and get
\begin{eqnarray*}
&&-\int_{B_\delta(0)}\left(a_\alpha+\frac{x^i\partial_ia_\alpha}{2}\right)\ua^2\, dx = \int_{\partial B_\delta(0)}B_\alpha(x)\, dx.
\end{eqnarray*}
The limit \eqref{lim:ua:H} yields
\begin{eqnarray}\label{poho:H}
&&\qquad \lim_{\alpha\to +\infty}\ma^{-(\bp-\bm)}\int_{\partial B_\delta(0)}B_\alpha(x)\, d\sigma\\
&&=\int_{\partial B_\delta(0)}\left[(x,\nu)\left(\frac{|\nabla H|^2}{2}-\left(\frac{\gamma}{|x|^2}+a_\infty\right)\frac{H^2}{2}\right)-\left(x^i\partial_i H+\frac{n-2}{2}H\right)\partial_\nu H \right]\, d\sigma.\nonumber
\end{eqnarray}
\medskip\noindent 
Assuming now that $\bp-\bm>2$, we show that
\begin{equation}\label{1a=0}
a_\infty(0)=0.
\end{equation}
Indeed, note first that in this case, $\bp>\frac{n}{2}$. It follows from \eqref{est:c0} that
$$\lim_{R\to +\infty}\lim_{\alpha\to 0}\ma^{-2}\int_{B_\delta(0)\setminus B_{R\ma}(0)}\left(a_\alpha+\frac{x^i\partial_ia_\alpha}{2}\right)\ua^2\, dx =0.$$
With a change of variable, we get that
$$\int_{ B_{R\ma}(0)}\left(a_\alpha+\frac{x^i\partial_ia_\alpha}{2}\right)\ua^2\, dx =\ma^2\int_{ B_{R}(0)}\left(a_\alpha+\frac{x^i\partial_ia_\alpha}{2}\right)(\ma x)\va^2\, dx.$$
The limit \eqref{lim:va:1} and the compactness of the embedding $H^1\hookrightarrow L^2$ yields the convergence of $\va$ to $U$ in $L^2_{loc}(\rn)$. It follows from \eqref{asymp:U} that $U\in L^2(\rn)$, the two preceding identities therefore yield
$$\lim_{\alpha\to 0}\ma^{-2}\int_{B_\delta(0)}\left(a_\alpha+\frac{x^i\partial_ia_\alpha}{2}\right)\ua^2\, dx =a_\infty(0)\int_{\rn}U^2\, dx.$$
Plugging this limit in the Pohozaev identity \eqref{id:poho:1} and using the limit above yields that $a_\infty(0)=O(\ma^{\bp-\bm-2})+o(1)$ as $\alpha\to +\infty$, and therefore $a_\infty(0)=0$. 

\medskip\noindent 
We now assume that $\bp-\bm=2$, and we show again that
\begin{equation}\label{2a=0}
a_\infty(0)=0.
\end{equation}
Indeed, assume that $a_\infty(0)\neq 0$. Without loss of generality, we can suppose that $a_\infty(0)>0$. Up to taking $\delta>0$ smaller and $\alpha$ large, we have that  $a_\alpha(x)+\frac{x^i\partial_ia_\alpha(x)}{2}\geq \frac{a_\infty(0)}{2}$ for $x\in B_\delta(0)$. It then follows from \eqref{id:poho:1} and \eqref{poho:H} that there exists $C>0$ such that
$$\int_{B_{R\ma}(0)}\ua^2\, dx\leq \int_{B_\delta(0)}\ua^2\, dx\leq C\ma^2\hbox{ for all }\alpha\in\nn.$$
With a change of variable, the limit \eqref{lim:va:1}, letting $\alpha\to +\infty$ and then $R\to +\infty$, we get that $U\in L^2(\rn)$, which is impossible due to \eqref{asymp:U} and  $2\bp=n$. Therefore $a_\infty(0)=0$.

\medskip\noindent
Finally, we show that if $\bp-\bm<2$, then 
\begin{equation}\label{masszero}
m_{\gamma,a_\infty}(\Omega)=0,
\end{equation}
where $m_{\gamma,a_\infty}(\Omega)$ is the Hardy-singular mass as defined in Proposition \ref{prop:def:mass}. \\ 
Indeed, since $2\bp<n$, we have that
\begin{eqnarray*}
\qquad \int_{B_\delta(0)}\left(a_\alpha+\frac{x^i\partial_ia_\alpha}{2}\right)\ua^2\, dx&=&O\left(\int_{B_\delta(0)}\ma^{\bp-\bm}|x|^{-2\bp}\, dx\right)\\
&=&O\left(\ma^{\bp-\bm}\delta^{n-2\bp}\right), \nonumber
\end{eqnarray*}
uniformly with respect to $\alpha$ and $\delta>0$. Combining with \eqref{poho:H}, we get that
\begin{equation}\label{lim:H:0}
\lim_{\delta\to 0}\int_{\partial B_\delta(0)}\left[(x,\nu)\left(\frac{|\nabla H|^2}{2}-\left(\frac{\gamma}{|x|^2}+a_\infty\right)\frac{H^2}{2}\right)-\left(x^i\partial_i H+\frac{n-2}{2}H\right)\partial_\nu H \right]\, d\sigma=0.
\end{equation}
Since $\bp-\bm<2$, it follows from the definition of the mass that there exists $c>0$ such that
$$H(x)=c\left(\frac{1}{|x|^{\bp}}+\frac{m_{\gamma,a_\infty}(\Omega)}{|x|^{\bm}}+o\left(\frac{1}{|x|^{\bm}}\right)\right)\hbox{ as }x\to 0.$$
Since $H$ solve the equation \eqref{eq:H}, standard elliptic theory yields that this estimate can be differentiated. Therefore, putting it into \eqref{lim:H:0} yields $m_{\gamma,a_\infty}(\Omega)=0$. 

\medskip\noindent Theorem \ref{th:blowup:1} is a consequence of  \eqref{1a=0}, \eqref{2a=0}, and \eqref{masszero}.
\qed

\section{\, Blow-Up analysis in the merely singular case}
In this section, we perform the blow-up analysis in the merely singular case, that is when 
$$s=0\hbox{ and }\gamma<0.$$
We let again $(a_\alpha)_{\alpha\in\nn}\in C^{1}(\Omegabar)$, $a_\infty\in C^{1}(\Omegabar)$, $(\la)_\alpha\in (0,+\infty)$ such that \eqref{lim:a:alpha} and \eqref{lim:la} hold.
We let $(u_\alpha)_\alpha\in \huno$ be a sequence of weak solutions to \eqref{eq:ua} such that \eqref{mass:ua} holds. In this case, \eqref{eq:ua} and \eqref{mass:ua} rewrite as:
\begin{equation}\label{eq:ua:bis}
\left\{\begin{array}{ll}
-\Delta \ua +\left(\frac{|\gamma|}{|x|^2}-\aa\right) \ua = \la \ua^{\crit-1}& \hbox{ in }\Omega\\
\ua\geq 0& \hbox{ a.e. in }\Omega\\
\ua=0& \hbox{ on }\partial \Omega
\end{array}\right.
\end{equation}
and
\begin{equation}\label{mass:ua:bis}
\int_\Omega \ua^{\crit}\, dx=1.
\end{equation}
We suppose that
\begin{equation}\label{lim:weak:bis}
\ua\rightharpoonup 0\hbox{ as }\alpha\to +\infty\hbox{ weakly in }\huno.
\end{equation}
We let $\Omega'$ be a smooth bounded domain of $\rn$ such that $\Omega\subset\subset \Omega'$ is relatively compact in $\Omega'$. We extend $(\aa)_\alpha$ and $a_\infty$ on $\Omega'$ such that \eqref{lim:a:alpha} holds on $\Omega'$ and that the operator $-\Delta-(\gamma|x|^{-2}+a_\infty)$ is coercive on $\Omega'$. This assumption is equivalent to saying that there exists $c>0$ such that for $\alpha\in\nn$ large enough, we have
\begin{equation}\label{hyp:coer:2}
\lambda_1(-\Delta-(\gamma|x|^{-2}+a_\alpha))=\inf_{\varphi\in \huno\setminus\{0\}}\frac{\int_\Omega\left(|\nabla\varphi|^2-\left(\frac{\gamma}{|x|^2}+a_\alpha\right)\varphi^2\right)\, dx}{\int_\Omega \varphi^2\, dx}\geq c>0.
\end{equation}

This section is devoted to the proof of the following result:
\begin{theorem}\label{th:blowup:2} Let $\Omega$ be a smooth bounded domain of $\rn$, $n\geq 3$, such that $0\in \Omega$ is an interior point. Fix $\gamma<0$ and  let $(\aa)_\alpha\in C^1(\Omegabar)$, $(\la)_\alpha\in (0,+\infty)$ and $(\ua)_\alpha\in \huno$ be such that \eqref{lim:a:alpha}, \eqref{lim:la}, \eqref{eq:ua:bis} and \eqref{mass:ua:bis} hold. We let $(\xa)_\alpha\in \Omega$ and $(\ma)_\alpha\in (0,+\infty)$ be such that $\ua(\xa):=\sup_{\Omega}\ua=\ma^{-\frac{n-2}{2}}$.\\ Then $\lim_{\alpha\to +\infty}\xa=x_0\in \Omegabar$, $\lim_{\alpha\to +\infty}\ma=0$ and
\begin{enumerate}[i)]
\item If $n\geq 4$, then $x_0\neq 0$ and $a_\infty(x_0)+\frac{\gamma}{|x_0|^2}=0$;
\item If $n=3$, then $x_0\in\Omega\setminus\{0\}$ and $R_{\gamma,a_\infty}(\Omega,x_0)=0$ (see \eqref{def:R} for the definition).
\end{enumerate}
In addition, there exists $C>0$ such that
$$\ua(x)\leq C\left(\frac{\ma}{\ma^2+|x-\xa|^2}\right)^{\frac{n-2}{2}}\hbox{ for all }x\in\Omega\hbox{ and }\alpha\in\nn.$$
\end{theorem}
Before delving into the proof, it is important to note a few observations that are relevant for the case $s=0$ and $\gamma<0$. First note that in this case $\bm<0$, and therefore, it follows from \eqref{asymp:ua:0} that for any $\alpha\in\nn$, $\ua$ can be extended continuously at $0$ by $0$, which means that we can and will consider $\ua\in C^0(\Omegabar)$. In the definition \eqref{def:xa}, we shall take $\tau:=0$ and therefore, the sequence $(\xa)_\alpha\in\Omega$ will be such that
\begin{equation}\label{def:xa:2}
\ua(\xa):=\sup_{x\in\Omega}\ua(x)\hbox{ and }\ma:=\ua(\xa)^{-2/(n-2)}.
\end{equation}
It then follows from Proposition \ref{prop:useful} that
\begin{equation}\label{ppty:xa}
\lim_{\alpha\to+\infty}\ma=0\hbox{ and }\lim_{\alpha\to +\infty}\frac{d(\xa,\partial\Omega)}{\ma}=+\infty.
\end{equation}
Another remark is that  \eqref{eq:ua:bis} implies
\begin{equation*}
-\Delta \ua -\aa \ua \leq \la \ua^{\crit-1} \hbox{ in }\Omega,
\end{equation*}
which means that $\ua$ is a subsolution of a nonlinear elliptic inequation with no Hardy potential term. We shall then be able to perform a blow-up analysis in the spirit of Druet-Hebey-Robert \cite{dhr} to obtain a pointwise control of $\ua$ by a standard bubble. The conclusion of Theorem \ref{th:blowup:2} will then follow from classical arguments via the Pohozaev identity and the analysis on the boundary in the spirit of Druet \cite{d2}.

\medskip\noindent Set
$$\va(x):=\ma^{\frac{n-2}{2}}\ua(\xa+\ma x)\hbox{ for all }x\in \ma^{-1}(\Omega-\xa).$$
Equation \eqref{eq:ua:bis} and \eqref{mass:ua:bis} then rewrites
\begin{equation}\label{eq:v:2}
-\Delta \va +\left(\frac{|\gamma|}{|\frac{\xa}{\ma}+x|^2}-\ma^2\aa(\xa+\ma x)\right) \va = \la \va^{\crit-1} \hbox{ in }\frac{\Omega-\xa}{\ma}
\end{equation}
\begin{equation}\label{eq:mass:va}
\int_{\frac{\Omega-\xa}{\ma}}\va^{\crit}\, dx=1.
\end{equation}
\medskip\noindent We first claim that 
\begin{equation}\label{result:step1prime}
\lim_{\alpha\to +\infty}\frac{|\xa|}{\ma}=+\infty,
\end{equation}
and
\begin{equation}\label{result:step1}
\lim_{\alpha\to +\infty}\va=v:=\left(\frac{1}{1+\frac{|\cdot|^{2}}{n(n-2)K(n,2)^2}}\right)^{\frac{n-2}{2}}\hbox{ in }C^2_{loc}(\rn)\hbox{ with }\int_{\rn}v^{\crit}\, dx=1.
\end{equation}
Indeed, it follows from \eqref{ppty:xa} that for any $R>0$, there exists $\alpha_0>0$ such that $B_R(0)\subset \frac{\Omega-\xa}{\ma}$ for all $\alpha>\alpha_0$. Since $(\ua)_\alpha$ is uniformly bounded in $\huno$, then $(\va)_\alpha$ is bounded in $H^1_{loc}(\rn)$. Up to extracting a subsequence, there exists $v\in H_{loc}^1(\rn)$ such that $\va\rightharpoonup v$ as $\alpha\to +\infty$ weakly in $H^1_{loc}(\rn)$ and strongly in $L^2_{loc}(\rn)$. Since 
$$-\Delta \va -\ma^2\aa(\xa+\ma x) \va \leq  \la \va^{\crit-1} \hbox{ in }B_R(0),$$
and $0\leq \va\leq 1$, it follows from DeGiorgi-Nash-Moser iterative scheme (see for instance Theorem 4.1 in Han-Lin \cite{hl}), that there exists $C>0$ such that for all $\alpha> \alpha_0$,
$$1=|\va(0)|\leq C \Vert\va\Vert_{L^2(B_1(0))}$$
and therefore, passing to the strong limit in $L^2$, we get that $1\leq C\Vert v\Vert_{L^2(B_1(0))}$, and hence $v\not\equiv 0$. Since $0<\va\leq \va(0)=1$, equation \eqref{eq:v:2} and elliptic theory yields $v\in C^2(\rn\setminus \{\theta_\infty\})$ and $\va\to v$ in $C^2_{loc}(\rn\setminus\{\theta_\infty\})$ with
\begin{equation}\label{eq:v}
-\Delta v +\frac{|\gamma|}{|x-\theta_\infty|^2}v = \mu_{\gamma,s,0}(\rn) v^{\crit-1} \hbox{ in }\rn\setminus\{\theta_\infty\}
\end{equation}
where $\theta_\infty:=-\lim_{\alpha\to +\infty}\ma^{-1}\xa$ if this limit is finite. Otherwise $\theta_\infty:=\infty$, in which case $\rn\setminus \{\theta_\infty\}:=\rn$. In addition, passing to the weak limit in \eqref{eq:mass:va} yields
\begin{equation*}
\int_{\rn}v^{\crit}\, dx=\lim_{R\to +\infty}\int_{B_R(0)}v^{\crit}\, dx\leq \lim_{R\to +\infty}\lim_{\alpha\to +\infty}\int_{B_R(0)}\va^{\crit}\, dx\leq 1. 
\end{equation*}
Since $\int_{B_R(0)}|\nabla\va|^2\, dx=\int_{B_{R\ma}(\xa)}|\nabla\ua|^2\, dx\leq C$ uniformy for all $R>0$ and $\alpha>0$ large enough, passing to the weak limit yields $|\nabla v|\in L^2(\rn)$. Since $v\in L^{\crit}(\rn)$, classical arguments yield that $v\in \dundeuxrn$. Multiplying \eqref{eq:v} by $v$ and integrating, we obtain
$$\int_{\rn}|\nabla v|^2\, dx\leq \int_{\rn}\left(|\nabla v|^2 +\frac{|\gamma|}{|x-\theta_\infty|^2}v^2 \right)\, dx=  \mu_{\gamma,s,0}(\rn) \int_{\rn} v^{\crit}\, dx.$$
Since $v\not\equiv 0$, the Sobolev inequality yields
$$\frac{\int_{\rn}|\nabla v|^2\, dx}{\left(\int_{\rn} v^{\crit}\, dx\right)^{\frac{2}{\crit}}}\geq \mu_{0,0}(\rn)=\mu_{\gamma,0}(\rn).$$
Since $\int_{\rn}v^{\crit}\, dx\leq 1$ and $|\gamma|>0$, putting these latest inequalities together yields 
$$\theta_\infty=\infty\hbox{ and }\int_{\rn}v^{\crit}\, dx=1.$$
We then get
\begin{equation}\label{lim:va}
\lim_{\alpha\to+\infty}\frac{|\xa|}{\ma}=+\infty \hbox{ and }\lim_{\alpha\to +\infty}\va=v\hbox{ in }C^2_{loc}(\rn)
\end{equation}
where $v\in \dundeuxrn\cap C^2(\rn)$ is such that
$$-\Delta v=\mu_{0,0}(\rn)v^{\crit-1}\hbox{ in }\rn\; ;\; \int_{\rn}v^{\crit}\, dx=1\; ;\; 0\leq v\leq v(0)=1.$$
Then \eqref{result:step1prime} and \eqref{result:step1}  follow from \eqref{lim:va}, this latest assertion and the classification of Caffarelli-Gidas-Spruck \cite{cgs}. 
 
\medskip\noindent
We now claim that there exists $C>0$ such that 
\begin{equation}\label{est:ua:point}
\ua(x)\leq\left(\frac{\ma}{\ma^2+|x-\xa|^2}\right)^{\frac{n-2}{2}}\hbox{ for all }\alpha\in\nn\hbox{ and }x\in \Omega.
\end{equation}
This estimate is by now standard and is in the spirit of similar results obtained by several authors. See for instance Druet-Hebey-Robert \cite{dhr} and the several references therein. When possible, we shall only sketch an outline to the proof and we refer to these references for details.

Note first that 
\begin{equation}\label{2.2} \lim_{R\to +\infty}\lim_{\alpha\to +\infty}\int_{\Omega\setminus B_{R\ma}(\xa)}\ua^{\crit}\, dx=0.
\end{equation}
Indeed, the convergence of $(\va)_\alpha$ to $v$ in \eqref{result:step1} yields that asymptotically, $B_{R\ma}(\xa)$ exhausts almost all the energy in \eqref{mass:ua:bis}.  

Next, we claim that there exists $C>0$ such that
\begin{equation*}
|x-\xa|^{\frac{n-2}{2}}\ua(x)\leq C\hbox{ for all }\alpha\in\nn\hbox{ and }x\in\Omega.
\end{equation*}
Indeed, if not we find $(\ya)_\alpha\in\Omega$ that achieve the supremum of the left-hand-side and which go to $+\infty$ as $\alpha\to +\infty$. The same blow-up procedure as above at $\ya$ yields that asymptotically, $B_{\ua(\ya)^{-2/(n-2)}}(\ya)$ carries a nonzero mass of $\ua^{\crit}\, dx$, contradicting \eqref{2.2}, since this ball is disjoint from $B_{R\ma}(\xa)$ for $R$ and $\alpha$ large. \\ 
A similar argument --that we omit-- yields that 
\begin{equation}\label{est:100}
\lim_{R\to +\infty}\lim_{\alpha\to +\infty}\sup_{x\in \Omega\setminus B_{R\ma}(\xa)}|x-\xa|^{\frac{n-2}{2}}\ua(x)=0.
\end{equation}
 Let now  $\eta_0\in C^\infty(\rr)$ be such that $0\leq\eta_0\leq 1$, $\eta_0(t)=0$ if $t\leq 1$ and $\eta_0(t)=1$ if $t\geq 2$. We define $\eta_\eps(x):=\eta_0(|x|/\eps)$ for $x\in\rn$. We claim that there exists $\epsilon>0$ such that 
\begin{equation}\label{coer:op}
-\Delta-\eta_\eps(x)\gamma|x|^{-2}-a_\infty-c/2\hbox{ is coercive.}
\end{equation}
To prove this claim, we shall need the following continuity lemma for the first eigenvalue. Recall that for any $V:\Omega\to\rr$ measurable such that for some $C>0$, we have $|x|^2|V(x)|\leq C$ for a.e. $x\in\Omega$, the following ratio
$$\lambda_1(-\Delta+V):=\inf_{\varphi\in \huno\setminus\{0\}}\frac{\int_\Omega (|\nabla \varphi|^2+V\varphi^2)\, dx}{\int_\Omega \varphi^2\, dx}$$
 is well defined and is finite.

\begin{lem}\label{lem:vp} Let $\Omega\subset\rn$, $n\geq 3$, be a smooth bounded domain. Let $(V_k)_k:\Omega\to \rr$ and $V_\infty:\Omega\to\rr$ be measurable functions and let $(x_k)_k\in\Omega$ be a sequence of points. We assume that
\begin{enumerate}[i)]
\item $\lim_{k\to +\infty}V_k(x)=V_\infty(x)$ for a.e. $x\in\Omega$,
\item There exists $C>0$ such that $|V_k(x)|\leq C|x-x_k|^{-2}$ for all $k\in\nn$ and $x\in\Omega$.
\item $\lim_{k\to +\infty}x_k=0\in\Omega.$
\item For some $\gamma_0<(n-2)^2/4$, there exists $\delta>0$ such that $|V_k(x)|\leq \gamma_0|x-x_k|^{-2}$ for all $k\in\nn$ and $x\in B_\delta(0)\subset\Omega$.  
\item The first eigenvalue $\lambda_1(-\Delta+V_k)$ is achieved for all $k\in\nn$.
\end{enumerate}
Then
\begin{equation*}
\lim_{k\to +\infty}\lambda_1(-\Delta+V_k)=\lambda_1(-\Delta+V_\infty).
\end{equation*}
\end{lem}
\smallskip\noindent{\it Proof of Lemma \ref{lem:vp}:} We first claim that $(\lambda_1(-\Delta+V_k))_k$ is bounded. Indeed, fix $\varphi\in \huno\setminus\{0\}$ and use the Hardy inequality to write for all $k\in\nn$, 
$$\lambda_1(-\Delta+V_k)\leq \frac{\int_\Omega (|\nabla \varphi|^2+V_k\varphi^2)\, dx}{\int_\Omega \varphi^2\, dx}\leq \frac{\int_\Omega (|\nabla \varphi|^2+C|x-x_k|^{-2}\varphi^2)\, dx}{\int_\Omega \varphi^2\, dx}:=M<+\infty$$
For the lower bound, we have for any $\varphi\in \huno$,
\begin{eqnarray}
\int_\Omega (|\nabla \varphi|^2+V_k\varphi^2)\, dx&= & \int_\Omega |\nabla \varphi|^2\, dx+\int_{B_\delta(0)}V_k\varphi^2\, dx+\int_{\Omega\setminus B_\delta(0)}V_k\varphi^2\, dx\nonumber \\
&\geq & \int_\Omega |\nabla \varphi|^2\, dx-\gamma_0\int_{B_\delta(0)}|x-x_k|^{-2}\varphi^2\, dx\nonumber\\
&&-4C\delta^{-2}\int_{\Omega\setminus B_\delta(0)}\varphi^2\, dx\nonumber \\
&\geq & \left(1-4\gamma_0/(n-2)^2\right)\int_\Omega |\nabla \varphi|^2\, dx-4C\delta^{-2}\int_{\Omega}\varphi^2\, dx.\label{ineq:123}
\end{eqnarray}
Since $\gamma_0<(n-2)^2/4$, we then get that $\lambda_1(-\Delta+V_k)\geq -4C\delta^{-2}$ for large $k$, which proves the lower bound.  

\smallskip\noindent Up to a subsequence, we can now assume that $(\lambda_1(-\Delta+V_k))_k$ converges as $k\to +\infty$. We now show that 
\begin{equation}\label{ineq:vp:100}
\liminf_{k\to +\infty}\lambda_1(-\Delta+V_k)\geq \lambda_1(-\Delta+V_\infty).
\end{equation}
For $k\in\nn$, we let $\varphi_k\in \huno$ be a minimizer of $\lambda_1(-\Delta+V_k)$ such that $\int_\Omega\varphi_k^2\, dx=1$. In particular,
\begin{equation}\label{eq:phi:vp}
-\Delta\varphi_k+V_k\varphi_k=\lambda_1(-\Delta+V_k)\varphi_k\hbox{ weakly in }\huno.
\end{equation}
Inequality \eqref{ineq:123} above yields the boundedness of $(\varphi_k)_k$ in $\huno$. Up to a subsequence, we let $\varphi\in \huno$ such that, as $k\to +\infty$, $\varphi_k\rightharpoonup \varphi$ weakly in $\huno$, $\varphi_k\to \varphi$ strongly in $L^2(\Omega)$ (then $\int_\Omega\varphi^2\, dx=1$) and $\varphi_k(x)\to\varphi(x)$ for a.e. $x\in \Omega$. Letting
 $k\to +\infty$ in \eqref{eq:phi:vp}, the hypothesis on $(V_k)$ allow us to conclude that 
$$-\Delta\varphi+V_\infty\varphi_k=\lim_{k\to +\infty}\lambda_1(-\Delta+V_k)\varphi\hbox{ weakly in }\huno.$$
Since $\int_\Omega\varphi^2\, dx=1$ and we have extracted subsequences, we then get \eqref{ineq:vp:100}.

\medskip\noindent Finally, we prove the reverse inequality. For $\eps>0$, let $\varphi\in \huno$ be such that 
$$ \frac{\int_\Omega (|\nabla \varphi|^2+V_\infty\varphi^2)\, dx}{\int_\Omega \varphi^2\, dx}\leq \lambda_1(-\Delta+a_\infty)+\eps.$$
We have 
$$\lambda_1(-\Delta+V_k)\leq \lambda_1(-\Delta+V_\infty)+\eps +\frac{\int_\Omega |V_k-V_\infty|\varphi^2\, dx}{\int_\Omega \varphi^2\, dx}.$$
The hypothesis of Lemma \ref{lem:vp} allow us to conclude that $\int_\Omega |V_k-V_\infty|\varphi^2\, dx\to 0$ as $k\to +\infty$. Therefore $\limsup_{k\to +\infty}\lambda_1(-\Delta+V_k)\leq \lambda_1(-\Delta+V_\infty)+\eps$ for all $\eps>0$. Letting $\eps\to 0$, we get the reverse inequality and 
 the conclusion of Lemma \ref{lem:vp}.\qed

\medskip\noindent We now prove \eqref{coer:op}. First note that the coercivity property \eqref{hyp:coer:2} yield
$$\lambda_1(-\Delta-\gamma|x|^{-2}-a_\infty-c/2)\geq c/2.$$
Define $V_\epsilon(x):=-\eta_\epsilon(x)\gamma|x|^{-2}-a_\infty-c/2$ for all $x\in\Omega$ and $\epsilon>0$. Since $V_\epsilon\in C^0(\overline{\Omega})$, the eigenvalue $\lambda_1(-\Delta+V_\epsilon)$ is achieved. It then follows from Lemma \ref{lem:vp} that $\lambda_1(-\Delta+V_\epsilon)\to \lambda_1(-\Delta-\gamma|x|^{-2}-a_\infty-c/2)\geq c/2$ as $\epsilon\to 0$. Therefore, there exists $\epsilon>0$ such that $\lambda_1(-\Delta-\eta_\epsilon(x)\gamma|x|^{-2}-a_\infty-c/2)>0$. This proves \eqref{coer:op}.

\medskip\noindent Fix now $\nu\in (0,1)$. We claim that there exists $C_\nu,R_\nu>0$ such that 
\begin{equation}\label{est:nu}
\ua(x)\leq C_\nu \ma^{\frac{n-2}{2}-\nu(n-2)}|x-\xa|^{-(n-2)(1-\nu)}\hbox{ for all }x\in\Omega\setminus B_{R_\nu\ma}(\xa).
\end{equation}
Since the proof is similar to Step 6.3 (p1228) in Ghoussoub-Robert \cite{gr1}, we just give the main points and leave the details to the reader. 
We let $G_\epsilon$ be the Green's function of the operator $-\Delta-\eta_\epsilon(x)\gamma|x|^{-2}-a_\infty-c/2$ with Dirichlet boundary condition on $\Omega'$. Since $\xa\in \Omega\subset\subset \Omega'$ for all $\alpha\in\nn$, it follows from classical properties of the Green's function (see for instance \cite{r1}) that there exists $c_1>0$ and $\delta>0$ such that
\begin{equation}\label{low:G}
\frac{|\nabla G_\epsilon(x,\xa)|}{G_\epsilon(x,\xa)}\geq \frac{c_1}{|x-\xa|}\hbox{ for all }x\in B_\delta(\xa)\subset \Omega'.
\end{equation}
and
\begin{equation}\label{low:G:2}
c_1^{-1}\geq |x-\xa|^{n-2}G_\epsilon(x,\xa)\geq c_1\hbox{ for all }x\in B_\delta(\xa)\subset \Omega'.
\end{equation}
Consider the operator $L_\alpha:=-\Delta-(\gamma|x|^{-2}+a_\alpha)-\la \ua^{\crit-2}$. Straightforward computations yield
\begin{eqnarray*}
\frac{L_\alpha G_\epsilon^{1-\nu}(\xa,\cdot)}{ G_\epsilon^{1-\nu}(\xa,\cdot)}
&=&\nu(1-\nu)\frac{|\nabla G_\epsilon(x,\xa)|^2}{G_\epsilon(x,\xa)^2}+\frac{|\gamma|(1-\eta_\epsilon)+\nu|\gamma|\eta_\epsilon}{|x|^2}\\
&&+\frac{c(1-\nu)}{2}+a_\infty-a_\alpha-\nu a_\infty-\la \ua^{\crit-2}\nonumber \\
&\geq& \nu(1-\nu)\frac{|\nabla G_\epsilon(x,\xa)|^2}{G_\epsilon(x,\xa)^2}+\frac{c(1-\nu)}{2}\nonumber\\
&&-\Vert a_\alpha-a_\infty\Vert_\infty-\nu \Vert a_\infty\Vert_\infty-\la \ua^{\crit-2}.\nonumber
\end{eqnarray*}
Writing $\Omega\setminus B_{R\ma}(\xa)$ as a subset of the union of $B_\delta(\xa)\setminus B_{R\ma}(\xa)$ and $\Omega\setminus B_\delta(\xa)$, and using \eqref{est:100} and \eqref{low:G}, we get that there exists $R_\nu>0$ such that $L_\alpha G_\epsilon^{1-\nu}(\xa,\cdot)>0$ in $\Omega\setminus B_{R_\nu\ma}(\xa)$. It follows from the convergence \eqref{result:step1} and \eqref{low:G:2} that there exists $C_\nu>0$ such that $\ua\leq C_\nu G_\epsilon^{1-\nu}(\xa,\cdot)\ma^{\frac{n-2}{2}-\nu(n-2)}$ on $\partial (\Omega\setminus B_{R_\nu\ma}(\xa))$. Since $L_\alpha\ua=0$, it then follows from the comparison principle of Beresticky-Nirenberg-Varadhan \cite{bnv} (with an extra care for the singular point $0$) that $\ua\leq C_\nu G_\epsilon^{1-\nu}(\xa,\cdot)\ma^{\frac{n-2}{2}-\nu(n-2)}$ on $\Omega\setminus B_{R_\nu\ma}(\xa)$. This combined with \eqref{low:G:2}, yield \eqref{est:nu}.  

\medskip\noindent 
We now prove the pointwise control claimed in \eqref{est:ua:point}.

\smallskip\noindent As a preliminary remark, we note that \eqref{est:nu} and the convergence \eqref{result:step1} yield that for any $\nu\in (0,1)$, there exists $C_\nu>0$ such that
\begin{equation}\label{eq:110}
\ua(x)\leq C_\nu \frac{\ma^{\frac{n-2}{2}-\nu(n-2)}}{\left(\ma+|x-\xa|\right)^{n-2-\nu(n-2)}}\hbox{ for all }x\in\Omega.
\end{equation}
Letting $G_\alpha$ be the Green's function for $-\Delta-(\gamma|x|^{-2}+a_\alpha)$ with Dirichlet boundary condition on $\Omega$, we get from \eqref{est:G:up} that there exists $C>0$ such that
\begin{equation*}
0<G_\alpha(x,y)\leq C \left(\frac{\min\{|x|,|y|\}}{\max\{|x|,|y|\}}\right)^{|\bm|}|x-y|^{2-n}\leq C|x-y|^{2-n}
\end{equation*}
for all $x,y\in\Omega$, $x\neq y$. Here, we have used that $\bm<0$ since $\gamma<0$. Green's representation formula in Appendix A and \eqref{eq:110} then yield
\begin{eqnarray*}
\ua(x)&=& \int_\Omega G_\alpha(x,y)\la \ua^{\crit-1}(y)\, dy\\
&\leq &C \int_\Omega |x-y|^{2-n}\left(\frac{\ma^{\frac{n-2}{2}-\nu(n-2)}}{\left(\ma+|x-\xa|\right)^{n-2-\nu(n-2)}}\right)^{\crit-1}\, dy.\nonumber
\end{eqnarray*}
By estimating this integral as in \cite{gr1}, one gets the pointwise control \eqref{est:ua:point}. 

\medskip\noindent 
Now assume that $\lim_{\alpha\to +\infty}\xa=0$, and let $\ra:=|\xa|$. We claim that there exists $\chi>0$ such that
\begin{equation}\label{step3one}
\lim_{\alpha\to 0}\frac{\ra^{n-2}}{\ma^{\frac{n-2}{2}}}\ua(\ra x)= \chi G_{\theta_\infty}(x) \quad \hbox{ for all }x\in \rn\setminus\{0,\theta_\infty\},
\end{equation}
where $\theta_\infty:=\lim_{\alpha\to +\infty}\frac{\xa}{|\xa|}$. Moreover, this convergence holds in $C^2_{loc}(\rn\setminus\{0,\theta_\infty\})$. Here, $G_{\theta_\infty}$ satisfies properties (i) to (iv) of Theorem \ref{th:green:gamma}.

Indeed, our assumptions and \eqref{result:step1prime} yield that
$$\lim_{\alpha\to +\infty}\ra=0\hbox{ and }\lim_{\alpha\to+\infty}\frac{\ra}{\ma}=+\infty.$$
Define for $\alpha\in\nn$, the functions
\begin{equation}\label{def:tua}
\tua(x):=\frac{\ra^{n-2}}{\ma^{\frac{n-2}{2}}}\ua(\ra x) \quad \hbox{for $x\in \frac{\Omega}{\ra}.$}
\end{equation}
It follows from \eqref{eq:ua:bis} and the pointwise control \eqref{est:ua:point} that
$$-\Delta \tua +\left(\frac{|\gamma|}{|x|^2}-\ra^2\aa(\ra x)\right) \tua = \la \left(\frac{\ma}{\ra}\right)^2\tua^{\crit-1} \quad \hbox{ in }\frac{\Omega}{\ma},$$
and
\begin{equation}\label{def:thetaa}
0<\tua(x)\leq C|x-\theta_\alpha|^{2-n}\hbox{ with }\theta_\alpha:=\frac{\xa}{|\xa|}.
\end{equation}
Elliptic theory yields the existence of $\tilde{u}\in C^2(\rn\setminus\{0, \theta_\infty\})$ such that
\begin{eqnarray}
&&\tua\to \tilde{u}\hbox{ in }C^2_{loc}(\rn\setminus\{0, \theta_\infty\})\label{lim:tua}\\
&&-\Delta \tilde{u} +\frac{|\gamma|}{|x|^2}\tilde{u} = 0\hbox{ in }\rn\setminus\{0, \theta_\infty\}\label{eq:tu:2}\\
&&0\leq \tilde{u}(x)\leq C|x-\theta_\infty|^{2-n}.\label{bnd:tu}
\end{eqnarray}
We are now aiming for a more precise control of $\tilde{u}$. For that, we consider $G_\alpha$,  the Green's function for $-\Delta -(\gamma|x|^{-2}+\aa)$ in $\Omega$ with Dirichlet boundary condition. For $x\in\rnp$, we have for all $\alpha\in\nn$, 
\begin{equation}\label{green:ra}
\ua(\ra x)=\int_\Omega G_\alpha(\ra x,y)\la \ua^{\crit-1}(y)\, dy.
\end{equation}
Since $\bm<0$, the estimate \eqref{est:G:up} yields
\begin{eqnarray*}
\tua(x)&\leq& C\frac{\ra^{n-2}}{\ma^{\frac{n-2}{2}}}\int_{B_R(0)}  \left(\frac{\min\{|\ra x|,|y|\}}{\max\{|\ra x|,|y|\}}\right)^{|\bm|}|\ra x-y|^{2-n}\left(\frac{\ma}{\ma^2+|y-\xa|^2}\right)^{\frac{n+2}{2}}\, dy\nonumber\\
 &\leq& C\frac{\ra^{n-2}}{\ma^{\frac{n-2}{2}}}\int_{D_\alpha^1}F_\alpha(x,y)\, dy+ C\frac{\ra^{n-2}}{\ma^{\frac{n-2}{2}}}\int_{D_\alpha^2}F_\alpha(x,y)\, dy,
\end{eqnarray*}
where
$$F_\alpha(x,y):= \left(\frac{\min\{|\ra x|,|y|\}}{\max\{|\ra x|,|y|\}}\right)^{|\bm|}|\ra x-y|^{2-n}\left(\frac{\ma}{\ma^2+|y-\xa|^2}\right)^{\frac{n+2}{2}},$$
$$D_\alpha^1:=B_R(0)\cap\left\{|\ra x-y|>\frac{|\ra x-\xa|}{2}\right\}\hbox{ and }D_\alpha^2:=B_R(0)\cap\left\{|\ra x-y|\leq \frac{|\ra x-\xa|}{2}\right\}.$$
We estimate these two integrals separately. With a change of variable $y=\xa+\ma z$, we get that 
\begin{eqnarray*}
\frac{\ra^{n-2}}{\ma^{\frac{n-2}{2}}}\int_{D_\alpha^1}F_\alpha(x,y)\, dy \leq |x-\theta_\alpha|^{2-n}\int_{B_{2R\ma^{-1}}(0)}\left(\frac{\min\{|x|,|\theta_\alpha+\frac{\ma}{\ra} z|\}}{\max\{| x|,|\theta_\alpha+\frac{\ma}{\ra} z|\}}\right)^{|\bm|}\left(\frac{1}{1+|y|^2}\right)^{\frac{n+2}{2}}\, dy.
\end{eqnarray*}
It follows from Lebesgue's convergence theorem that
\begin{equation}\label{est:f:1}
\limsup_{\alpha\to +\infty}\frac{\ra^{n-2}}{\ma^{\frac{n-2}{2}}}\int_{D_\alpha^2}F_\alpha(x,y)\, dy\leq C|x-\theta_\infty|^{2-n}\left(\frac{\min\{|x|,1\}}{\max\{| x|,1\}}\right)^{|\bm|}
\end{equation}
For the second integral, we use that $|y-\xa|\geq |\ra x-\xa|/2$ for all $y\in D_\alpha^2$ to write  
\begin{eqnarray}\label{est:f:2}
&&\frac{\ra^{n-2}}{\ma^{\frac{n-2}{2}}}\int_{D_\alpha^1}F_\alpha(x,y)\, dy\\
&&\qquad \qquad \leq \frac{\ma^2}{\ra^2}|x-\theta_\alpha|^{-n-2}\int_{B_{C_x}(0)}\left(\frac{\min\{| x|,|z|\}}{\max\{| x|,|z|\}}\right)^{|\bm|}| x-z|^{2-n}\, dz,\nonumber 
\end{eqnarray}
for some $C_x>0$. Putting together \eqref{est:f:1} and \eqref{est:f:2} and letting $\alpha\to +\infty$ yields
\begin{equation}\label{est:123}
\tilde{u}(x)\leq C|x-\theta_\infty|^{2-n}\left(\frac{\min\{|x|,1\}}{\max\{| x|,1\}}\right)^{|\bm|}\hbox{ for all }x\in \rn\setminus\{0,\theta_\infty\}.
\end{equation}
We now prove a local reverse inequality. Since $G_\alpha\geq 0$, Green's representation \eqref{green:ra}, the lower bound \eqref{est:G:low}, the limit $\lim_{\alpha\to +\infty}\xa=0$ and a change of variable yield
\begin{eqnarray*}
\tua(x)&\geq& \frac{\ra^{n-2}}{\ma^{\frac{n-2}{2}}}\int_{B_{\ma}(\xa)} G_\alpha(\ra x,y)\la \ua^{\crit-1}(y)\, dy\\
&\geq& c\frac{\ra^{n-2}}{\ma^{\frac{n-2}{2}}} \int_{B_{\ma}(\xa)} \left(\frac{\min\{|\ra x|,|y|\}}{\max\{|\ra x|,|y|\}}\right)^{|\bm|}|\ra x-y|^{2-n}\ua^{\crit-1}(y)\, dy\\
&\geq &c \int_{B_{1}(0)} \left(\frac{\min\{|\ra x|,|\xa+\ma z|\}}{\max\{|\ra x|,|\xa+\ma z|\}}\right)^{|\bm|}| x-\theta_\alpha-\ma \ra^{-1}z|^{2-n} \va^{\crit-1}(z)\, dz
\end{eqnarray*}
Since $\ma=o(\ra)$ as $\alpha\to +\infty$, then, for $x\in \rn\setminus\{0,\theta_\infty\}$, as $\alpha\to +\infty$, we get that
\begin{equation}\label{ineq:tu}
\tilde{u}(x)\geq c\left(\frac{\min\{| x|,1\}}{\max\{| x|,1\}}\right)^{|\bm|}| x-\theta_\infty|^{2-n}\hbox{ for all }x\in \rn\setminus\{0,\theta_\infty\}.
\end{equation}
In particular, around $\theta_\infty\neq 0$, $\tilde{u}$ is controled from above and below by $|\cdot-\theta_\infty|^{2-n}$. It then follows from equation \eqref{eq:tu:2} and the classical classification of singular solutions of elliptic equations that there exists $\chi>0$ such that
$\tilde{u}(x)\sim_{x\to \theta_\infty}\frac{\chi}{(n-2)\omega_{n-1}|x-\theta_\infty|^{n-2}}$.
Integrating by parts, it follows from the pointwise control \eqref{est:123} and the equation \eqref{eq:tu:2} that
$$\chi \varphi(\theta_\infty)=\int_{\rn}\tilde{u}(x)\left(-\Delta \varphi-\frac{\gamma }{|x|^2}\varphi\right)\, dx \qquad \hbox{ for all }\varphi\in C^\infty_c(\rn).$$
The uniqueness result of Theorem \ref{th:green:gamma} then yields that $\tilde{u}=\chi\cdot G_{\theta_\infty}$. This complete the proof of \eqref{step3one}. 

\medskip\noindent
Next, we show that 
\begin{equation}\label{lim:nonzero}
\lim_{\alpha\to +\infty}\xa=x_0\neq 0.
\end{equation}
Indeed, otherwise we can assume that $\ra:=|\xa|\to 0$ as $\alpha\to +\infty$, so that \eqref{step3one} applies. We define $\tua$ as in \eqref{def:tua}. For $\delta\in (0,1)$, the Pohozaev identity \eqref{gen:p:id} applied on $B_{\delta\ra}(\xa)\subset\subset \Omega\setminus\{0\}$ with $p:=\xa$ and combined with \eqref{eq:ua:bis} yield
\begin{eqnarray}\label{poho:4}
 &-&  \int_{B_{\delta\ra}(\xa)}\left(\gamma \frac{(x,\xa)}{|x|^4}+a_\alpha+\frac{(x-\xa)^i\partial_i a_\alpha}{2}\right)\ua^2\, dx\nonumber \\
&=& \int_{\partial B_{\delta\ra}(\xa)}(x-\xa,\nu)\left(\frac{|\nabla \ua|^2}{2}-\left(\frac{\gamma}{|x|^2}+a_\alpha\right)\frac{\ua^2}{2}-\frac{\la\ua^{\crit}}{\crit}\right)\, d\sigma\nonumber\\
&& - \int_{\partial B_{\delta\ra}(\xa)}\left((x-\xa)^i\partial_i \ua+\frac{n-2}{2}\ua\right)\partial_\nu\ua\, d\sigma\nonumber\\
&=& \left(\frac{\ma}{\ra}\right)^{n-2}\int_{\partial B_{\delta}(\theta_\alpha)}(x-\theta_\alpha,\nu)\left(\frac{|\nabla \tua|^2}{2}-\left(\frac{\gamma}{|x|^2}+\ra^2a_\alpha(\ra x)\right)\frac{\tua^2}{2}\right)\, d\sigma \nonumber\\
&& -\left(\frac{\ma}{\ra}\right)^{n}\int_{\partial B_{\delta}(\theta_\alpha)}(x-\theta_\alpha,\nu)\frac{\la\tua^{\crit}}{\crit}\, d\sigma\nonumber\\
&&-\left(\frac{\ma}{\ra}\right)^{n-2}\int_{\partial B_{\delta}(\theta_\alpha)}\left((x-\theta_\alpha)^i\partial_i \tua+\frac{n-2}{2}\tua\right)\partial_\nu\ua\, d\sigma.
\end{eqnarray}
where $\theta_\alpha$ is defined in \eqref{def:thetaa}. In particular $|\theta_\alpha|=1$.

\medskip\noindent We first assume that $n\geq 4$. The convergence \eqref{lim:tua} of $\tua$ and $\delta<1$ yield
\begin{equation}\label{est:poho:1}
\int_{B_{\delta\ra}(\xa)}\left(\gamma \frac{(x,\xa)}{|x|^4}+a_\alpha+\frac{(x-\xa)^i\partial_i a_\alpha}{2}\right)\ua^2\, dx=O\left(\left(\frac{\ma}{\ra}\right)^{n-2}\right)
\end{equation}
The change of variable $x=\xa+\ma y$ yield
\begin{eqnarray}\label{est:poho:3}
&&\int_{B_{\delta\ra}(\xa)}\left(\gamma \frac{(x,\xa)}{|x|^4}+a_\alpha+\frac{(x-\xa)^i\partial_i a_\alpha}{2}\right)\ua^2\, dx\\
&&=\frac{\ma^2}{\ra^2}\int_{B_{\delta\ra/\ma}(0)}\left(\gamma \frac{(\theta_\alpha+\ma\ra^{-1} y,\theta_\alpha)}{|\theta_\alpha+\ma\ra^{-1} y|^4}+\ra^2 a_\alpha(\xa+\ma y)\right)\va^2\, dx.\nonumber\\
&&\quad+\frac{\ma^2}{\ra^2}\int_{B_{\delta\ra/\ma}(0)}\ra^2\frac{\ma y^i\partial_i a_\alpha(\xa+\ma y)}{2}\va^2\, dx.\nonumber
\end{eqnarray}
Since $\va(x)\leq C(1+|x|^2)^{-1-n/2}$ from \eqref{est:ua:point}, and when $n\geq 5$, Lebesgue's convergence theorem yields
\begin{equation}\label{est:poho:2}
\int_{B_{\delta\ra}(\xa)}\left(\gamma \frac{(x,\xa)}{|x|^4}+a_\alpha+\frac{(x-\xa)^i\partial_i a_\alpha}{2}\right)\ua^2\, dx=\left(\gamma \int_{\rn}v^2\, dx+o(1)\right)\frac{\ma^2}{\ra^2}
\end{equation}
as $\alpha\to +\infty$. Hence $\gamma\neq 0$, \eqref{est:poho:1} and \eqref{est:poho:2} yield 
$$1=O\left(\frac{\ma}{\ra}\right)^{n-4}\hbox{ as }\alpha\to +\infty,$$
which is a contradiction.\\
Now if $n=4$, we use that $|\gamma|>0$, \eqref{est:poho:1} and \eqref{est:poho:3} to get that there exists $C>0$ such that for any $R>0$, we have that
$$\int_{B_R(0)}\va^2\, dx\leq C\hbox{ for all }\alpha.$$
Letting $\alpha\to+\infty$ and then $R\to +\infty$ yields $\int_{\rn}v^2\, dx<+\infty$, a contradiction with \eqref{result:step1} that settles the case $n=4$.  

\medskip\noindent We now deal with the case $n=3$. With the pointwise control \eqref{step3one}, we have that
\begin{eqnarray*}
&&\left|\int_{B_{\delta\ra}(\xa)}\left(\gamma \frac{(x,\xa)}{|x|^4}+a_\alpha+\frac{(x-\xa)^i\partial_i a_\alpha}{2}\right)\ua^2\, dx\right|\\
&&\qquad \leq C\ra^{-2}\int_{B_{\delta\ra}(\xa)}\ma|x-\xa|^{-2}\, dx\leq C\delta\frac{\ma}{\ra} \nonumber
\end{eqnarray*}
Plugging this inequality in \eqref{poho:4}, using the convergence \eqref{lim:tua} and letting $\alpha\to +\infty$ yield that the expression 
\begin{equation}\label{poho:G}
\qquad  \int_{\partial B_{\delta}(\theta_\infty)}\left[(x-\theta_\infty,\nu)\left(\frac{|\nabla G_{\theta_\infty}|^2}{2}-\frac{\gamma G_{\theta_\infty}^2}{2|x|^2}\right)-\left((x-\theta_\infty)^i\partial_i G_{\theta_\infty}+\frac{n-2}{2}G_{\theta_\infty}\right)\partial_\nu G_{\theta_\infty}\right]\, d\sigma
\end{equation}
is $O(\delta)$ as $\delta\to 0$. 
Since $n=3$, there exists $\beta_{\theta_\infty}\in C^2(\rn\setminus\{0\})$ such that
$$G_{\theta_\infty}(x)=\frac{1}{4\pi|x-\theta_\infty|}+\beta_{\theta_\infty}(x)\hbox{ for all }x\in \rn\setminus\{0,\theta_\infty\}.$$
Letting $\delta\to 0$ in \eqref{poho:G}, classical computations then yield
\begin{equation}\label{zero:mass}
\beta_{\theta_\infty}(\theta_\infty)=0.
\end{equation}
We shall give an integral expression for $\beta_{\theta_\infty}$. Since $n=3$, it follows from the pointwise control \eqref{est:G:glob} and from the definition that $\beta_{\theta_\infty}\in \dundeuxrn$ and is controled at $\infty$ by $x\mapsto |x|^{-1}$. Since $-\Delta\beta_{\theta_\infty}-\gamma|x|^{-2}\beta_{\theta_\infty}=-\gamma|x|^{-2}|x-\theta_\infty|^{-1}/{4\pi}$, integrating by parts yields
$$\beta_\infty(x)=\frac{-\gamma}{4\pi}\int_{\rr^3}\frac{G_{x}(y)}{|y|^2|y-\theta_\infty|}\, dy$$
for all $x\in \rn$. Since $\gamma<0$, we then get that $\beta_\infty>0$, contradicting \eqref{zero:mass}. This proves \eqref{lim:nonzero} also when $n=3$.

\medskip\noindent
We now show that 
\begin{equation}\label{precise}
\hbox{If $x_0\in\partial\Omega$, then $n\geq 4\hbox{ and }\hat{a}_\infty(x_0)=0,$ where $\hat{a}_\infty:=a_\infty+\gamma|\cdot|^{-2}$.}
\end{equation}
Indeed, 
 let $U,V\subset \rn$ be open sets such that $0\in U$, $x_0\in V$ and $\varphi: U\to V$ a smooth diffeomorphism such that
$$\varphi(0)=x_0\; ; \; \varphi(U\cap\rnm)=\varphi(U)\cap\Omega\hbox{ and }\varphi(U\cap\partial\rnm)=\varphi(U)\cap\partial\Omega.$$
Up to a rotation, we can assume that the differential of $\varphi$ at $0$ is $d\varphi_0=Id_{\rn}$. Let $(x_{\alpha,1}, \bar{x}_\alpha)\in U\cap\rnm=U\cap ((-\infty,0)\times\rr^{n-1})$ be such that $x_\alpha=\varphi(x_{\alpha,1},\bar{x}_\alpha)$. In particular, we have that $d(\xa,\partial\Omega)=(1+o(1))|x_{\alpha,1}|$ as $\alpha\to +\infty$. Set $d_\alpha:=|x_{\alpha,1}|=-x_{\alpha,1}$, and then $x_0\in\partial\Omega$ and \eqref{lim:bndy} yield
\begin{equation}\label{lim:d}
\lim_{\alpha\to +\infty}\frac{d_\alpha}{\ma}=+\infty\hbox{ and }\lim_{\alpha\to +\infty}d_\alpha=0.
\end{equation}
Let
$$\tua(x):=\frac{\da^{n-2}}{\ma^{\frac{n-2}{2}}}\ua\circ\varphi((0,\bar{x}_\alpha)+\da x)\quad \hbox{ for }x\in \frac{(U-(0,\bar{x}_\alpha))\cap\rnm}{\da},$$
and for convenience, define for each $\alpha\in\nn$, the function 
$$\hat{a}_\alpha(x):=a_\alpha(x)+\gamma|x|^{-2}\quad \hbox{ for }x\in \overline{\Omega}\setminus\{0\}.$$
Equation \eqref{eq:ua:bis} rewrites
\begin{equation*}
\left\{\begin{array}{ll}
-\Delta_{g_\alpha}\tua-\da^2\hat{a}_\alpha\circ\varphi((0,\bar{x}_\alpha)+\da x) \tua=\left(\frac{\ma}{\da}\right)^2\tua^{\crit-1}&\hbox{ in }\frac{(U-(0,\bar{x}_\alpha))\cap\rnm}{\da}\\
\hfill \tua=0 \qquad \qquad \quad &\hbox{ in }\frac{(U-(0,\bar{x}_\alpha))\cap\partial\rnm}{\da}.
\end{array}\right.
\end{equation*}
Moreover, the pointwise control \eqref{est:ua:point} reads 
$$\tua(x)\leq C\left|x-\frac{(x_{\alpha,1},0)}{\da}\right|^{2-n}\hbox{ for }x\in\frac{(U-(0,\bar{x}_\alpha))\cap\rnm}{\da}.$$
It then follows from classical elliptic theory that there exists $\tilde{u}\in C^2(\overline{\rnm}\setminus\{(-1,0)\})$ such that
\begin{eqnarray*}
&&\lim_{\alpha\to +\infty}\tua=\tilde{u} \hbox{ in }C^2_{loc}(\overline{\rnm}\setminus\{(-1,0)\})\\
&&\Delta \tilde{u}=0 \hbox{ in }\overline{\rnm}\setminus\{(-1,0)\}\; ;\; \tilde{u}_{|\partial\rnm}\equiv 0\\
&& 0\leq \tilde{u}(x)\leq C |x-(-1,0)|^{2-n}\hbox{ for all }x\in \rnm\setminus\{(-1,0)\}
\end{eqnarray*}
By reflecting $\tilde{u}$ along the hyperplane $\{x_1=0\}$, we get a harmonic function on $\rn\setminus\{(\pm1,0)\}$, which is nonnegative for $x_1<0$, nonpositive for $x_1>0$, and vanishing for $x_1=0$. Therefore, there exists $c\geq 0$ such that
$$\tilde{u}(x)=c\left(|x-(-1,0)|^{2-n}-|x-(1,0)|^{2-n}\right)\hbox{ for all }x\in \rnm\setminus\{(-1,0)\}.$$
A proof that is similar to the one for \eqref{ineq:tu} and using the pointwise control of the Green's function in Robert \cite{r1} --and that we omit it here-- gives that $c>0$. 
Fix now $0<\delta<1$ and define
$$\hat{u}_\alpha(x):=\frac{\da^{n-2}}{\ma^{\frac{n-2}{2}}}\ua(\xa+\da x)\hbox{ for }x\in B_\delta(0).$$
It follows from the convergence result above that 
\begin{equation}\label{cv:hat:u}
\lim_{\alpha\to +\infty}\hat{u}_\alpha=\hat{u}:=c|\cdot|^{2-n}+\phi\hbox{ in }C^2_{loc}(B_\delta(0)\setminus\{0\}),
\end{equation}
where $\phi(x)=-c|x-(2,0)|^{2-n}$ for all $x\in B_\delta(0)$. We now use the Pohozaev identity \eqref{gen:p:id} on $B_{\delta\da}(\xa)$,  equation \eqref{eq:ua:bis} and an integration by parts, to obtain
\begin{eqnarray}\label{id:poho:5}
&&-\int_{B_{\delta\da}(\xa)}\left(\hat{a}_\alpha+\frac{(x-p)^i\partial_i\hat{a}_\alpha}{2}\right) \ua^2\, dx\\
&&\qquad \qquad= \int_{\partial B_{\delta\da}(\xa)}(x-p,\nu)\left(\frac{|\nabla \ua|^2}{2}-\frac{\hat{a}_\alpha \ua^{2}}{2}-\frac{\la \ua^{\crit}}{\crit}\right)d\sigma.\nonumber\\
&&\qquad \qquad \quad-\int_{\partial B_{\delta\da}(\xa)}\left((x-p)^i\partial_i \ua+\frac{n-2}{2}\ua\right)\partial_\nu \ua\, d\sigma.\nonumber
\end{eqnarray}
Taking $p:=\xa$ in this identity yields
\begin{eqnarray*}
&&-\int_{B_{\delta\da}(\xa)}\left(\hat{a}_\alpha+\frac{(x-\xa)^i\partial_i\hat{a}_\alpha}{2}\right) \ua^2\, dx\\
&&= \int_{\partial B_{\delta\da}(\xa)}\left[(x-\xa,\nu)\left(\frac{|\nabla \ua|^2}{2}-\frac{\hat{a}_\alpha \ua^{2}}{2}-\frac{\la \ua^{\crit}}{\crit}\right)-\left((x-\xa)^i\partial_i \ua+\frac{n-2}{2}\ua\right)\partial_\nu \ua\right]\, d\sigma.\nonumber
\end{eqnarray*}
With the change of variable $x=\xa+\ma y$ in the first integral, and $x=\xa+\da z$ in the second integral, we get that
\begin{eqnarray*}
&&\qquad -\ma^2\int_{B_{\delta\da\ma^{-1}}(0)}\left(\hat{a}_\alpha(\xa+\ma y)+\frac{\ma y^i\partial_i\hat{a}_\alpha(\xa+\ma y)}{2}\right) \va^2\, dx\\
&& \qquad = \left(\frac{\ma}{\da}\right)^{n-2}\int_{\partial B_{\delta}(0)}(z,\nu)\left(\frac{|\nabla \hat{u}_\alpha|^2}{2}-\frac{\da^2\hat{a}_\alpha (\xa+\da z)\hat{u}_\alpha^{2}}{2}-\frac{\la\ma^2 \hat{u}_\alpha^{\crit}}{\crit\da^2}\right)\, d\sigma\nonumber\\
&&\quad \qquad -\left(\frac{\ma}{\da}\right)^{n-2}\int_{\partial B_{\delta}(0)}\left(z^i\partial_i \hat{u}_\alpha+\frac{n-2}{2}\hat{u}_\alpha\right)\partial_\nu \hat{u}_\alpha\, d\sigma.\nonumber
\end{eqnarray*}
Fix $i\in \{1,...,n\}$ and differentiate \eqref{id:poho:5} with respect to the $i^{th}$ variable $p_i$ to obtain
\begin{eqnarray*}
&&-\int_{B_{\delta\da}(\xa)}\left(\frac{\partial_i\hat{a}_\alpha}{2}\right) \ua^2\, dx = \int_{\partial B_{\delta\da}(\xa)}\nu_i\left(\frac{|\nabla \ua|^2}{2}-\frac{\hat{a}_\alpha \ua^{2}}{2}-\frac{\la \ua^{\crit}}{\crit}\right)\, d\sigma\nonumber\\
&& \qquad -\int_{\partial B_{\delta\da}(\xa)}\partial_i \ua\partial_\nu \ua\, d\sigma.\nonumber
\end{eqnarray*}
Performing the same changes of variables as above yields
\begin{eqnarray}\label{id:67}
&& -\ma^2\int_{B_{\delta\da\ma^{-1}}(0)}\partial_i\hat{a}_\alpha(\xa+\ma y)\va^2\, dx\\
&&= \da^{-1}\left(\frac{\ma}{\da}\right)^{n-2}\int_{\partial B_{\delta}(0)}\nu_i\left(\frac{|\nabla \hat{u}_\alpha|^2}{2}-\frac{\da^2\hat{a}_\alpha (\xa+\da z)\hat{u}_\alpha^{2}}{2}-\frac{\la\ma^2 \hat{u}_\alpha^{\crit}}{\crit\da^2}\right)\, d\sigma\nonumber \\
&& \qquad  -\da^{-1}\left(\frac{\ma}{\da}\right)^{n-2}\int_{\partial B_{\delta}(0)}\partial_i \hat{u}_\alpha\partial_\nu \hat{u}_\alpha\, d\sigma.\nonumber
\end{eqnarray}
With the convergence \eqref{cv:hat:u} of $\hat{u}_\alpha$ and an explicit computation, we get that
\begin{eqnarray*}
&&\lim_{\alpha\to +\infty}\int_{\partial B_{\delta}(0)}\left[(z,\nu)\left(\frac{|\nabla \hat{u}_\alpha|^2}{2}-\frac{\da^2\hat{a}_\alpha (\xa+\da z)\hat{u}_\alpha^{2}}{2}-\frac{\la\ma^2 \hat{u}_\alpha^{\crit}}{\crit\da^2}\right)-\left(z^i\partial_i \hat{u}_\alpha+\frac{n-2}{2}\hat{u}_\alpha\right)\partial_\nu \hat{u}_\alpha\right]\, d\sigma\\
&=&\int_{\partial B_{\delta}(0)}\left[(z,\nu)\frac{|\nabla \hat{u}|^2}{2}-\left(z^i\partial_i \hat{u}+\frac{n-2}{2}\hat{u}\right)\partial_\nu \hat{u}\right]\, d\sigma\\
&=&\frac{(n-2)^2c}{2}\phi(0)=-(n-2)^2 2^{1-n}c^2.
\end{eqnarray*}
Indeed, the limit is independent of $\delta$ since $\phi$ is harmonic. Similarly,
\begin{eqnarray}\label{id:68}
 &&\lim_{\alpha\to +\infty}\int_{\partial B_{\delta}(0)}\left[\nu_i\left(\frac{|\nabla \hat{u}_\alpha|^2}{2}-\frac{\da^2\hat{a}_\alpha (\xa+\da z)\hat{u}_\alpha^{2}}{2}-\frac{\la\ma^2 \hat{u}_\alpha^{\crit}}{\crit\da^2}\right)-\partial_i \hat{u}_\alpha\partial_\nu \hat{u}_\alpha\right]\, d\sigma\nonumber\\
&&\qquad =\int_{\partial B_{\delta}(0)}\left[\nu_i\frac{|\nabla \hat{u}|^2}{2}-\partial_i \hat{u}\partial_\nu \hat{u}\right]\, d\sigma\\
&&\qquad =(n-2)c\omega_{n-1}\partial_i\phi(0)=-2^{1-n}(n-2)^2c^2\omega_{n-1}\delta_{i,1}.\nonumber
\end{eqnarray}

\medskip\noindent We now divide the analysis in three cases.

\smallskip\noindent{\it Case 1: $n\geq 5$.} Since $\va\leq C(1+|x|^{2})^{1-n/2}$ from \eqref{est:ua:point}, and $\va\to v$ in $C^2_{loc}(\rn)$, then Lebesgue's theorem applied to the identities above yields
\begin{equation}\label{eq:165}
\ma^2\left( \hat{a}_\infty(x_0)\int_{\rn}v^2\, dx+o(1)\right)= \left(\frac{\ma}{\da}\right)^{n-2}\cdot\left((n-2)^2 2^{1-n}c^2+o(1)\right)
\end{equation}
and, with $i=1$,
\begin{equation}\label{eq:166}
\ma^2\left( \partial_1 \hat{a}_\infty(x_0)\int_{\rn}v^2\, dx+o(1)\right)=\da^{-1} \left(\frac{\ma}{\da}\right)^{n-2}\cdot\left(2^{1-n}(n-2)^2c^2\omega_{n-1}
+o(1)\right)
\end{equation}
In particular, we get that $\hat{a}_\infty(x_0)=0$.

\smallskip\noindent{\it Case 2: $n=4$.} Arguing as in the case $n\geq 5$, we get that
\begin{equation*}
\ma^2\left( (\hat{a}_\infty(x_0)+o(1))\int_{B_{\delta\da\ma^{-1}}(0)}\va^2\, dx+O(1)\right)= \left(\frac{\ma}{\da}\right)^{2}\cdot\left((n-2)^2 2^{1-n}c^2+o(1)\right)
\end{equation*}
and for $i=1$,
\begin{eqnarray*}
&&\ma^2\left( (\partial_1 \hat{a}_\infty(x_0)+o(1))\int_{B_{\delta\da\ma^{-1}}(0)}\va^2\, dx+O(1)\right) =\da^{-1} \left(\frac{\ma}{\da}\right)^{2}\cdot\left(2^{1-n}(n-2)^2c^2\omega_{n-1}
+o(1)\right).\nonumber
\end{eqnarray*}
Since $\int_{\rn}v^2\, dx=+\infty$ when $n=4$, here again, we get that $\hat{a}_\infty(x_0)=0$.

\smallskip\noindent{\it Case 3: $n=3$.} Here we need to show that $x_0\notin \partial \Omega$. Indeed, the uniform control $\va\leq C(1+|x|^{2})^{1-n/2}$, the estimates \eqref{id:67} and \eqref{id:68} yield
$O(\ma \da)=\frac{\ma}{\da}\cdot \left(-(n-2)^2 2^{1-n}c^2\right)$ and therefore $1=O(\da^2)$, contradicting \eqref{lim:d}.

The proof of \eqref{precise} is complete.

\medskip\noindent
Assume now that $n\geq 4$ and $x_0\in \Omega$, set for convenience $\hat{a}_\alpha(x):=a_\alpha(x)+\gamma|x|^{-2}$. Performing the Pohozaev identity \eqref{gen:p:id} on $B_{\delta}(\xa)$ assuming that $B_{2\delta}(\xa)\subset \Omega$, we get that
\begin{eqnarray}\label{step6}
&&-\int_{B_{\delta}(\xa)}\left(\hat{a}_\alpha+\frac{(x-\xa)^i\partial_i \hat{a}_\alpha}{2}\right)\ua^2\, dx\\
&&\quad= \int_{\partial B_{\delta}(\xa)}(x-\xa,\nu)\left(\frac{|\nabla \ua|^2}{2}-\hat{a}_\alpha\frac{\ua^2}{2}-\frac{\la\ua^{\crit}}{\crit}\right)\, d\sigma \nonumber\\
&&\qquad - \int_{\partial B_{\delta}(\xa)}\left((x-\xa)^i\partial_i \ua+\frac{n-2}{2}\ua\right)\partial_\nu\ua\, d\sigma. \nonumber\end{eqnarray}
The pointwise control \eqref{est:ua:point} and elliptic theory yield $\ua(x)+|\nabla\ua(x)|\leq C\ma^{\frac{n-2}{2}}$ for $x\in \partial B_{\delta}(\xa)$, and therefore as $\alpha\to +\infty$, 
$$\int_{B_{\delta}(\xa)}\left(\hat{a}_\alpha+\frac{(x-\xa)^i\partial_i \hat{a}_\alpha}{2}\right)\ua^2\, dx=O(\ma^{n-2}).$$
Arguing as in the cases $n\geq 5$ and $n=4$ in the proof of \eqref{eq:165} and \eqref{eq:166} above, we then get that
$$\hat{a}_\infty(x_0)=a_\infty(x_0)+\gamma|x_0|^{-2}=0.$$

\medskip\noindent
Finally, assume that $n=3$. It follows from the above that $x_0\neq 0$ and $x_0\not\in\partial\Omega$. Therefore $(\hat{a}_\alpha)$ converges to $\hat{a}_\infty$ in $C^1(B_{2\delta}(x_0))$ for some small $\delta>0$. Passing to the limit as $\alpha\to +\infty$ and $\delta\to 0$ in \eqref{step6} above, and performing standard computations (see for instance Druet \cite{d2}), we get that the mass of the operator $-\Delta-(a_\infty+\gamma|x|^{-2})$ vanishes at $x_0$. In other words, $R_{\gamma,a_\infty}(x_0)=0$. This completes the proof of Theorem \ref{th:blowup:2}.

\section{\, Proof of Theorem \ref{gamma.star}}

Again, we start with the truly singular case and prove the following.

\begin{proposition}Let $\Omega$ be a smooth bounded domain in $\rn$ ($n\geq 3$) such that $0\in  \Omega$. Assume that either $s>0$ or $\gamma>0$. 
  If $\frac{(n-2)^2}{4}-1  <\gamma <\frac{(n-2)^2}{4}$, then
  \begin{enumerate}
\item  
$\lambda^*(\Omega)>0.$
 \item Moreover, if  $\mu_{\gamma, s, \lambda^*}(\Omega)$ is not achieved, then $m_{\gamma, \lambda^*}(\Omega)= 0$, and 
 \begin{equation*}
 \lambda^*(\Omega)=\sup \{\lambda;\,  m_{\gamma, \lambda}(\Omega)\leq 0\}.
 \end{equation*}  
  \end{enumerate}
\end{proposition}
\noindent{\it Proof:} For $\lambda>\lambda^*$, the infimum is achieved, and therefore, there exists $u_\lambda\in \huno$ such that
\begin{equation}\label{eq:ul}
\left\{\begin{array}{ll}
-\Delta u_\lambda -\left(\frac{\gamma}{|x|^2}+\lambda\right) u_\lambda = \mu_{\gamma,s,\lambda}(\Omega) \frac{u_\lambda^{\crits-1}}{|x|^s}& \hbox{ in }\Omega\\
u_\lambda\geq 0& \hbox{ a.e. in }\Omega\\
u_\lambda=0& \hbox{ on }\partial \Omega
\end{array}\right.
\end{equation}
and
\begin{equation}\label{mass:ul}
\int_\Omega \frac{u_\lambda^{\crits}}{|x|^s}\, dx=1.
\end{equation}
As one checks, $(u_\lambda)_{\lambda>\lambda^*}$ is bounded in $\huno$, and therefore, up to extracting a sub-family, it has a weak limit $u_{\lambda^*}$ as $\lambda\to \lambda^*$. If $u_{\lambda^*}\not\equiv 0$, then classical arguments yield it is a minimizer for $\mu_{\gamma,s,\lambda^*}(\Omega)$.

Suppose now $\lambda^*=0$, this means that $u_{\lambda^*}$ is a minimizer for $\mu_{\gamma,s, 0}(\Omega)=\mu_{\gamma,s, 0}(\rn)$, which is impossible since $u_{\lambda^*}$ has compact support, hence $u_{\lambda^*}\equiv 0$. It then follows from Theorem \ref{th:blowup:1} that $m_{\gamma,\lambda^*}(\Omega)=m_{\gamma, 0}(\Omega)=0.$ To get to a contradiction, we shall now prove that $m_{\gamma,0}(\Omega)<0$. Indeed, let $H\in C^2(\Omegabar\setminus\{0\})$ be as in Proposition \ref{prop:def:mass} for $h\equiv 0$. It follows from the definition of $H$ and the expansion \eqref{exp:H} that $x\mapsto H'(x):=H(x)-|x|^{-\bp}\in H_1^2(\Omega)\cap C^2(\Omegabar\setminus\{0\})$, it satisfies $-\Delta H'-\gamma|x|^{-2}H'=0$ in $\Omega\setminus\{0\}$ and $H'(x)<0$ for $x\in \partial\Omega$. It then follows from the comparison principle that $H'<0$ in $\Omega\setminus\{0\}$. Therefore, the expression \eqref{exp:H} yields that $c_2<0$, and therefore $m_{\gamma,0}(\Omega)<0$. A contradiction that yields that  $\lambda^*>0$.  

We now show (2) under the hypothesis that $\mu_{\gamma, s, \lambda^*}(\Omega)$ is not achieved. Indeed, under such an assumption, the weak limit $u_{\lambda^*}$ as $\lambda\to \lambda^*$ is necessarily identically zero. It then follows from Theorem \ref{th:blowup:1} that $m_{\gamma,\lambda^*}(\Omega)=0.$

Finally, let $\bar \lambda:=\sup \{\lambda;\,  m_{\gamma, \lambda}(\Omega)\leq 0\},$ and note that if $\lambda >\bar \lambda$, then $m_{\gamma, \lambda}(\Omega) >0$ and $\mu_{\gamma, s, \lambda^*}(\Omega)$ is achieved in view of Theorem \ref{positivemass:2}, which means that $\lambda \geq \lambda^*$. In other words, $\bar \lambda \geq \lambda^*$. On the other hand, from the strict monotonicity of the mass, if $\bar \lambda > \lambda^*$, then $m_{\gamma, \bar \lambda}(\Omega) > m_{\gamma,  \lambda^*}(\Omega)=0$, which is a contradiction, hence  $\bar \lambda = \lambda^*$. $\Box$.

An identical proof that uses Theorem \ref{th:blowup:2} as opposed to Theorem \ref{th:blowup:1}, and the mass $R_{\gamma, \lambda}(\Omega)$ as opposed to $m_{\gamma, \lambda}(\Omega)$ gives the analogous result in the merely singular case. Note that in this case, the argument of Druet \cite{d2} yields that $\mu_{\gamma, 0, \lambda^*}(\Omega)$ is not achieved for $n=3$, and therefore this hypothesis is readily satisfied. In summary, we have shown the following result.

\begin{proposition} \label{prop:smalldim:1}
Assume $n=3$, $s=0$ and $\gamma \leq 0$.  Then
\begin{equation*}
\lambda^*(\Omega)=\sup \{\lambda;\,  R_{\gamma, \lambda}(\Omega)\leq 0\}>0.
\end{equation*}
\end{proposition}

\section*{\, Appendix A: Green's function for $-\Delta-\gamma|x|^{-2}-h(x)$ on a bounded domain}

\begin{theorem}\label{th:green:gamma:domain} Let $\Omega$ be a smooth bounded domain of $\rn$ such that $0\in\Omega$ is an interior point. We fix $\gamma<\frac{(n-2)^2}{4}$. We let $h\in C^{0,\theta}(\Omegabar)$ be such that $-\Delta-\gamma|x|^{-2}-h$ is coercive. Then there exists $G: (\Omega\setminus\{0\})^2\setminus \{(x,x)/\, x\in \Omega\setminus\{0\}\}\to \rr$ such that for all $p\in\Omega\setminus\{0\}$, 

\smallskip\noindent{\bf (i)} For any $p\in\Omega\setminus\{0\}$, $G_p:=G(p,\cdot)\in H_{1}^2(\Omega\setminus B_\delta(p))$ for all $\delta>0$, $G_p\in C^{2,\theta}(\overline{\Omega}\setminus\{0,p\})$

\smallskip\noindent{\bf (ii)} For all $f\in L^{\frac{2n}{n+2}}(\Omega)\cap L^p_{loc}(\Omegabar-\{0\})$, $p>n/2$, and all $\varphi\in \huno$ such that
$$-\Delta\varphi-\left(\frac{\gamma}{|x|^2}+h(x)\right)\varphi=f\hbox{ in }\Omega\; ; \; \varphi_{|\partial\Omega}=0,$$
then we have that
\begin{equation}\label{id:172}
\varphi(p)=\int_{\Omega}G(p,x)f(x)\, dx
\end{equation}

\medskip\noindent In addition, $G>0$ is unique and 

\smallskip\noindent{\bf (iii)} For all $p\in\Omega\setminus\{0\}$, there exists $c_0(p)>0$ such that
\begin{equation}\label{asymp:G}
G_p(x)\sim_{x\to 0} \frac{c_0(p)}{|x|^{\bm}}\hbox{ and }G_p(x)\sim_{x\to p}\frac{1}{(n-2)\omega_{n-1}|x-p|^{n-2}}
\end{equation}
\smallskip\noindent{\bf (iv)}  There exists $c>0$ such that
\begin{equation}\label{est:G:up}
0< G_p(x)\leq c \left(\frac{\max\{|p|,|x|\}}{\min\{|p|,|x|\}}\right)^{\bm}|x-p|^{2-n}\hbox{ for }x\in\Omega-\{0,p\}.
\end{equation}
\smallskip\noindent{\bf (v)}  For all $\omega\Subset\Omega$, there exists $c(\omega)>0$ such that
\begin{equation}\label{est:G:low}
c(\omega) \left(\frac{\max\{|p|,|x|\}}{\min\{|p|,|x|\}}\right)^{\bm}|x-p|^{2-n}\leq G_p(x)\hbox{ for all }p,x\in\omega\setminus\{0\}.
\end{equation}
\end{theorem}

\medskip\noindent{\it Proof:} Fix $\delta_0>0$ such that $B_{\delta_0}(0)\subset \Omega$. We let $\eta_\eps(x):=\tilde{\eta}(\eps^{-1}|x|)$ for all $x\in\rn$ and $\eps>0$, where $\tilde{\eta}\in C^\infty(\rr)$ is nondecreasing and such that $\tilde{\eta}(t)=0$ for $t<1$ and $\tilde{\eta}(t)=1$ for $t>1$. Set
$$L_\eps:=-\Delta-\left(\frac{\gamma\eta_\eps}{|x|^2}+h(x)\right).$$
It follows from Lemma \ref{lem:vp} and the coercivity of $-\Delta-\left(\gamma|x|^{-2}+h\right)$ that there exists $\eps_0>0$ and $c>0$ such that such that for all $\varphi\in \huno$ and $\eps\in (0,\eps_0)$, 
$$\int_\Omega\left(|\nabla \varphi|^2-\left(\frac{\gamma\eta_\eps}{|x|^2}+h(x)\right)\varphi^2\right)\, dx\geq c\int_\Omega\varphi^2\, dx.$$
As a consequence, there exists $c>0$ such that for all $\varphi\in \huno$ and $\eps\in (0,\eps_0)$, 
\begin{equation}\label{bnd:coer}
\int_\Omega\left(|\nabla \varphi|^2-\left(\frac{\gamma\eta_\eps}{|x|^2}+h(x)\right)\varphi^2\right)\, dx\geq c\Vert\varphi\Vert_{D^{1,2}}^2.
\end{equation}
Let $G_\eps>0$ be the Green's function of $-\Delta-\left(\gamma\eta_\eps|x|^{-2}+h\right)$ on $\Omega$ with Dirichlet boundary condition. The existence follows from the coercivity and the $C^{0,\theta}$ regularity of the potential for any $\eps>0$.

\medskip\noindent{\bf Step 1: Integral bounds for $G_\eps$.} We claim that for all $\delta>0$ and $1<q<\frac{n}{n-2}$ and $\delta'\in (0,\delta)$, there exists $C(\delta,q)>0$ and $C(\delta,\delta')>0$ such that 
\begin{equation}\label{int:bnd:G}
\Vert G_\eps(x,\cdot)\Vert_{L^q(\Omega)}\leq C(\delta,q)\hbox{ and }\Vert G_\eps(x,\cdot)\Vert_{L^{\frac{2n}{n-2}}(\Omega\setminus B_{\delta'}(x))}\leq C(\delta,\delta')
\end{equation}
for all $x\in \Omega$, $|x|>\delta$.

\medskip\noindent Indeed, fix $f\in C^\infty_c(\Omega)$ and let $\varphi_\eps\in C^{2,\theta}(\overline{\Omega})$ be the solution to the boundary value problem
\begin{equation}
\label{eq:phi:eps}\left\{\begin{array}{ll}
L_\eps\varphi_\eps=-\Delta \varphi_\eps-\left(\frac{\gamma\eta_\eps}{|x|^2}+h(x)\right)\varphi_\eps= f &\hbox{ in }\Omega\\
\quad \varphi_\eps=0&\hbox{ on }\partial\Omega
\end{array}\right.
\end{equation}
Multiplying the equation by $\varphi_\eps$, integrating by parts on $\Omega$, using \eqref{bnd:coer} and H\"older's inequality, we get that 
$$\int_\Omega |\nabla\varphi|^2\, dx\leq C\Vert f\Vert_{\frac{2n}{n+2}}\Vert\varphi_\eps\Vert_{\frac{2n}{n-2}}$$
where $C>0$ is independent of $\epsilon$, $f$ and $\varphi_\epsilon$. The Sobolev inequality $\Vert\varphi\Vert_{\frac{2n}{n-2}}\leq C\Vert \nabla\varphi\Vert_2$ for $\varphi\in \huno$ then yields
\begin{equation*}
\Vert\varphi_\eps\Vert_{\frac{2n}{n-2}}\leq C\Vert f\Vert_{\frac{2n}{n+2}}
\end{equation*}
where $C>0$ is independent of $\epsilon$, $f$ and $\varphi_\epsilon$. 

Fix $p>n/2$ and $\delta\in (0,\delta_0)$ and $\delta_1,\delta_2>0$ such that $\delta_1+\delta_2<\delta$, and $x\in\Omega$ such that $|x|>\delta$. It follows from standard elliptic theory that 
\begin{eqnarray*}
|\varphi_\eps(x)|&\leq & \Vert \varphi\Vert_{C^0(B_{\delta_1}(x))}\\
&\leq & C\left(\Vert \varphi_\epsilon\Vert_{L^{\crit}(B_{\delta_1+\delta_2}(x))}+\Vert f\Vert_{L^{p}(B_{\delta_1+\delta_2}(x))}\right)\\
&\leq & C\left(\Vert f\Vert_{L^{\frac{2n}{n+2}}(\Omega)}+\Vert f\Vert_{L^{p}(B_{\delta_1+\delta_2}(x))}\right)
\end{eqnarray*}
where $C>0$ depends on $p,\delta,\delta_1,\delta_2$, $\gamma$ and $\Vert h\Vert_\infty$. Therefore, Green's representation formula yields
\begin{equation}\label{rep:G:f}
\left|\int_\Omega G_\eps(x,\cdot)f\, dy\right|\leq C\left(\Vert f\Vert_{L^{\frac{2n}{n+2}}(\Omega)}+\Vert f\Vert_{L^{p}(B_{\delta_1+\delta_2}(x))}\right)
\end{equation}
for all $f\in C^\infty_c(\Omega)$.  It follows from \eqref{rep:G:f} that
$$\left|\int_\Omega G_\eps(x,\cdot)f\, dy\right|\leq C\cdot\Vert f\Vert_{L^{p}(\Omega)}$$
for all $f\in C^\infty_c(\Omega)$ where $p>n/2$. It then follows from duality arguments that for any $q\in (1, n/(n-2))$ and any $\delta>0$, there exists $C(\delta,q)>0$ such that $\Vert G_\eps(x,\cdot)\Vert_{L^q(\Omega)}\leq C(\delta,q)$ for all $\eps<\eps_0$ and $x\in \Omega\setminus B_\delta(0)$.

\medskip\noindent Let $\delta'\in (0,\delta)$ and $\delta_1,\delta_2>0$ such that $\delta_1+\delta_2<\delta'$. We get from \eqref{rep:G:f}  that 
\begin{equation}\label{rep:G:f:2}
\left|\int_\Omega G_\eps(x,\cdot)f\, dy\right|\leq C\Vert f\Vert_{L^{\frac{2n}{n+2}}(\Omega\setminus B_{\delta'}(x))}
\end{equation}
for all $f\in C^\infty_c(\Omega\setminus B_{\delta'}(x))$. Here again, a duality argument yields \eqref{int:bnd:G}, which proves the claim in  Step 1.

\medskip\noindent{\bf Step 2: Convergence of $G_\epsilon$.} Fix $x\in \Omega\setminus\{0\}$. For $0<\eps<\eps'$, since $G_\eps(x,\cdot)$, $G_{\eps'}(x,\cdot)$ are $C^2$ outside $x$, we have  
$$-\Delta(G_\eps(x,\cdot)-G_{\eps'}(x,\cdot))-\left(\frac{\gamma\eta_\eps}{|\cdot|^2}+h\right)(G_\eps(x,\cdot)-G_{\eps'}(x,\cdot))= \frac{\gamma(\eta_\eps-\eta_{\eps'})}{|\cdot|^2}G_{\eps'}(x,\cdot)$$
in the strong sense. The coercivity \eqref{bnd:coer} then yields
$$G_\eps(x,\cdot)\geq G_{\eps'}(x,\cdot)\hbox{ for }0<\eps<\eps'\hbox{ if }\gamma\geq 0,$$
and the reverse inequality if $\gamma<0$. It then follows from the integral bound \eqref{int:bnd:G} and elliptic regularity that there exists $G(x,\cdot)\in C^{2,\theta}(\overline{\Omega}\setminus\{0,x\})$ such that
$$\lim_{\eps\to 0}G_\eps(x,\cdot)=G(x,\cdot)\hbox{ in }C^2_{loc}(\overline{\Omega}-\{0,x\}).$$
In particular, $G$ is symmetric and 
\begin{equation}\label{eq:G:x}
-\Delta G(x,\cdot)-\left(\frac{\gamma}{|\cdot|^2}+h\right)G(x,\cdot)=0\hbox{ in }\overline{\Omega}\setminus\{0,x\}. 
\end{equation}
Moreover, passing to the limit $\eps\to 0$ in \eqref{int:bnd:G} and using elliptic regularity, we get that for all $\delta>0$, $1<q<\frac{n}{n-2}$ and $\delta'\in (0,\delta)$, there exist $C(\delta,q)>0$ and $C(\delta,\delta')>0$ such that for all $x\in \Omega$, $|x|>\delta$,
\begin{equation}\label{int:bnd:G:bis}
\Vert G(x,\cdot)\Vert_{L^q(\Omega)}\leq C(\delta,q)\hbox{ and }\Vert G(x,\cdot)\Vert_{L^{\frac{2n}{n-2}}(\Omega\setminus B_{\delta'}(x))}\leq C(\delta,\delta').
\end{equation}
Moreover, for any $f\in L^p(\Omega)$, $p>n/2$, let $\varphi_\eps\in C^2(\overline{\Omega})$ be such that \eqref{eq:phi:eps} holds, and fix $x\in \Omega\setminus \{0\}$. Passing to the limit $\eps\to 0$ in the Green identity $\varphi_\eps(x)=\int_\Omega G_\eps(x,\cdot)f\, dy$ yields
\begin{equation}\label{id:regul}
\varphi(x)=\int_\Omega G(x,\cdot)f\, dy\hbox{ for all }x\in\Omega\setminus\{0\}
\end{equation}
where $\varphi\in \huno\cap C^{0}(\overline{\Omega}\setminus\{0\})$ is the only weak solution to 
$$\left\{\begin{array}{ll}
-\Delta \varphi-\left(\frac{\gamma}{|x|^2}+h(x)\right)\varphi= f &\hbox{ in }\Omega\\
\varphi=0&\hbox{ on }\partial\Omega
\end{array}\right.$$
In particular, the strong comparision principle yields $G(x,\cdot)>0$ for $x\in \Omega\setminus\{0\}$.

\medskip\noindent{\bf Step 3: Upper bound for $G(x,y)$ when one variable is far from $0$.}

 It follows from \eqref{eq:G:x}, elliptic theory and \eqref{int:bnd:G:bis} that for any $\delta>0$, there exists $C(\delta)>0$ such that
\begin{equation}\label{est:1}
0<G(x,y)\leq C(\delta)\hbox{ for }x,y\in\Omega\hbox{ such that }|x|>\delta,\, |y|>\delta,\, |x-y|>\delta.
\end{equation}
We claim that for any $\delta>0$, there exists $C(\delta)>0$ such that
\begin{equation}\label{est:2}
0<|x-y|^{n-2}G(x,y)\leq C(\delta)\hbox{ for }x,y\in\Omega\hbox{ such that }|x|>\delta\hbox{ and } |y|>\delta.
\end{equation}
Indeed, with no loss of generality, we can assume that $\delta\in (0,\delta_0)$. Define now $\Omega_\delta:=\Omega\setminus B_{\delta/2}(0)$, and fix $x\in\Omega$ such that $|x|>\delta$. Let $H_x$ be the Green's function for $-\Delta -\left(\frac{\gamma}{|x|^2}+h(x)\right)$ in $\Omega_\delta$ with Dirichlet boundary condition. Classical estimates (see \cite{r1}) yield the existence of $C(\delta)>0$ such that $|x-y|^{n-2}H_x(y)\leq C(\delta)$ for all $x,y\in \Omega_\delta$. It is easy to check that
$$\left\{\begin{array}{ll}
-\Delta (G_x-H_x)-\left(\frac{\gamma}{|x|^2}+h\right)(G_x-H_x)= 0 &\hbox{ weakly in }\Omega_\delta\\
G_x-H_x=0&\hbox{ on }\partial\Omega\\
G_x-H_x=G_x&\hbox{ on }\partial B_{\delta/2}(0).
\end{array}\right.$$
Regularity theory then yields that $G_x-H_x\in C^{2,\theta}(\overline{\Omega_\delta})$. It follows from \eqref{est:1} that $G_x$ is bounded by a constant depending only on $\delta$ on $\partial B_{\delta/2}(0)$ for $|x|>\delta$. The comparison principle then yields $|G_x(y)-H_x(y)|\leq C(\delta)$ for $y\in \Omega_\delta$ and $|x|>\delta$. The above bound for $H_x$ and \eqref{est:1} then yields \eqref{est:2}.

We now claim that for any $0<\delta'<\delta$, there exists $C(\delta,\delta')>0$ such that
\begin{equation}\label{est:3}
0<|y|^{\bm}G(x,y)\leq C(\delta,\delta')\hbox{ for }x,y\in\Omega\hbox{ such that }|x|>\delta>\delta'>|y|>0.
\end{equation}
Indeed, fix $\delta_1<\delta$ and use \eqref{est:1} to deduce that $G_x(y)\leq C(\delta,\delta_1)$ for all $x\in \Omega\setminus B_\delta(0)$ and $y\in \partial B_{\delta_1}(0)$. Since $\delta_1<|x|$, we have that
 $$\left\{\begin{array}{ll}
-\Delta G_x-\left(\frac{\gamma}{|x|^2}+h\right)G_x= 0 &\hbox{  in }H_1^2(B_{\delta_1}(0))\\
0<G_x\leq C(\delta,\delta')&\hbox{ on }\partial B_{\delta_1}(0).
\end{array}\right.$$
It follows from \eqref{ppty:super:2} below that for $\delta_1>0$ small enough, there exists $u_{\beta_-}\in H_1^2(B_{\delta_1}(0))$ such that $c_1\leq |z|^{\bm}u_{\beta_-}(z)<c_2$ for all $z\in B_{\delta_1}(0)$, and
$$-\Delta u_{\beta_-}-\left(\frac{\gamma}{|x|^2}+h\right)u_{\beta_-}\geq 0 \hbox{  in }H_1^2(B_{\delta_1}(0)).
$$
Therefore, there exists $C(\delta,\delta')>0$ such that $G_x(z)\leq C(\delta,\delta')u_{\beta_-}(z)$ for all $z\in \partial B_{\delta_1}(0)$. It then follows from the comparison principle that $G_x(y)\leq C(\delta,\delta')u_{\beta_-}(y)$ for all $y\in B_{\delta_1}(0)\setminus\{0\}$. Combining this with \eqref{est:1}, we obtain \eqref{est:3}.

Note that by symmetry, we also get that for any $0<\delta'<\delta$, there exists $C(\delta,\delta')>0$ such that
\begin{equation}\label{est:4}
|x|^{\bm}G(x,y)\leq C(\delta,\delta')\hbox{ for }x,y\in\Omega\hbox{ such that }|y|>\delta>\delta'>|x|>0.
\end{equation}

\medskip\noindent{\bf Step 4: Upper bound for $G(x,y)$ when both variables approach $0$.}\\

We claim first that for all $c_1,c_2,c_3>0$, there exists $C(c_1,c_2,c_3)>0$ such that for $x,y\in\Omega$ such that $c_1|x|<|y|<c_2|x|$ and $|x-y|>c_3|x|$, we have
 \begin{equation}\label{est:5}
|x-y|^{n-2}G(x,y)\leq C(c_1,c_2,c_3).\end{equation}
Indeed, fix $x\in B_{\delta_0/2}(0)\setminus\{0\}\subset \Omega\setminus\{0\}$, and define
$$H(z):=G_x(|x|z)\hbox{ for }z\in B_{\delta_0/|x|}(0)\setminus \left\{0,\frac{x}{|x|}\right\},$$
so that 
$$-\Delta H-\left(\frac{\gamma}{|z|^2}+|x|^2h(|x|z)\right)H=0\hbox{ in }B_{\delta_0/|x|}(0)\setminus\left \{0,\frac{x}{|x|}\right\}.$$
Since $H>0$, it follows from the Harnack inequality that for all $R>0$ large enough and $\alpha>0$ small enough, there exist $\delta_1>0$ and $C>0$ independent of $|x|<\delta_1$ such that
\begin{equation*}
H(z)\leq C H(z')\hbox{ for all }z,z'\in B_R(0)\setminus \left(B_\alpha(0)\cup B_\alpha\left(\frac{x}{|x|}\right)\right),
\end{equation*}
which rewrites as:
\begin{equation}\label{harnack:g}
G_x(y)\leq C G_x(y')\hbox{ for all }y,y'\in B_{R|x|}(0)\setminus \left(B_{\alpha|x|}(0)\cup B_{\alpha|x|}(x)\right).
\end{equation}
Let $u_{\beta_+}$ be a sub-solution to \eqref{ppty:super:3}. In particular, for $|x|<\delta_2$ small, there exists $C>0$ such that
$$G_x(z)\geq c|x|^{\bp}\left(\inf_{\partial B_{R|x|}(0)}G_x\right) u_{\beta_+}(z)\hbox{ for all }z\in \partial B_{R|x|}(0).$$
Since $-\Delta G_x-(\gamma|\cdot|^{-2}+h)G_x=0$ outside $0$, it follows from coercivity and the comparison principle that
$$G_x(z)\geq c|x|^{\bp}\left(\inf_{\partial B_{R|x|}(0)}G_x\right) u_{\beta_+}(z)\hbox{ for all }z\in \Omega\setminus B_{R|x|}(0).$$
Fix $z_0\in \Omega\setminus\{0\}$. Then for $\delta_3$ small enough, it follows from \eqref{est:4} and the Harnack inequality \eqref{harnack:g} that there exists $C>0$ independent of $x$ such that
$$G_x(y)\leq C |x|^{-\bp-\bm}\hbox{ for all }y\in B_{R|x|}(0)\setminus \left(B_{\alpha|x|}(0)\cup B_{\alpha|x|}(x)\right)$$
Taking $\alpha>0$ small enough and $R>0$ large enough, we then get \eqref{est:5} for $|x|<\delta_3$. The general case for arbitrary $x\in \Omega\setminus \{0\}$ then follows from \eqref{est:2}. This prove \eqref{est:5}.  

Next we claim that for all $c_1,c_2>0$, there exists $C(c_1,c_2)>0$ such that
 \begin{equation}\label{est:6}
|x-y|^{n-2}G(x,y)\leq C(c_1,c_2)\hbox{ for }x,y\in\Omega\hbox{ such that }c_1|x|<|y|<c_2|x|.
\end{equation}
For that, we fix $x\in B_{\delta_0/2}(0)\setminus\{0\}$ and set
$$H(z):=|x|^{n-2}G_x(x+|x|z)\hbox{ for all }z\in B_{1/2}(0)\setminus \{0\}.$$
We have that $H\in C^2(\overline{B_{1/2}(0)}\setminus\{0\})$ and satisfies
$$-\Delta H-\left(\frac{\gamma}{\left|\frac{x}{|x|}+z\right|^2}+|x|^2h(x+|x|z)\right)H=\delta_0\hbox{ weakly in }B_{1/2}(0).$$
We now argue as in the proof of \eqref{est:2}. From \eqref{est:5}, we have that $|H(z)|\leq C$ for all $z\in \partial B_{1/2}(0)$ where $C$ is independent of $x\in B_{\delta_0/2}(0)\setminus\{0\}$. Let $\Gamma_0$ be the Green's function of $-\Delta -\left(\frac{\gamma}{\left|\frac{x}{|x|}+z\right|^2}+|x|^2h(x+|x|z)\right)$ at $0$ on $B_{1/2}(0)$ with Dirichlet boundary condition. Therefore, $H-\Gamma_0\in C^2(\overline{B_{1/2}(0)})$ and, via the comparison principle, it is bounded by its supremum on the boundary. Therefore $|z|^{n-2}H(z)\leq C$ for all $B_{1/2}(0)\setminus\{0\}$ where $C$ is independent of $x\in B_{\delta_0/2}(0)\setminus\{0\}$. Scaling back and using \eqref{est:5}, we get \eqref{est:6} for $x\in B_{\delta_0/2}(0)\setminus\{0\}$. The general case is a consequence of \eqref{est:2}. This ends the proof of \eqref{est:6}. 

We now show that there exists $C>0$ such that
 \begin{equation}\label{est:7}
|y|^{\bm}|x|^{\bp}G(x,y)\leq C\hbox{ for }x,y\in\Omega\hbox{ such that }|y|<\frac{1}{2}|x|.
\end{equation}
Indeed, the proof goes essentially as in \eqref{est:3}. Fix $x\in B_{\delta_0/2}(0)$, $x\neq 0$, and set $H(z):=|x|^{n-2}G_x(|x|z)$ for $z\in B_{1/2}(0)\setminus\{0\}$. We have that
$$-\Delta H-\left(\frac{\gamma}{|z|^2}+|x|^2h(|x|z)\right)H=0\hbox{ in }H_1^2(B_{1/2}(0)).$$
Moreover, it follows from \eqref{est:5} that there exists $C>0$ such that $|H(z)|\leq C$ for all $z\in \partial B_{1/2}(0)$. Then, as above, using a super-solution, we get that there exists $C>0$ such that $0<H(z)\leq C|z|^{-\bm}$ for all $z\in B_{1/2}(0)\setminus\{0\}$. Scaling back yields \eqref{est:7} when $x\in B_{\delta_0/2}(0)$. The general case follows from \eqref{est:3}. This proves \eqref{est:7}. 

Again, by symmetry, we conclude that there exists $C>0$ such that
 \begin{equation}\label{est:8}
|x|^{\bm}|y|^{\bp}G(x,y)\leq C\hbox{ for }x,y\in\Omega\hbox{ such that }|x|<\frac{1}{2}|y|.
\end{equation}
Finally, one easily checks that \eqref{est:G:up} is a direct consequence of \eqref{est:7}, \eqref{est:8} and \eqref{est:6}. When $f\in C^\infty_c(\Omega)$, identity \eqref{id:172} is a consequence of \eqref{id:regul}. The general case follows from density and the integral controls on $G$. The behavior \eqref{asymp:G} is a consequence of the classification of solutions to harmonic equations and Theorem \ref{th:sing}.

To conclude, we shall briefly sketch the proof of the lower bound  \eqref{est:G:low}. Indeed, in Steps 3 and 4, we repeatedly used the comparison principle to get the upper bound for $G$ by considering domains on the boundary of which $G$ was bounded from above. As one checks, in the case when $x,y$ are in $\omega\subset\subset\Omega$, $G$ is also bounded from below by some positive constant on the boundary of these domains. This yields the lower bound \eqref{est:G:low}, and completes 
 the proof of Theorem \ref{th:green:gamma:domain}.

\section*{\, Appendix B: Green's function for $-\Delta-\gamma|x|^{-2}$ on $\rn$}
In this section, we prove the following:
\begin{theorem}\label{th:green:gamma} Fix $\gamma<\frac{(n-2)^2}{4}$. For all $p\in\rnp$, there exists $G:\rn\setminus\{0,p\}\to\rr$ such that 

\smallskip\noindent{\bf (i)} $G\in H_{1,loc}^2(\rn\setminus\{p\})$,

\smallskip\noindent{\bf (ii)} For all $\varphi\in C^\infty_c(\rn)$, we have that
\begin{equation}\label{eq:kernel}
\varphi(p)=\int_{\rn}G(x)\left(-\Delta \varphi-\frac{\gamma }{|x|^2}\varphi\right)\, dx\hbox{ for all }\varphi\in C^\infty_c(\rn)
\end{equation}
\medskip\noindent Moreover, if $G,G'$ satisfy $(i)$ and $(ii)$ and are positive, then there exists $C\in\rr$ such that $G(x)-G'(x)=C|x|^{-\bm}$ for all $x\in\rn\setminus\{0,p\}$.


\medskip\noindent In addition, there exists one and only one function $G:=G_p>0$ such that (i) and (ii) hold and

\smallskip\noindent{\bf (iii)} For all $p\in\rnp$, there exists $c_0(p),c_\infty(p)>0$ such that
\begin{equation*}
G_p(x)\sim_{x\to 0} \frac{c_0(p)}{|x|^{\bm}}\hbox{ and }G_p(x)\sim_{x\to \infty} \frac{c_\infty(p)}{|x|^{\bp}}
\end{equation*}
and \begin{equation}
G_p(x)\sim_{x\to p}\frac{1}{(n-2)\omega_{n-1}|x-p|^{n-2}}.
\end{equation}
\smallskip\noindent{\bf (iv)}  There exists $c>0$ independent of $p$ such that
\begin{equation}\label{est:G:glob}
c^{-1} \left(\frac{\max\{|p|,|x|\}}{\min\{|p|,|x|\}}\right)^{\bm}|x-p|^{2-n}\leq G_p(x)\leq c \left(\frac{\max\{|p|,|x|\}}{\min\{|p|,|x|\}}\right)^{\bm}|x-p|^{2-n}
\end{equation}
\end{theorem}

\medskip\noindent {\it Remark:} Note that when $\gamma=0$, we have $\bm=0$, $\bp=n-2$ and $G(p,x)=\frac{1}{(n-2)\omega_{n-1}}|x-p|^{2-n}$ for all $x,p\in\rn$, $x\neq p$.

\medskip\noindent{\it Proof:} 
We shall again proceed with several steps. 

\smallskip\noindent{\bf Step 1: Construction of a positive kernel at a given point:} For a fixed $p_0\in\rn\setminus\{0\}$, we show that there exists $G\in C^2(\rn\setminus\{0,p_0\})$ such that
\begin{equation}\label{lim:G:bis}
\left\{\begin{array}{ll}
-\Delta G-\frac{\gamma}{|x|^2}G=0&\hbox{ in }\rn\setminus\{0,p_0\}\\
G>0&\\G\in L^{\frac{2n}{n-2}}(B_\delta(0))&\hbox{ with }\delta:=|p_0|/4\\
G\hbox{ satisfies }(ii).
\end{array}\right.
\end{equation}
Indeed, let $\tilde{\eta}\in C^\infty(\rr)$ be a nondecreasing function such that $0\leq \tilde{\eta}\leq 1$, $\tilde{\eta}(t)=0$ for all $t\leq 1$ and $\tilde{\eta}(t)=1$ for all $t\geq 2$. For $\eps>0$, set $\eta_\eps(x):=\tilde{\eta}\left(\frac{|x|}{\eps}\right)$ for all $x\in \rn$.
For $R>0$, we argue as in the proof of \eqref{bnd:coer} to deduce that the operator $-\Delta-\frac{\gamma\eta_\eps}{|x|^2}$ is coercive on $B_R(0)$ and that there exists $c>0$ independent of $R,\eps>0$ such that
\begin{equation*}
\int_{B_R(0)}\left(|\nabla \varphi|^2-\frac{\gamma\eta_\eps}{|x|^2}\varphi^2\right)\, dx\geq c\int_{B_R(0)}|\nabla \varphi|^2\, dx \quad \hbox{for all $\varphi\in C^\infty_c(B_R(0))$.}
\end{equation*}

\medskip\noindent Consider $R,\eps>0$ such that $R>2|p_0|$ and $\eps<\frac{|p_0|}{6}$, and let $G_{R,\eps}$ be the Green's function of $-\Delta-\frac{\gamma\eta_\eps}{|x|^2}$ in $B_R(0)$ at the point $p_0$ with Dirichlet boundary condition. We have that $G_{R,\eps}>0$ since the operator is coercive.

\medskip\noindent Fix $R_0>0$ and $q'\in (1,\frac{n}{n-2})$, then by arguing as in the proof of \eqref{int:bnd:G}, we get that there exists $C=C(\gamma,p_0, q', R_0)$ such that 
 \begin{equation}\label{bnd:G:1}
 \Vert G_{R,\eps}\Vert_{L^{q'}(B_{R_0}(0))}\leq C\hbox{ for all }R>R_0\hbox{ and }0<\eps<\frac{|p_0|}{6}, 
 \end{equation} 
 and
 \begin{equation}\label{bnd:G:0}
 \Vert G_{R,\eps}\Vert_{L^{\frac{2n}{n-2}}(B_{\delta_0}(0))}\leq C\hbox{ for all }R>R_0\hbox{ and }0<\eps<\frac{|p_0|}{6},
 \end{equation} 
where $\delta:=|p_0|/4$. Arguing again as in Step 2 of the proof of Theorem \ref{th:green:gamma:domain}, there exists $G\in C^2(\rn\setminus\{0,p_0\})$ such that
\begin{equation}\label{lim:G}
\left\{\begin{array}{ll}
G_{R,\eps}\to G\geq 0&\hbox{ in }C^2_{loc}(\rn\setminus\{0,p_0\})\hbox{ as }R\to +\infty,\; \eps\to 0\\
-\Delta G-\frac{\gamma}{|x|^2}G=0&\hbox{ in }\rn\setminus\{0,p_0\}\\
G\in L^{\frac{2n}{n-2}}(B_\delta(0))
\end{array}\right.
\end{equation}
Fix $\varphi\in C^\infty_c(\rn)$. For $R>0$ large enough, we have that $\varphi(p_0)=\int_{\rn}G_{R,\eps}(-\Delta\varphi-\gamma\eta_\eps|x|^{-2}\varphi)\, dx$. With the integral bounds above, we then get that $x\mapsto G(x) |x|^{-2}\in L^1_{loc}(\rn)$. Therefore, we get 
\begin{equation}\label{formula:G}
\varphi(p_0)=\int_{\rn}G(x)\left(-\Delta \varphi-\frac{\gamma }{|x|^2}\varphi\right)\, dx\hbox{ for all }\varphi\in C^\infty_c(\rn).
\end{equation}
As a consequence, $G>0$.

\medskip\noindent{\bf Step 2: Asymptotic behavior at $0$ and $p$ for solutions to \eqref{lim:G:bis}.}  It  follows from Theorem \ref{th:sing} below that either $G$ behaves like $|x|^{-\bm}$ or $|x|^{-\bp}$ at $0$. Since $G\in L^{\frac{2n}{n-2}}(B_\delta(0))$ for some small $\delta>0$ and $\bm<\frac{n-2}{2}<\bp$, we  get that there exists $c>0$ such that 
\begin{equation}\label{asymp:G:0}
\lim_{x\to 0}|x|^{\bm}G(x)=c.
\end{equation}
In addition, Theorem \ref{th:sing} yields $G\in H_{1,loc}^2(\rn\setminus\{p_0\})$. Since $G$ is positive and smooth in a neighborhood of $p$, it follows from \eqref{formula:G} and the classification of 
solutions to harmonic equations that
\begin{equation}\label{asymp:p}
G(x)\sim_{x\to p_0}\frac{1}{(n-2)\omega_{n-1}|x-p_0|^{n-2}}.
\end{equation}

\medskip\noindent{\bf Step 3: Asymptotic behavior at $\infty$ for solutions to \eqref{lim:G:bis}:} We let 
$$\tilde{G}(x):=\frac{1}{|x|^{n-2}}G\left(\frac{x}{|x|^2}\right)\hbox{ for all }x\in\rn\setminus\left\{0, \frac{p_0}{|p_0|^2}\right\},$$ 
be the Kelvin's transform of $G$. We have that
 $$-\Delta \tilde{G}-\frac{\gamma}{|x|^2}\tilde{G}=0\hbox{ in }\rn\setminus\left\{0, \frac{p_0}{|p_0|^2}\right\}.$$
Since $\tilde{G}>0$, it follows from Theorem \ref{th:sing} that there exists $c_1>0$ such that
\begin{equation*}
\hbox{either }\tilde{G}(x)\sim_{x\to 0}\frac{c_1}{|x|^{\bm}}\hbox{ or }\tilde{G}(x)\sim_{x\to 0}\frac{c_1}{|x|^{\bp}}.
\end{equation*}
Coming back to $G$, we get that
\begin{equation*}
\hbox{either }G(x)\sim_{x\to \infty}\frac{c_1}{|x|^{\bp}}\hbox{ or }G(x)\sim_{|x|\to \infty}\frac{c_1}{|x|^{\bm}}.
\end{equation*}
Assuming we are in the second case, for any $c\leq c_1$, we define
$$\bar{G}_c(x):=G(x)-\frac{c}{|x|^{\bm}}\hbox{ in }\rn\setminus\{0, p_0\},$$
which satisfy $-\Delta \bar{G}-\frac{\gamma}{|x|^2}\bar{G}=0$ in $\rn\setminus\{0, p_0\}$. It follows from \eqref{asymp:G:0} and \eqref{asymp:p} that for $c<c_1$, $\bar{G}_c>0$ around $p_0$ and $\infty$. It then follows from the coercivity of $-\Delta-\gamma|x|^{-2}$ that $\bar{G}_c>0$ in $\rn\setminus\{0, p\}$ for $c<c_1$. Letting $c\to c_1$ yields $\bar{G}_{c_1}\geq 0$, and then $\bar{G}_{c_1}> 0$ since it is positive around $p_0$.  Since $\bar{G}_{c_1}(x)=o(|x|^{-\bm})$ as $|x|\to \infty$, performing again a Kelvin transform and using Theorem \ref{th:sing}, we get that $|x|^{\bp}\bar{G}_{c_1}(x)\to c_2>0$ as $|x|\to \infty$. Then there exists $c_3>0$ such that

$$\lim_{x\to 0}|x|^{\bm}\bar{G}_{c_1}(x)=c_3>0\hbox{ and }\lim_{x\to \infty}|x|^{\bp}\bar{G}_{c_1}(x)=c_2.$$
Since $x\mapsto |x|^{-\bm}\in H_{1,loc}^2(\rn)$, we get that $\varphi(p)=\int_{\rn}\bar{G}_{c_1}(x)\left(-\Delta \varphi-\frac{\gamma }{|x|^2}\varphi\right)\, dx$ for all $\varphi\in C^\infty_c(\rn)$. 

\medskip\noindent{\bf Step 4: Uniqueness:} Let $G_1,G_2>0$ be 2 functions such that $(i),(ii)$ hold for $p:=p_0$, and set $H:=G_1-G_2$. It follows from Steps 2 and 3 that there exists $c\in\rr$ such that $H'(x):=H(x)-c|x|^{-\bm}$ satisfies
\begin{equation}\label{bnd:H:prime}
H'(x)=_{x\to 0}O\left(|x|^{-\bm}\right)\hbox{ and }H'(x)=_{|x|\to \infty}O\left(|x|^{-\bp}\right).
\end{equation}
We then have that $H\in H_{1,loc}^2(\rn\setminus\{p_0\})$ is such that
$$\int_{\rn}H'(x)\left(-\Delta\varphi-\frac{\gamma}{|x|^2}\varphi\right)\, dx=0 \quad \hbox{for all $\varphi\in C^\infty_c(\rn)$.}$$
The ellipticity of the Laplacian then yields that $H'\in C^\infty(\rnp)$. The pointwise bounds \eqref{bnd:H:prime} yield that $H'\in \dundeuxrn$. Multiplying $-\Delta H'-\frac{\gamma}{|x|^2}H'=0$ by $H'$, integrating by parts and using the coercivity yields that $H'\equiv 0$, and therefore, $G_1-G_2=c|\cdot|^{-\bm}$. This proves uniqueness.

\medskip\noindent{\bf Step 5: Existence.}  It follows from Step 3 that, up to substracting a multiple of $|\cdot|^{-\bm}$,  
there exists $G_{p_0}>0$ satisfying (i), (ii) and the pointwise controls (iii) at $p_0$. It is a consequence of (iii) that there exists $c>0$ such that
\begin{equation*}
c^{-1} \left(\frac{\max\{1,|x|\}}{\min\{1,|x|\}}\right)^{\bm}|x-p_0|^{2-n}\leq G_{p_0}(x)\leq c \left(\frac{\max\{1,|x|\}}{\min\{1,|x|\}}\right)^{\bm}|x-p_0|^{2-n}
\end{equation*}
for all $x\in \rn\setminus\{0,p_0\}$, $c$ depending on $p_0$. For $p\in\rnp$, consider $\rho_p: \rn\to\rn$ a linear isometry such that $\rho_p(\frac{p_0}{|p_0|})=\frac{p}{|p|}$, and define
$$G_p(x):=\left(\frac{|p_0|}{|p|}\right)^{n-2}G_{p_0}\left(\left(\rho_p^{-1}\left(\frac{|p_0|}{|p|}x\right)\right)\right)\hbox{ for all }x\in \rn\setminus\{0,p\}.$$
It is easy to check that $G_p>0$ and that it satisfies (i), (ii), (iii) and (iv).

%

\section*{\, Appendix C: Singular solutions to $-\Delta u-c(x) |x|^{-2}u=0$}
We collect here a few results that should be classical, but quite difficult to find in the literature. These results and their proofs are closely related to the work of the authors in \cite{gr4}, to which we shall frequently refer for details.

\begin{theorem}[Optimal regularity and Generalized Hopf's Lemma]\label{th:hopf} Fix $\gamma<\frac{(n-2)^2}{4}$ and let $f: \Omega\times \rr\to \rr$ be a Caratheodory function such that
$$|f(x, v)|\leq C|v| \left(1+\frac{|v|^{\crits-2}}{|x|^s}\right)\hbox{\rm  for all }x\in \Omega\hbox{ \rm and }v\in \rr.$$
Let $u\in H_1^2(B_1(0))$ be a weak solution of 
\begin{equation}\label{regul:eq}
-\Delta u-\frac{\gamma+O(|x|^\theta)}{|x|^2}u=f(x,u)\,\, \hbox{\rm  in }H_1^2(B_{1/2}(0))
\end{equation}
for some $\theta>0$. Then, there exists $K\in\rr$ such that
\begin{equation}\label{eq:hopf}
\lim_{x\to 0}\frac{u(x)}{|x|^{-\bm}}=K.
\end{equation}
Moreover, if $u\geq 0$ and $u\not\equiv 0$, we have that $K>0$.
\end{theorem}

\begin{theorem}\label{th:sing} Let $u\in C^2(B_1(0)\setminus \{0\})$ be a positive solution to 
$$-\Delta u-\frac{c(x)}{|x|^2}u=0\hbox{ in }B_1(0)\setminus\{0\}$$
where $c(x)=\gamma+O(|x|^\theta)$ as $x\to 0$ with $\gamma<(n-2)^2/4$ and $\theta\in (0,1)$. Then there exists $\alpha>0$ such that
\begin{equation*}
\hbox{either }u(x)\sim_{x\to 0}\frac{\alpha}{|x|^{\bm}}\hbox{ or }u(x)\sim_{x\to 0}\frac{\alpha}{|x|^{\bp}}.
\end{equation*}
In particular, $u\in H_1^2(B_{1/2}(0))$ if and only if the first case holds.
\end{theorem}

\begin{proposition}\label{prop:funda:2} Let $u\in C^2(\rn\setminus \{0\})$ be a nonnegative function such that
\begin{equation}\label{eq:u:36}
-\Delta u-\frac{\gamma}{|x|^2}u=0\hbox{ in }\rn\setminus \{0\}.
\end{equation}
Then there exist $\lambda_-,\lambda_+\geq 0$ such that
$$u(x)=\lambda_- |x|^{-\bm}+\lambda_+ |x|^{-\bp}\hbox{ for all }x\in \rn\setminus \{0\}.$$
\end{proposition}

\medskip\noindent{\it Proofs:} The proofs of these results follow closely the proofs of Theorems 6.1 and 7.1 and Proposition 7.4 of \cite{gr4}. Here are the ingredients to adapt:

\smallskip\noindent{\it Sub- and super-solutions:} The first step is the following result:
\begin{proposition}\label{prop:sub:super} Fix $\gamma<(n-2)^2/4$ and $\theta\in (0,1)$ and let $c:B_1(0)\to\rr$ be such that $c(x)=\gamma+O(|x|^\theta)$ as $x\to 0$. We choose $\beta\in \{\bm,\bp\}$. Then there exists $u_\beta^{(+)}, u_\beta^{(-)}\in C^2(B_1(0))$ such that for $\delta>0$ small enough, 
\begin{equation}\label{ppty:super:2}\left\{\begin{array}{ll}
-\Delta u_{\beta}^{(+)}-\frac{c(x)}{|x|^2}u_{\beta}^{(+)}>0&\hbox{ in }B_\delta(0)\\
u_{\beta}^{(+)}(x)\sim|x|^{-\beta}& \hbox{ as }x\to 0
\end{array}\right.\hbox{ ; }\left\{\begin{array}{ll}
-\Delta u_{\beta}^{(-)}-\frac{c(x)}{|x|^2}u_{\beta}^{(-)}<0&\hbox{ in }B_\delta(0)\\
u_{\beta}^{(-)}(x)\sim|x|^{-\beta}&\hbox{ as }x\to 0
\end{array}
\right.\end{equation}
\end{proposition}
The proof is as follows. For $\beta\in \{\bm,\bp\}$, we define $u_{\beta}:\, x\mapsto |x|^{-\beta}+\lambda |x|^{-\beta'}$. A straightforward computation yields
$$-\Delta u_{\beta}-\frac{c(x)}{|x|^2}u_{\beta}=|x|^{-\beta'-2}\left(\lambda (\beta'(n-2-\beta')-\gamma)+O(|x|^\theta)+O(|x|^{\theta-(\beta-\beta')}\right)$$
as $x\to 0$. Then, choosing $\beta'\in \rr$ such that $0<\beta-\beta'<\theta$ and $\beta'(n-2-\beta')-\gamma\neq 0$, we get either a sub- or a supersolution taking $\lambda$ positive or negative. This proves the proposition.

\smallskip\noindent{\it Sub-solution with Dirichlet boundary condition:} We let $u_{\bp}$ as above be a super-solution on $B_\delta(0)\setminus\{0\}$. Take $\eta\in C^\infty(\rn)$ such that $\eta(x)=0$ for $x\in B_{\delta/4}(0)$ and $\eta(x)=1$ for $x\in \rn\setminus B_{\delta/3}(0)$. Define on $B_\delta(0)$ the function 
$$f(x):=\left(-\Delta-\frac{c(x)}{|x|^2}\right)(\eta u_{\bp}), $$
Note that $f$ vanishes around $0$ and that it is in $C^\infty(\overline{B_\delta(0)})$. Let $v\in D^{1,2}(B_\delta(0))$ be such that
$$\left\{\begin{array}{ll}
-\Delta v-\frac{c(x)}{|x|^2}v=f&\hbox{ in }B_\delta(0)\\
\hfill v=0&\hbox{ on }\partial B_\delta(0).
\end{array}\right.$$
Note that for $\delta>0$ small enough, $-\Delta-(\gamma+O(|x|^\theta))|x|^{-2}$ is coercive on $B_\delta(0)$, and therefore, the existence of $v$ is ensured for small $\delta$. Define
$$u_{\bp}^{(d)}:=u_{\bp}-\eta u_{\bp}+v.$$
The definition of $\eta$ and $v$ yields
\begin{equation}\label{ppty:super:3}\left\{\begin{array}{ll}
-\Delta u_{\bp}^{(d)}-\frac{c(x)}{|x|^2}u_{\bp}^{(d)}>0&\hbox{ in }B_\delta(0)-\{0\}\\
\hfill u_{\bp}^{(d)}=0 &\hbox{ in }\partial B_\delta(0)
\end{array}
\right.\end{equation}
Moreover, since $-\Delta v-c(x)|x|^{-2}v=0$ around $0$ and $v\in D^{1,2}(B_\delta(0))$, it follows from Theorem \ref{th:hopf} that there exists $C>0$ such that $|v(x)|\leq C |x|^{-\bm}$ for all $x\in B_\delta(0)$. Then it follows from the expression of $u_{\bp}$ that
$$u_{\bp}(x)\sim_{x\to 0}|x|^{-\bp}.$$
We then get a supersolution satisfying \eqref{ppty:super:3} with the above behavior at $0$. This is similar for a subsolution.

\smallskip\noindent In the proof above, it is important that the operator $-\Delta-c(x)|x|^{-2}$ is coercive on $B_\delta(0)$ for $\delta>0$ small enough. Now let $\Omega$ be a smooth bounded domain of $\rn$ and let $\gamma<(n-2)^2/4$ and $h\in C^{0,\theta(\overline{\Omega})}$ be such that  $-\Delta-(\gamma|x|^{-2}+h)$ is coercive on $\Omega$. Arguing as above, we get that there exists $u_{\bp}^{(d,\Omega)}\in C^{2,\theta}(\overline{\Omega}\setminus\{0\})$ such that
\begin{equation}\label{ppty:super}\left\{\begin{array}{ll}
-\Delta u_{\bp}^{(d,\Omega)}-\left(\frac{\gamma}{|x|^2}+h\right)u_{\bp}^{(d)}>0&\hbox{ in }\Omega\setminus\{0\}\\
u_{\bp}^{(d,\Omega)}>0 &\hbox{ in }\Omega\setminus\{0\}\\
u_{\bp}^{(d,\Omega)}=0 &\hbox{ in }\partial \Omega
\end{array}
\right.\end{equation}
and
$$u_{\bp}^{(d,\Omega)}(x)\sim_{x\to 0}|x|^{-\bp}.$$
Similarly, we get a subsolution.

\smallskip\noindent These points are enough to adapt the proofs of the above-mentioned results of \cite{gr4} to our context.\qed

\end{document}